\documentclass[preprint,3p,12pt]{elsarticle}

\usepackage{microtype}
\usepackage[font=footnotesize,labelfont=bf,labelformat=simple]{subcaption}
\usepackage{fancyhdr}
\usepackage{amsmath}
\usepackage{amssymb}
\usepackage{amsthm}
\usepackage{hyperref}
\usepackage{graphicx}
\usepackage{import}
\usepackage{xcolor}
\usepackage{multirow}
\usepackage{multicol}
\usepackage{csquotes}
\usepackage{newtxtext,newtxmath}
\usepackage{textcomp}
\usepackage{psfrag}
\usepackage{pgfplots}
\usepackage{epsfig}
\usepackage{pst-all}
\usepackage{eqnarray}
\usepackage{stmaryrd}
\usepackage{tcolorbox}
\usepackage{setspace}
\usepackage{xfrac}
\usepackage{relsize}
\usepackage{etoolbox}
\usepackage{hyperref}
\usepackage{longtable}
\usepackage{fancyvrb}
\usepackage{booktabs}
\usepackage{tikz}
\usepackage{tikz-3dplot}
\usepackage{bm}

\usetikzlibrary{positioning,shapes,arrows,backgrounds,calc,fit,perspective,spy,patterns,external,decorations,decorations.markings}

\psset{xunit=1cm, yunit=1cm}

\pgfplotsset{compat=newest} 
\pgfplotsset{plot coordinates/math parser=false}

\DeclareMathOperator*{\A}{\mathlarger{\mathlarger{{\mathsf{A}}}}}

\biboptions{numbers,sort&compress}

\definecolor{blue}{RGB}{31, 119, 180}
\definecolor{lightblue}{RGB}{214, 199, 232}

\definecolor{orange}{RGB}{255, 127, 14}
\definecolor{lightorange}{RGB}{255, 187, 120}

\definecolor{green}{RGB}{44, 160, 44}
\definecolor{lightgreen}{RGB}{152, 223, 138}

\definecolor{red}{RGB}{214, 39, 40}
\definecolor{lightred}{RGB}{255, 152, 150}

\definecolor{purple}{RGB}{148, 103, 189}
\definecolor{lightpurple}{RGB}{197, 176, 213}

\definecolor{brown}{RGB}{140, 86, 75}
\definecolor{lightbrown}{RGB}{196, 156, 148}

\definecolor{pink}{RGB}{227, 119, 194}
\definecolor{lightpink}{RGB}{247, 182, 210}

\definecolor{gray}{RGB}{127, 127, 127}
\definecolor{lightgray}{RGB}{199, 199, 199}

\definecolor{yellow}{RGB}{189, 189, 34}
\definecolor{lightyellow}{RGB}{219, 219, 141}

\definecolor{cyan}{RGB}{186, 190, 207}
\definecolor{lightcyan}{RGB}{158, 218, 229}

\definecolor{myblueD}{RGB}{25 102 179}
\definecolor{mygreenD}{RGB}{26 153 26}
\definecolor{myorangeD}{RGB}{190 100 0}

\definecolor{cFirst}{HTML}{0000E6}
\definecolor{cSecond}{HTML}{E60000}
\definecolor{cThird}{HTML}{2C8F00}
\definecolor{cFourth}{HTML}{FF751A}

\definecolor{cTab20c1}{HTML}{3182bd}
\definecolor{cTab20c2}{HTML}{6baed6}
\definecolor{cTab20c3}{HTML}{9ecae1}
\definecolor{cTab20c4}{HTML}{c6dbef}
\definecolor{cTab20c5}{HTML}{e6550d}
\definecolor{cTab20c6}{HTML}{fd8d3c}
\definecolor{cTab20c7}{HTML}{fdae6b}
\definecolor{cTab20c8}{HTML}{fdd0a2}
\definecolor{cTab20c9}{HTML}{31a354}
\definecolor{cTab20c10}{HTML}{74c476}
\definecolor{cTab20c11}{HTML}{a1d99b}
\definecolor{cTab20c12}{HTML}{c7e9c0}
\definecolor{cTab20c13}{HTML}{756bb1}
\definecolor{cTab20c14}{HTML}{9e9ac8}
\definecolor{cTab20c15}{HTML}{bcbddc}
\definecolor{cTab20c16}{HTML}{dadaeb}
\definecolor{cTab20c17}{HTML}{636363}
\definecolor{cTab20c18}{HTML}{969696}
\definecolor{cTab20c19}{HTML}{bdbdbd}
\definecolor{cTab20c20}{HTML}{d9d9d9}

\definecolor{cTab1}{HTML}{1f77b4}
\definecolor{cTab2}{HTML}{ff7f0e}
\definecolor{cTab3}{HTML}{2ca02c}
\definecolor{cTab4}{HTML}{d62728}
\definecolor{cTab5}{HTML}{9467bd}
\definecolor{cTab6}{HTML}{8c564b}
\definecolor{cTab7}{HTML}{e377c2}
\definecolor{cTab8}{HTML}{7f7f7f}
\definecolor{cTab9}{HTML}{bcbd22}
\definecolor{cTab10}{HTML}{17becf}

\colorlet{bulk}{black!20!white}

\DeclareFontFamily{U}{MnSymbolC}{}
\DeclareSymbolFont{MnSyC}{U}{MnSymbolC}{m}{n}
\DeclareFontShape{U}{MnSymbolC}{m}{n}{
    <-6> MnSymbolC5
   <6-7> MnSymbolC6
   <7-8> MnSymbolC7
   <8-9> MnSymbolC8
   <9-10> MnSymbolC9
  <10-12> MnSymbolC10
  <12->   MnSymbolC12}{}
\DeclareMathSymbol{\filledsquare}{\mathbin}{MnSyC}{104}

\pgfplotscreateplotcyclelist{color_list}{
	blue \\
	orange \\
	green \\
	red \\
	purple \\
	brown \\
	pink \\
	gray \\
	yellow \\
	cyan \\
	lightblue \\
	lightorange \\
	lightgreen \\
	lightred \\
	lightpurple \\
	lightbrown \\
	lightpink \\
	lightgray \\
	lightyellow \\
	lightcyan \\
}

\pgfplotscreateplotcyclelist{color_list_mark}{
	blue,        mark = square*, mark size = 1.5pt, mark options = {fill = blue!50!white, solid} \\
	orange,      mark = *, mark size = 1.5pt, mark options = {fill = orange!50!white, solid} \\
	green,       mark = triangle*, mark size = 1.5pt, mark options = {fill = green!50!white, solid} \\
	red,         mark = diamond*, mark size = 1.5pt, mark options = {fill = red!50!white, solid} \\
	purple,      mark = pentagon*, mark size = 1.5pt, mark options = {fill = purple!50!white, solid} \\
	brown,       mark = square*, mark size = 1.5pt, mark options = {fill = brown!50!white, solid} \\
	pink,        mark = *, mark size = 1.5pt, mark options = {fill = pink!50!white, solid} \\
	gray,        mark = triangle*, mark size = 1.5pt, mark options = {fill = gray!50!white, solid} \\
	yellow,      mark = diamond*, mark size = 1.5pt, mark options = {fill = yellow!50!white, solid} \\
	cyan,        mark = pentagon*, mark size = 1.5pt, mark options = {fill = cyan!50!white, solid} \\
	lightblue,   mark = square*, mark size = 1.5pt, mark options = {fill = lightblue!50!white, solid} \\
	lightorange, mark = *, mark size = 1.5pt, mark options = {fill = lightorange!50!white, solid} \\
	lightgreen,  mark = triangle*, mark size = 1.5pt, mark options = {fill = lightgreen!50!white, solid} \\
	lightred,    mark = diamond*, mark size = 1.5pt, mark options = {fill = lightred!50!white, solid} \\
	lightpurple, mark = pentagon*, mark size = 1.5pt, mark options = {fill = lightpurple!50!white, solid} \\
	lightbrown,  mark = square*, mark size = 1.5pt, mark options = {fill = lightbrown!50!white, solid} \\
	lightpink,   mark = *, mark size = 1.5pt, mark options = {fill = lightpink!50!white, solid} \\
	lightgray,   mark = triangle*, mark size = 1.5pt, mark options = {fill = lightgray!50!white, solid} \\
	lightyellow, mark = diamond*, mark size = 1.5pt, mark options = {fill = lightyellow!50!white, solid} \\
	lightcyan,   mark = pentagon*, mark size = 1.5pt, mark options = {fill = lightcyan!50!white, solid} \\
}

\pgfplotsset{
	colormap={TemperatureMap}{%
		rgb = (0.14901961, 0.28235295, 0.92941177)
		rgb = (0.1764706, 0.30588236, 0.93333334)
		rgb = (0.18039216, 0.30980393, 0.93333334)
		rgb = (0.19215687, 0.32156864, 0.9372549)
		rgb = (0.19607843, 0.3254902, 0.9372549)
		rgb = (0.20784314, 0.3372549, 0.9372549)
		rgb = (0.20784314, 0.3372549, 0.9372549)
		rgb = (0.21176471, 0.34117648, 0.9372549)
		rgb = (0.22352941, 0.3529412, 0.9372549)
		rgb = (0.22745098, 0.35686275, 0.9372549)
		rgb = (0.23921569, 0.36862746, 0.9372549)
		rgb = (0.23921569, 0.36862746, 0.9372549)
		rgb = (0.2509804, 0.38039216, 0.9411765)
		rgb = (0.25490198, 0.3882353, 0.9411765)
		rgb = (0.25490198, 0.3882353, 0.9411765)
		rgb = (0.26666668, 0.4, 0.9411765)
		rgb = (0.27058825, 0.40392157, 0.9411765)
		rgb = (0.28235295, 0.41568628, 0.9411765)
		rgb = (0.28627452, 0.41960785, 0.9411765)
		rgb = (0.28627452, 0.41960785, 0.9411765)
		rgb = (0.3019608, 0.43529412, 0.9411765)
		rgb = (0.3019608, 0.43529412, 0.9411765)
		rgb = (0.31764707, 0.4509804, 0.94509804)
		rgb = (0.31764707, 0.4509804, 0.94509804)
		rgb = (0.31764707, 0.4509804, 0.94509804)
		rgb = (0.3372549, 0.46666667, 0.94509804)
		rgb = (0.3372549, 0.46666667, 0.94509804)
		rgb = (0.3529412, 0.48235294, 0.94509804)
		rgb = (0.3529412, 0.48235294, 0.94509804)
		rgb = (0.35686275, 0.4862745, 0.94509804)
		rgb = (0.37254903, 0.49803922, 0.94509804)
		rgb = (0.37254903, 0.49803922, 0.94509804)
		rgb = (0.3764706, 0.5019608, 0.94509804)
		rgb = (0.3882353, 0.5137255, 0.9490196)
		rgb = (0.39607844, 0.5176471, 0.9490196)
		rgb = (0.40784314, 0.5254902, 0.9490196)
		rgb = (0.40784314, 0.5254902, 0.9490196)
		rgb = (0.41960785, 0.5372549, 0.9490196)
		rgb = (0.42352942, 0.5411765, 0.9490196)
		rgb = (0.43529412, 0.5529412, 0.9529412)
		rgb = (0.44313726, 0.5568628, 0.9529412)
		rgb = (0.44313726, 0.5568628, 0.9529412)
		rgb = (0.45490196, 0.5686275, 0.9529412)
		rgb = (0.45882353, 0.57254905, 0.9529412)
		rgb = (0.47843137, 0.5882353, 0.9529412)
		rgb = (0.48235294, 0.5921569, 0.9529412)
		rgb = (0.48235294, 0.5921569, 0.9529412)
		rgb = (0.5019608, 0.6117647, 0.95686275)
		rgb = (0.5058824, 0.6117647, 0.95686275)
		rgb = (0.5058824, 0.6117647, 0.95686275)
		rgb = (0.5294118, 0.627451, 0.95686275)
		rgb = (0.5294118, 0.627451, 0.95686275)
		rgb = (0.5529412, 0.64705884, 0.9607843)
		rgb = (0.5529412, 0.64705884, 0.9607843)
		rgb = (0.5568628, 0.6509804, 0.9607843)
		rgb = (0.5764706, 0.6666667, 0.9607843)
		rgb = (0.5803922, 0.67058825, 0.9607843)
		rgb = (0.59607846, 0.6862745, 0.9647059)
		rgb = (0.59607846, 0.6862745, 0.9647059)
		rgb = (0.6039216, 0.6901961, 0.9647059)
		rgb = (0.61960787, 0.7019608, 0.9647059)
		rgb = (0.627451, 0.70980394, 0.9647059)
		rgb = (0.6431373, 0.72156864, 0.96862745)
		rgb = (0.6431373, 0.72156864, 0.96862745)
		rgb = (0.654902, 0.73333335, 0.96862745)
		rgb = (0.6627451, 0.7372549, 0.96862745)
		rgb = (0.6627451, 0.7372549, 0.96862745)
		rgb = (0.6784314, 0.7490196, 0.972549)
		rgb = (0.6862745, 0.7529412, 0.972549)
		rgb = (0.7019608, 0.76862746, 0.972549)
		rgb = (0.7058824, 0.77254903, 0.972549)
		rgb = (0.7058824, 0.77254903, 0.972549)
		rgb = (0.7254902, 0.78431374, 0.9764706)
		rgb = (0.7294118, 0.7882353, 0.9764706)
		rgb = (0.7490196, 0.8039216, 0.9764706)
		rgb = (0.7490196, 0.8039216, 0.9764706)
		rgb = (0.7490196, 0.8039216, 0.9764706)
		rgb = (0.77254903, 0.8235294, 0.98039216)
		rgb = (0.77254903, 0.8235294, 0.98039216)
		rgb = (0.7921569, 0.8392157, 0.98039216)
		rgb = (0.7921569, 0.8392157, 0.98039216)
		rgb = (0.79607844, 0.84313726, 0.98039216)
		rgb = (0.8156863, 0.85490197, 0.9843137)
		rgb = (0.8156863, 0.85490197, 0.9843137)
		rgb = (0.81960785, 0.85882354, 0.9843137)
		rgb = (0.83137256, 0.8666667, 0.9843137)
		rgb = (0.8352941, 0.87058824, 0.9843137)
		rgb = (0.84313726, 0.8784314, 0.9882353)
		rgb = (0.84313726, 0.8784314, 0.9882353)
		rgb = (0.85490197, 0.8862745, 0.9882353)
		rgb = (0.85882354, 0.8901961, 0.9882353)
		rgb = (0.8666667, 0.8980392, 0.9882353)
		rgb = (0.87058824, 0.9019608, 0.9882353)
		rgb = (0.87058824, 0.9019608, 0.9882353)
		rgb = (0.88235295, 0.9098039, 0.99215686)
		rgb = (0.8862745, 0.9137255, 0.99215686)
		rgb = (0.8980392, 0.92156863, 0.99215686)
		rgb = (0.9019608, 0.9254902, 0.99215686)
		rgb = (0.9019608, 0.9254902, 0.99215686)
		rgb = (0.9137255, 0.93333334, 0.99215686)
		rgb = (0.9137255, 0.93333334, 0.99215686)
		rgb = (0.9137255, 0.93333334, 0.99215686)
		rgb = (0.92941177, 0.94509804, 0.99607843)
		rgb = (0.92941177, 0.94509804, 0.99607843)
		rgb = (0.9411765, 0.9529412, 0.99215686)
		rgb = (0.9411765, 0.9529412, 0.99215686)
		rgb = (0.9411765, 0.9529412, 0.99215686)
		rgb = (0.94509804, 0.9607843, 0.9843137)
		rgb = (0.94509804, 0.9607843, 0.98039216)
		rgb = (0.9529412, 0.9647059, 0.972549)
		rgb = (0.9529412, 0.9647059, 0.972549)
		rgb = (0.9529412, 0.9647059, 0.96862745)
		rgb = (0.95686275, 0.96862745, 0.9607843)
		rgb = (0.95686275, 0.96862745, 0.95686275)
		rgb = (0.9607843, 0.972549, 0.9529412)
		rgb = (0.9607843, 0.972549, 0.9529412)
		rgb = (0.9647059, 0.9764706, 0.94509804)
		rgb = (0.96862745, 0.9764706, 0.9411765)
		rgb = (0.96862745, 0.9764706, 0.9411765)
		rgb = (0.972549, 0.98039216, 0.9372549)
		rgb = (0.972549, 0.98039216, 0.93333334)
		rgb = (0.98039216, 0.9843137, 0.9254902)
		rgb = (0.98039216, 0.9843137, 0.92156863)
		rgb = (0.98039216, 0.9843137, 0.92156863)
		rgb = (0.9843137, 0.9882353, 0.9137255)
		rgb = (0.9843137, 0.9882353, 0.9137255)
		rgb = (0.9882353, 0.9882353, 0.8901961)
		rgb = (0.9882353, 0.9882353, 0.8901961)
		rgb = (0.9882353, 0.9882353, 0.8901961)
		rgb = (0.9882353, 0.99215686, 0.87058824)
		rgb = (0.9882353, 0.99215686, 0.87058824)
		rgb = (0.9882353, 0.99215686, 0.8666667)
		rgb = (0.9882353, 0.99215686, 0.84705883)
		rgb = (0.9882353, 0.99215686, 0.84313726)
		rgb = (0.99215686, 0.99215686, 0.8235294)
		rgb = (0.99215686, 0.99215686, 0.8235294)
		rgb = (0.99215686, 0.99215686, 0.81960785)
		rgb = (0.99215686, 0.99215686, 0.8039216)
		rgb = (0.99215686, 0.99215686, 0.79607844)
		rgb = (0.99215686, 0.99215686, 0.78039217)
		rgb = (0.99215686, 0.99215686, 0.78039217)
		rgb = (0.99215686, 0.99215686, 0.77254903)
		rgb = (0.99607843, 0.99215686, 0.75686276)
		rgb = (0.99607843, 0.99215686, 0.74509805)
		rgb = (0.99607843, 0.99215686, 0.7372549)
		rgb = (0.99607843, 0.99215686, 0.7411765)
		rgb = (0.99607843, 0.99215686, 0.7137255)
		rgb = (0.99607843, 0.99215686, 0.7058824)
		rgb = (0.99607843, 0.99215686, 0.70980394)
		rgb = (0.99607843, 0.99215686, 0.6745098)
		rgb = (0.99607843, 0.99215686, 0.6666667)
		rgb = (0.99607843, 0.99215686, 0.63529414)
		rgb = (0.99607843, 0.99215686, 0.6313726)
		rgb = (0.99607843, 0.99215686, 0.6313726)
		rgb = (0.99607843, 0.99215686, 0.59607846)
		rgb = (0.99607843, 0.99215686, 0.59607846)
		rgb = (0.99607843, 0.99215686, 0.5568628)
		rgb = (0.99607843, 0.99215686, 0.5568628)
		rgb = (0.99607843, 0.99215686, 0.5529412)
		rgb = (0.99607843, 0.9882353, 0.52156866)
		rgb = (0.99607843, 0.9882353, 0.5137255)
		rgb = (0.99607843, 0.9882353, 0.48235294)
		rgb = (0.99607843, 0.9882353, 0.4862745)
		rgb = (0.99607843, 0.9882353, 0.47843137)
		rgb = (0.99607843, 0.9882353, 0.44705883)
		rgb = (0.99607843, 0.9882353, 0.4509804)
		rgb = (0.99607843, 0.9882353, 0.4392157)
		rgb = (0.99215686, 0.9843137, 0.41960785)
		rgb = (0.9882353, 0.972549, 0.40784314)
		rgb = (0.9882353, 0.96862745, 0.40392157)
		rgb = (0.9882353, 0.96862745, 0.40392157)
		rgb = (0.9843137, 0.95686275, 0.3882353)
		rgb = (0.9843137, 0.9529412, 0.38431373)
		rgb = (0.98039216, 0.9411765, 0.37254903)
		rgb = (0.98039216, 0.9372549, 0.36862746)
		rgb = (0.98039216, 0.9372549, 0.36862746)
		rgb = (0.9764706, 0.9254902, 0.3529412)
		rgb = (0.9764706, 0.9254902, 0.34901962)
		rgb = (0.972549, 0.9098039, 0.33333334)
		rgb = (0.972549, 0.9098039, 0.33333334)
		rgb = (0.972549, 0.9098039, 0.33333334)
		rgb = (0.9647059, 0.89411765, 0.31764707)
		rgb = (0.9647059, 0.89411765, 0.31764707)
		rgb = (0.9647059, 0.89411765, 0.31764707)
		rgb = (0.9607843, 0.8784314, 0.29803923)
		rgb = (0.9607843, 0.8745098, 0.29411766)
		rgb = (0.95686275, 0.8627451, 0.28235295)
		rgb = (0.95686275, 0.8627451, 0.28235295)
		rgb = (0.95686275, 0.85882354, 0.28235295)
		rgb = (0.9529412, 0.8392157, 0.2784314)
		rgb = (0.9490196, 0.83137256, 0.27450982)
		rgb = (0.94509804, 0.8156863, 0.27058825)
		rgb = (0.94509804, 0.8156863, 0.27058825)
		rgb = (0.9411765, 0.80784315, 0.27058825)
		rgb = (0.9372549, 0.7921569, 0.26666668)
		rgb = (0.93333334, 0.77254903, 0.2627451)
		rgb = (0.93333334, 0.7647059, 0.2627451)
		rgb = (0.93333334, 0.7647059, 0.2627451)
		rgb = (0.92941177, 0.7490196, 0.25882354)
		rgb = (0.9254902, 0.7411765, 0.25882354)
		rgb = (0.9254902, 0.7411765, 0.25882354)
		rgb = (0.92156863, 0.72156864, 0.2509804)
		rgb = (0.92156863, 0.7176471, 0.2509804)
		rgb = (0.9137255, 0.69411767, 0.24705882)
		rgb = (0.9137255, 0.6901961, 0.24705882)
		rgb = (0.9137255, 0.6901961, 0.24705882)
		rgb = (0.9098039, 0.6666667, 0.24313726)
		rgb = (0.9098039, 0.6666667, 0.24313726)
		rgb = (0.9019608, 0.6431373, 0.23921569)
		rgb = (0.9019608, 0.6431373, 0.23921569)
		rgb = (0.9019608, 0.6392157, 0.23921569)
		rgb = (0.8980392, 0.6156863, 0.23137255)
		rgb = (0.8980392, 0.6117647, 0.23137255)
		rgb = (0.8901961, 0.5921569, 0.22745098)
		rgb = (0.8901961, 0.5921569, 0.22745098)
		rgb = (0.8901961, 0.58431375, 0.22745098)
		rgb = (0.8862745, 0.5686275, 0.22352941)
		rgb = (0.8862745, 0.5686275, 0.22352941)
		rgb = (0.88235295, 0.56078434, 0.22352941)
		rgb = (0.8784314, 0.5411765, 0.21960784)
		rgb = (0.8745098, 0.5254902, 0.21568628)
		rgb = (0.8745098, 0.5176471, 0.21568628)
		rgb = (0.8745098, 0.5176471, 0.21568628)
		rgb = (0.87058824, 0.49803922, 0.21176471)
		rgb = (0.87058824, 0.49019608, 0.21176471)
		rgb = (0.8627451, 0.47058824, 0.20392157)
		rgb = (0.8627451, 0.46666667, 0.20392157)
		rgb = (0.8627451, 0.47058824, 0.20392157)
		rgb = (0.85882354, 0.4392157, 0.2)
		rgb = (0.85882354, 0.43529412, 0.2)
		rgb = (0.85490197, 0.40392157, 0.19607843)
		rgb = (0.85490197, 0.4, 0.19607843)
		rgb = (0.85490197, 0.4, 0.19607843)
		rgb = (0.84705883, 0.36078432, 0.19215687)
		rgb = (0.84705883, 0.36078432, 0.19215687)
		rgb = (0.84705883, 0.35686275, 0.19215687)
		rgb = (0.84313726, 0.32156864, 0.1882353)
		rgb = (0.84313726, 0.31764707, 0.1882353)
		rgb = (0.8392157, 0.28235295, 0.18039216)
		rgb = (0.8392157, 0.28627452, 0.18039216)
		rgb = (0.8392157, 0.2784314, 0.18039216)
		rgb = (0.8352941, 0.24705882, 0.1764706)
		rgb = (0.83137256, 0.23529412, 0.1764706)
		rgb = (0.827451, 0.20392157, 0.17254902)
		rgb = (0.827451, 0.20784314, 0.17254902)
		rgb = (0.827451, 0.19607843, 0.17254902)
		rgb = (0.8235294, 0.17254902, 0.16862746)
		rgb = (0.81960785, 0.14509805, 0.16470589)
		rgb = (0.81960785, 0.13333334, 0.16470589)
		rgb = (0.80784315, 0.07450981, 0.10980392)
	},
}

\usetikzlibrary{arrows.meta}
\pgfplotsset{compat=1.18}
\definecolor{oiGreen}{HTML}{009E73}
\definecolor{oiBlue}{HTML}{0072B2}
\definecolor{oiSky}{HTML}{56B4E9}
\definecolor{oiPurple}{HTML}{CC79A7}
\definecolor{oiOrange}{HTML}{E69F00}
\definecolor{oiVermillion}{HTML}{D55E00}
\definecolor{oiGrey}{HTML}{4D4D4D}
\definecolor{oiYellow}{HTML}{F0E442}
\definecolor{oiLightGrey}{HTML}{BBBBBB}
\pgfplotsset{
  every axis/.append style={
    axis line style={oiGrey, line width=0.5pt},
    tick style={oiGrey, line width=0.5pt},
    label style={font=\small},
    tick label style={font=\footnotesize},
    grid style={black!12, line width=0.35pt},
  },
}

\journal{Finite Elements in Analysis \& Design}

\begin{document}

\begin{frontmatter}
	\title{Efficient Thermo-Viscoplastic Analysis Using a Multi-Level $hp-$Finite Cell Method with Non-Negative Moment Fitting}

	\author[hsd]{Jan Niklas Schmäke\corref{cor1}}
	\author[cats]{Oliver Wege}
	\author[btu,hsd]{Martin Ruess\corref{cor2}}

	\address[hsd]{Düsseldorf University of Applied Sciences, Münsterstr. 156, 40476 Düsseldorf, Germany}
	\address[cats]{Chair for Computational Analysis of Technical Systems (CATS), RWTH Aachen University, Schinkelstraße 2, 52056, Aachen, Germany}
	\address[btu]{Computer-Aided Methods in Civil Engineering, Brandenburg University of Technology Cottbus-Senftenberg, Konrad-Wachsmann-Allee 2, 03046 Cottbus }

	\cortext[cor1]{Corresponding author:\\
	Jan Niklas Schmäke, University of Applied Sciences, Münsterstr. 156, 40476 Düsseldorf, Germany; e-mail: jan.schmaeke@hs-duesseldorf.de}
	\cortext[cor2]{Corresponding author:\\
	Martin Ruess, Computer-Aided Methods in Civil Engineering, Brandenburg University of Technology Cottbus-Senftenberg, Konrad-Wachsmann-Allee 2, 03046 Cottbus, Germany; e-mail: martin.ruess@b-tu.de}

	\begin{abstract}
	An extension of the multi-level \(hp\)-Finite Cell Method is proposed for the simulation of thermo-viscoplastic problems with temperature-dependent material behavior. The approach combines hierarchical adaptive refinement with a non-negative moment fitting (NNMF) quadrature scheme for efficient and robust integration of non-linear, history-dependent constitutive models on cut cells. The NNMF formulation yields sparse, positive quadrature rules that significantly reduce the number of integration points while maintaining stability and accuracy. \color{mygreenD} An error-indicator-driven \(h\)-refinement strategy enables localized resolution of pronounced solution gradients through successive adaptive refinement cycles. \color{black} The framework is implemented within a partitioned thermo-mechanical scheme that accounts for temperature-dependent mechanical behavior while restricting the thermal model to stationary heat conduction and neglecting heat generation due to inelastic dissipation. The methodology is evaluated on benchmark and application-oriented examples. The results demonstrate improved accuracy and substantial computational savings compared to the established adaptive space-tree integration.
\end{abstract}

\end{frontmatter}

\vfill
\noindent\textcopyright~2026. This manuscript version is made available under the
CC-BY-NC-ND 4.0 license \url{https://creativecommons.org/licenses/by-nc-nd/4.0/}

\newpage
\addtocontents{toc}{\protect\sloppy}
\tableofcontents
\newpage

\section{Introduction}
\label{sec:introduction}
The numerical prediction of industrial production processes in the steel industry has become critically important, particularly in the context of energy efficiency, cost optimization, and stringent quality control requirements. In processes such as continuous casting and hot forming, the material response is governed predominantly by thermo-viscoplastic behavior, posing significant challenges for numerical simulation. These simulations are inherently demanding, constrained by issues of robustness, accuracy, and substantial computational complexity.

The Finite Cell Method (FCM) has established itself as a high-order fictitious domain method, which is particularly well suited to analysis problems involving highly complex geometric and physical domains \cite{Duester:08.1,Parvizian:07.1,Schillinger:14.3}. In particular, recent developments regarding numerical integration of the governing equations and hierarchic adaptive refinement have contributed to an increasing efficiency of the method and optimized its application with regard to a reduced modeling and analysis effort. The following two extensions are the focus of this publication: (i) The introduction of a moment fitting quadrature scheme to the FCM by \emph{Joulaian et al.} \cite{Joulaian:16.1} initially has provided a good alternative in the context of linear problems to the then established sub-cell integration scheme based on recursive bisection \cite{Zander:14.1,Ruess:13.1}. 
For non-linear problems, the extension to non-negative moment fitting \cite{Hubrich:2019a} including selectively enriched and adaptive versions \cite{Duester:20.1,Duester:20.2} was an important step to ensure the required numerical stability at a reasonable number of quadrature points \cite{Legrain:2021,Garhuom:2022}. A much higher impact of this quadrature scheme is expected for material non-linear problems which need to draw on the deformation history stored at the quadrature points and where model refinement requires additional integration effort apart from the expense of cut cells. 
(ii) The concept of overlay refinements by mesh superposition has proven to be a powerful refinement strategy which circumvents hanging nodes \cite{Rank:92.1,Duester:07.1}. The development of the technique started with a single overlay solution, constructed from linear functions and has later been extended to make use of multiple overlay meshes \cite{Schillinger:12.1,Schillinger:11.1} and high-order shape functions on all meshes \cite{Zander:16.1,Zander:15.1,Zander:17.1}. In the context of the Finite Cell Method, the concept has shown to open the door to independent, hierarchical and adaptive model refinement, supporting the key features that characterize the idea of an embedding domain, the use of higher-order approximation and adaptive integration \cite{Kopp:2022}. 
The original concept was recently extended with anisotropic refinement capabilities, providing greater flexibility and substantially improved efficiency for problems exhibiting pronounced directional solution behavior \cite{Zander:22.1}.\\
A powerful local adaptive refinement scheme is of particular interest in the context of non-linear material behavior, where the progressive gradual change of the material generally starts locally and the propagation is sensitive to local effects \cite{Pamin:2017,Borst:1993,Ozcan:2019aa}. In recent decades, numerous material models have been developed and implemented in the finite element method to capture rate-independent and rate-dependent plasticity at high temperatures and temperature gradients some of which are critically reviewed in \cite{Nahrmann:2020,Chaboche:2008}.
We use a \emph{Perzyna}-type viscoplastic model combined with a linear coupling to thermal expansion for metals over a wide temperature range, following the formulation introduced by \emph{Oppermann et al.} \cite{Oppermann:2022}. In contrast to the original fully non-linear thermo-viscoplastic coupling, which accounts for temperature evolution due to viscoplastic dissipation, the present approach is restricted to the linear thermal expansion component, reflecting the numerically motivated scope of this work in which material modeling plays a subordinate role.
The applied model takes into account temperature-dependent material properties, work hardening and softening, as they occur in applications such as hot forming and continuous casting.\\

In this paper, we combine and extend the latest developments in the Finite Cell Method to provide a robust, efficient and reliable analysis framework for thermo-viscoplastic problems in the steel industry, allowing automated error-driven model refinement through successive adaptive refinement cycles. The key focus of the presented work is on efficiency aspects and further aspects of a self-controlling model refinement in the course of propagating plasticity. We apply the non-negative moment fitting approach and demonstrate its superiority over the established adaptive space-tree quadrature for the numerical integration of the governing equations using a series of examples of increasing complexity. Furthermore, it is shown that the robustness and numerical complexity of the non-negative moment fitting approach is essential when used in the context of hierarchical multi-level \(hp\)-refinement of non-linear problems requiring a sophisticated iterative solution. 
We equip the refinement scheme with a robust gradient-based error indicator that systematically drives adaptive model refinement between successive adaptive refinement cycles, ensuring efficient resource allocation while reliably resolving localized features of the physical response. 

We use a thermo-viscoplastic model which is capable of reliably capturing temperature-dependent parameters for non-linear hardening and softening and which is fit for the analysis of models from the steel industry. Since the primary objective of this work is the assessment of adaptive refinement and integration strategies, transient heat transport and heat generation due to inelastic dissipation are neglected. The thermal field is therefore treated as stationary and acts only through its influence on the mechanical response.
We demonstrate various aspects in depth using a number of selected benchmarks and examples of increasing complexity to evaluate the potential and limitations of the proposed approach for problems relevant to both industry and science.\\

The publication is organized as follows: in the following second section, we give a brief summary of the main aspects of the Finite Cell Method in combination with a non-negative moment fitting quadrature scheme in order to provide the framework for section three in which we reconsider an  adaptive hierarchic multi-level \(hp\)-refinement and introduce a simple gradient based error indicator for non-linear analyses. In section four, we specify the thermo-viscoplastic model which we use for a number of numerical test cases presented in section five. Finally, we summarize and conclude the main findings of this work in section six.

\section{Immersed boundary method with refined quadrature rule}
\label{sec:ImmersedLocalRefinement}
Immersed boundary methods show their particular strength when conventional discretization methods reach their limits with regard to the efficient and reliable generation of body-fitted meshes of complex solution domains. A variety of methods have been developed over time to open a new level of simulation quality and efficiency.
In this section, we focus on the Finite Cell Method (FCM) \cite{Schillinger:14.3,Varduhn:16.1,Duester:08.1} as one of the methods that gained much attention over the past decade in various fields of applications. 
The proposed analysis framework combines the Finite Cell Method with high-order finite elements and hierarchical overlay-mesh refinement, enabling efficient local adaptive refinement while retaining the accuracy of high-order discretizations \cite{Zander:15.1,Schillinger:11.1}.

\subsection{Modeling with Finite Cells}
\label{sec:FCM}

\begin{figure}[h!]
	\centering
	\includegraphics[width=0.94\textwidth]{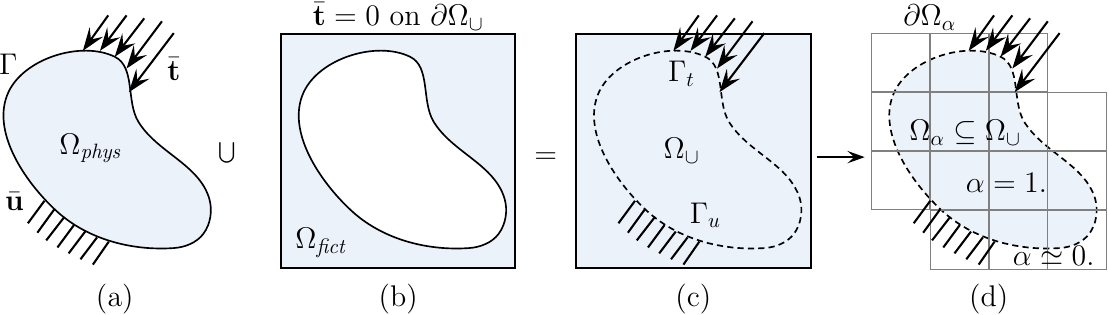}
	\caption{Concept of the Finite Cell Method: the boundary value problem with physical domain \(\Omega_\mathit{phys}\) and boundary conditions along \(\Gamma\) (a) is immersed in a fictitious extension domain of simple shape (b) forming the analysis domain \(\Omega_\cup\) (c). Discretization with finite cells reveals the finally considered analysis domain \(\Omega_\alpha\) (d).}
		\label{FCM_Schema}
\end{figure}

The discretization pipeline of the FCM is shown in Fig. \ref{FCM_Schema}. The problem domain \(\Omega_\mathit{phys}\) with natural and essential boundary conditions along \(\Gamma_t\) and  \(\Gamma_u\), respectively, is extended by a fictitious domain \(\Omega_\mathit{fict}\) (Figs. \ref{FCM_Schema}(a)-(c)), the unification of which results in the analysis domain $\Omega_\cup$ (Fig. \ref{FCM_Schema}(c)). Discretization of the analysis domain with high-order elements on a Cartesian grid (finite cells) reveals cells without contribution to the physical domain which are discarded from the analysis to obtain a better conditioned set of algebraic equations (Fig. \ref{FCM_Schema}(d))  \cite{Prenter:2017, Ruess:13.1}. 

\paragraph{The Finite Cell formulation in the context of stationary heat flow}
In analogy to structural mechanics, the Finite Cell Method is applied straightforwardly to stationary heat conduction problems by considering the temperature as the primary field variable and the heat flux as the constitutively related flux associated with the temperature gradient. For the purposes of the present study, the thermal problem is assumed to be quasi-stationary. This simplification allows the investigation to focus on the thermo-viscoplastic response and the performance of the proposed adaptive refinement and integration strategies rather than on transient heat transport effects.

The finite cells are essentially represented by high-order hexahedral elements which provide a numerical quadrature scheme allowing for integration of cells with cutting interface \(\Gamma\) between physical and fictitious domain. Our implementation utilizes the hierarchical $p$-version FEM approximation space formed by integrated Legendre polynomials \cite{Duester:08.1,Szabo:04.1}. The overall concept considers a penalization of the volume integrands in the fictitious extension domain with an indicator function \(\alpha(\mathbf{x})\) with value \(\alpha =1.0\) for points \(\mathbf{x}\in \Omega_{phys}\) and with a value \(\alpha \ll 1.0\) for points \(\mathbf{x}\in \Omega_\mathit{fict}\) \cite{Yang:2012b,Rank:12.1}. In the context of stationary heat flow, the weak form is modified according to conceptual changes of the FCM as follows. Find \(u\in H^1(\Omega)\), such that 
\begin{equationarray}{rclcl}
	\label{WFHeat}
	a(u,v) &=& l(v) &&\forall v\in H^1(\Omega)
\end{equationarray}
with 
\begin{equationarray}{rclcl}
	a(u,v) &=& \int_{\Omega_{\cup}} \alpha(\mathbf{x})\nabla v\cdot\kappa\nabla u\, d\Omega + \int_{\Gamma_{u}} \beta\, uv\,d\Gamma \\[4mm]
	l(v) &=& \int_{\Omega_{\cup}} v \bar{w}\,d\Omega  - \int_{\Gamma_{q}}v \bar{q}\,d\Gamma 
		+ \beta\int_{\Gamma_{u}} v \bar{u}\,d\Gamma 
\end{equationarray}
where \(u\) and \(v\) denote the unknown temperature and corresponding test function in the Sobolev space $H^1$, respectively, where \(\bar{w}\) denotes prescribed values for heat inflow per unit of time and volume and where \(\bar{q}\) and \(\bar{u}\) describe values for the heat outflow and temperature on the corresponding part of the domain surface \(\partial\Omega=\Gamma_u\cup\Gamma_q\) with \(\Gamma_u\cap\Gamma_q=\emptyset\), respectively. The thermal heat conduction coefficient is represented by \(\kappa\) and the value of \(\alpha\) for point locations belonging to \(\Omega_\mathit{fict}\) is chosen from the interval \([10^{-16},10^{-8}]\) depending on the problem under consideration (see e.g. \cite{Zander:12.1,Ruess:13.1}). The influence of $\alpha$ on the approximation error has been analyzed in detail in \cite{Dauge:2015}.
 
The approach requires a careful enforcement of essential boundary conditions on the surface \(\Gamma_u\) of the physical domain \(\Omega_\mathit{phys}\). Established concepts in this context comprise a penalty approach and a more sophisticated Nitsche-based concept \cite{Ruess:13.1}. In the present work, essential boundary conditions are imposed weakly using a penalty formulation with penalty parameter \(\beta\). For deeper insight to the details of the Finite Cell Method, the essential aspects and concepts, we refer the interested reader to \cite{Schillinger:14.3,Zander:14.1,Duester:08.1}.

\subsection{Non-negative moment fitting quadrature}
\label{sec:NNMFquadrature}
\begin{figure}[h!]
	\centering
	\includegraphics[width=0.8\textwidth]{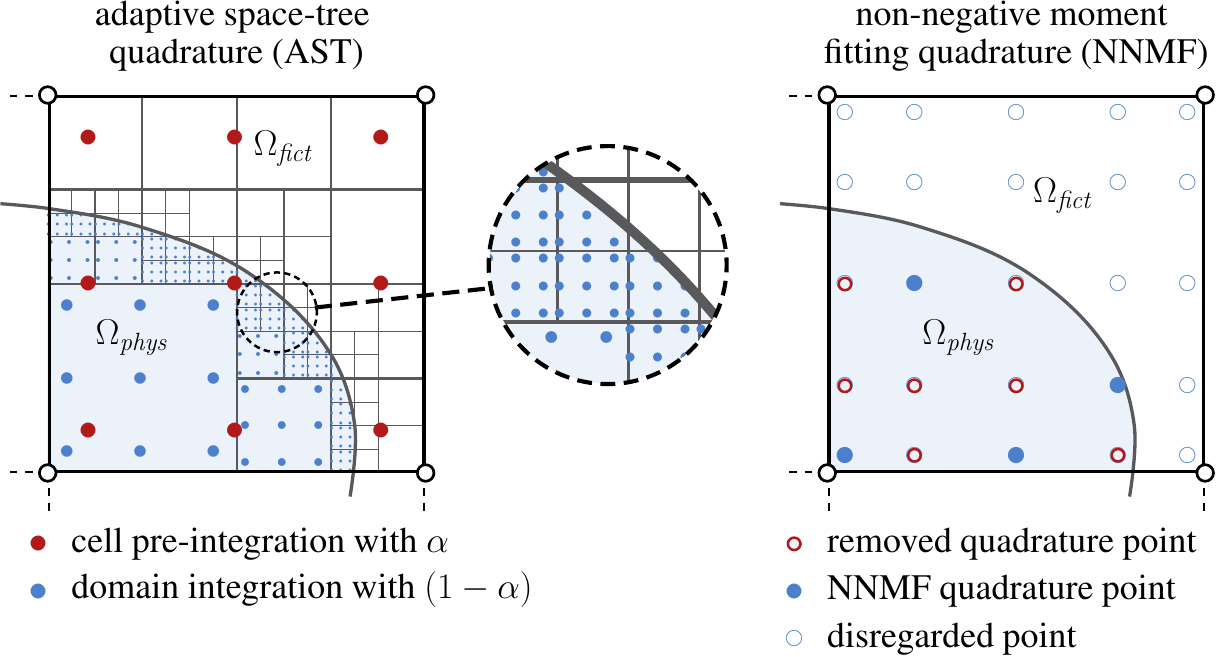}
	\caption{Quadrature scheme of the adaptive space-tree (AST) and the non-negative moment fitting (NNMF) approach, respectively. During the set up of the moment fitting quadrature, optimal quadrature points with corresponding weights are identified and remaining points are removed from an initial set of points.}
		\label{NNMF_Schema}
\end{figure}

The solution of the weak form \eqref{WFHeat} requires accurate integration of cut cells for which a number of reliable methods have been proposed over time, of which the adaptive space-tree technique (AST) has established itself as a popular and appropriate method \cite{Peto:20.1,Abedian:13.1,Ruess:13.1}. The  concept behind AST is a recursive bisection of the reference unit cell to create smaller and smaller integration sub-cells which adaptively capture the interface of the intersecting cell (Fig. \ref{NNMF_Schema}). Thus, the number of quadrature points increases in the direction of the intersection and refines the integration result. Although the AST approach has proven its reliability and accuracy for numerous problems, it suffers from its high computational complexity. The situation is even worse for non-linear problems, where an iterative solution is required at each step of the analysis to regain control of the equilibrium equations. 

In this context, quadrature schemes based on \emph{moment fitting} gained broad popularity \cite{Mueller:13.1,Sudhakar:13.1,Joulaian:16.1}.
The method reduces the number of quadrature points to a small but optimally distributed number of points without compromising the accuracy of the integration result. 
In Fig. \ref{NNMF_Schema}, the basic differences between the AST approach and the (non-negative) moment fitting are illustrated.

Following the general idea of numerical quadrature schemes, the moment fitting approach determines the location and corresponding weights of the quadrature points from the solution of a simpler secondary problem which applies pre-integration of an interpolation approach of corresponding order that fits the actual integration domain. The desired \(n\) quadrature weights  \(w_i\) and their locations \(\xi_i\) are determined from the non-linear system of equations: 
\begin{equationarray}{rcl}
	\label{MFEquation}
	\mathbf{A}\,\mathbf{w} &=& \mathbf{b}\\[2mm]\nonumber
	\begin{bmatrix}
		N_1(\xi_1) & \hdots & N_1(\xi_i) & \hdots & N_1(\xi_n)\\
		\vdots && \vdots && \vdots\\
		N_j(\xi_1) & \hdots & N_j(\xi_i) & \hdots & N_j(\xi_n)\\
		\vdots && \vdots && \vdots\\
		N_m(\xi_1) & \hdots & N_m(\xi_i) & \hdots & N_m(\xi_n)
	\end{bmatrix}
		
	\begin{bmatrix}
		w_1\\
		\vdots\\
		w_i\\
		\vdots\\
		w_n
	\end{bmatrix}

	&=&
		
	\begin{bmatrix}
		\int_{\Omega_\mathit{phys}}N_1(\xi)\,d\Omega\\
		\vdots\\
		\int_{\Omega_\mathit{phys}}N_j(\xi)\,d\Omega\\
		\vdots\\
		\int_{\Omega_\mathit{phys}}N_m(\xi)\,d\Omega
	\end{bmatrix} 
	
\end{equationarray}
in which \(N_j(\xi)\) denote a suitable set of \(m\) shape functions for appropriate domain interpolation, e.g. the set of integrated Legendre polynomials as proposed in \cite{Joulaian:16.1} and the integrals in \(\mathbf{b}\) are termed \emph{moments}, which gives the name of the method. The problem is linearized and thus significantly simplified by a basic selection of suitable quadrature point locations \(\xi_i\) which leads to a variety of versions, most of which are based on point locations of standard Gauss quadrature \cite{Duester:20.1,Hubrich:2019a,Duester:20.2,Joulaian:16.1}.

In the Finite Cell Method, a finite value of the indicator function \(\alpha(\mathbf{x})\) avoids vanishing stiffness contributions associated with basis functions supported predominantly in the fictitious domain and thereby improves the conditioning of the resulting algebraic system compared to the limiting case \(\alpha=0\). This conditioning effect is independent of the robustness of the non-linear solution procedure. While a finite fictitious-domain scaling parameter ensures a well-posed algebraic system, negative quadrature weights may still adversely affect the convergence behavior of Newton iterations in non-linear, history-dependent problems. This observation motivates the restriction of the moment-fitting problem to non-negative weights. 
The use of non-negative least-squares (NNLS) algorithms for the construction of cubature rules was originally proposed in \cite{Huybrechs:2009} and subsequently further developed in \cite{Sommariva:2015}. The approach was first applied to non-linear fictitious-domain analyses in \cite{Legrain:2021} and has since been adopted more broadly within the FCM community, e.g., in \cite{Garhuom:2022,Messmer:2024}.

According to \emph{Tchakaloff's theorem}, the non-linear moment fitting problem has solutions with $n=m$ positive weights, where $m$ denotes the number of moment-fitting constraints and $n$ the number of quadrature points \cite{Tchakaloff:1957,Davis:1967}. 
Due to the linearization with preselected point location, these solutions may not be found. Nevertheless, there is a set of \(n\geq m\) positive weights that can be found by minimizing the linearized equations \eqref{MFEquation} subject to the inequality constraint of non-negative weights.
\begin{equationarray}{rclclcl}
  \label{NNMFEquation}
	\mbox{min } \|\mathbf{A}\mathbf{w}-\mathbf{b}\|_2\qquad \mbox{subject to } \quad w_i\geq 0, \;i=1,\ldots, n
\end{equationarray}

\noindent It is worth noting the following properties:
\begin{itemize}
 \setlength{\parskip}{-3pt}
	\item[$\filledsquare$] In general, matrix \(\mathbf{A}\in\mathbb{R}^{m\times n}\) is non-symmetric and non-square.
	\item[$\filledsquare$] The existence of sparse solutions with a small fraction of non-zero weights has been proven in \cite{Davis:1967} if the initial set of quadrature points \(n\) is chosen sufficiently large.
	\item[$\filledsquare$] Optimization solvers like \emph{Trust Region Newton Methods} require many initial points \(n\gg m\) and often fail to provide sparse solutions.
\end{itemize} 

Herein, we follow the stable solution approach of \cite{Garhuom:2022}, in which the moment-fitting problem is solved using the non-negative least-squares algorithm of \emph{Lawson and Hanson} \cite{Lawson:1995}, yielding the desired sparse set of quadrature weights. 
The solver utilizes a QR-factorization to handle the rectangular matrix problem \cite{Ruess:2009a}. 
 
In the case of a cell that is intersected by the domain boundary or any other domain-separating interface, the following steps determine a set of optimally distributed quadrature points based on non-negative moment fitting:
\begin{enumerate}
	\item The moments \(b_i=\int_{\Omega_\mathit{phys}}N_j(\xi)\,d\Omega\) are approximated using the recursive AST approach in order to capture the physical domain integral within the finite cell accurately and robustly. 
Although recursive AST quadrature is computationally demanding, the associated cost is usually negligible compared to solving the NNLS problem (cf. Fig. \ref{PorousPlateRuntime} in section \ref{sec:porousPlate}).
	\item To set up the moment fitting equations, the physical cell domain is equipped with an initial choice of \(n\) Gauss-Legendre quadrature points. The number of points is increased step by step in various computations in order to improve the computed results for the weights until \(n>L m\) where \(L\) is empirically chosen from \(L=\{3,4,5,6\}\).
	Across the problems considered within our framework, stable solutions were obtained with \(L=3\) in the vast majority of cells. To rigorously limit computational cost, the parameter was capped at \(L=6\). In the rare instances where \(L=6\) failed to yield a convergent solution of the NNMF equations, the procedure reverted to the AST approach for the affected cell. 
	\item The solution of the moment fitting equations \eqref{NNMFEquation} using the solver described in \cite{Lawson:1995} provides the non-negative weights from which all non-zero weights are selected. The solution is accepted as converged if the residual norm of \eqref{NNMFEquation} satisfies a user-defined threshold: \(\|\mathbf{r}\|_2=\|\mathbf{b}- \mathbf{A}\,\mathbf{w}\|_2<\epsilon\) with \(\epsilon\) being the unit round-off value times a moderate constant of order \(\mathcal{O}(1)\). 
\end{enumerate}
The aforementioned steps consider the following aspects: (i) the selected choice of \(n\) Gauss-Legendre quadrature points is applied to the unit cell but can also be applied within an AST-sub-cell decomposition of the cut unit cell in order to densify quadrature points in the initialization step 2 and (ii) zero and negative weights from the solution step 2 are discarded from the process to ensure the desired non-negativity condition. 

Although the NNMF setup phase is computationally demanding, it is confined to cut cells only, thereby preserving favorable scalability in three-dimensional settings. Moreover, finite cell discretizations typically employ substantially larger elements with higher approximation order, in contrast to the much finer, low-order elements characteristic of conventional FEM, thereby reducing the number of degrees of freedom required to achieve a given approximation quality. 

\section{Error-indicator-controlled adaptive local refinement}
\label{sec:refinement}
Easy and efficient model refinement is a key factor for the success of any numerical approximation concept to ensure a proper convergence behavior for complex solution fields. In the context of the Finite Cell Method, a hierarchical approach of refinement according to the principles of superposition has proven to be a powerful concept in various respects \cite{Zander:16.1,Zander:17.2}. In the following, we apply and test this concept in the framework of plasticity and viscoplasticity analyses. We briefly summarize the concept implemented in our software framework, in which we introduce an error indicator based on the idea of \emph{Kelly} \cite{Kelly:1983} to control the entire refinement process.

\subsection{Local multi-level \(hp\)-refinement}
\label{sec:multi-levelhp}
The principles of the multi-level \(hp\)-refinement are illustrated in Fig. \ref{multilevelhp} revealing its basic characteristic which is the application of the superposition of locally refined and coupled overlay meshes. The shape functions of overlay meshes must (i) satisfy \(\mathcal{C}^0\)-continuity along the boundary at which the base and overlay meshes are coupled and (ii) ensure linear independence from the function space of the base mesh and the mesh of the underlying refinement level, respectively. Both requirements are met at the element level in order to ensure monotonic convergence. 

\begin{figure}[h!]
	\centering
	\begin{tabular}{c c c}
		\tikzsetnextfilename{scheme_multi_level_hp_1dR3}
\begin{tikzpicture}[font = \footnotesize]

	\begin{scope}[shift = {(0, 0)}, scale = 1.0]

		\draw[densely dotted] (0, 0) -- (0, 2.2);
		\draw[densely dotted] (1.5, 0) -- (1.5, 2.2);

		\draw[red, densely dashed] (0,0) -- (1.5,0);
		\draw[black] (1.5, 0) -- ++(1.5, 0);
		\node[anchor=east] at (0,0) {$k=0$};

		\begin{scope}[shift = {(0, 0)}]
			\draw[blue] (0, 0) plot[domain=0:1.5] ({\x}, {1.5 * 2/3 * \x});
		\end{scope}

		\begin{scope}[shift = {(1.5, 0)}]
			\draw[blue] (0, 0) plot[domain=0:1.5] ({\x}, {1.5 - 1.5 * 2/3 * \x});
			\draw[blue] (0, 0) plot[domain=0:1.5] ({\x}, {1.5 * 2/3 * \x});
			\draw[red] (0, 0) plot[domain=0:1.5] ({\x}, {1.5 * 1/4 * sqrt(6) * (pow((\x * 4/3 - 1),2) - 1)});
			\draw[red] (0, 0) plot[domain=0:1.5] ({\x}, {1.5 * 1/4 * sqrt(10) * (pow(\x * 4/3 - 1,2) - 1) * (\x * 4/3 - 1)});
			\draw[red] (0, 0) plot[domain=0:1.5] ({\x}, {1.5 * 1/16 * sqrt(14) * ( 5 * pow(\x * 4/3 - 1, 4) - 6 * pow(\x * 4/3 - 1, 2) + 1 )});
		\end{scope}

		\draw[red, fill=white] (0, 0) circle (1.5pt);
		\draw[black, fill=white] (1.5, 0) circle (1.5pt);
		\draw[black, fill=white] (3, 0) circle (1.5pt);


		\begin{scope}[shift = {(0, {0.2 + 2 * (1)})}]
			\draw[densely dotted] (1.5, 0) -- (1.5, .2+10/3-2.2);
			\draw[densely dotted] (0.75, 0) -- (0.75, .2+10/3-2.2);

			\draw[red, densely dashed] (0.75,0) -- (1.5,0);
			\draw[black] (0, 0) -- ++(0.75, 0);
			\node[anchor=east] at (0,0) {$k=1$};

			\begin{scope}[shift = {(0, 0)}, scale = 0.5]
				\begin{scope}[shift = {(0, 0)}]
					\draw[blue] (0, 0) plot[domain=0:1.5] ({\x}, {1.5 - 1.5 * 2/3 * \x});
					\draw[blue] (0, 0) plot[domain=0:1.5] ({\x}, {1.5 * 2/3 * \x});
					\draw[red] (0, 0) plot[domain=0:1.5] ({\x}, {1.5 * 1/4 * sqrt(6) * (pow((\x * 4/3 - 1),2) - 1)});
					\draw[red] (0, 0) plot[domain=0:1.5] ({\x}, {1.5 * 1/4 * sqrt(10) * (pow(\x * 4/3 - 1,2) - 1) * (\x * 4/3 - 1)});
				\end{scope}

				\begin{scope}[shift = {(1.5, 0)}]
					\draw[blue] (0, 0) plot[domain=0:1.5] ({\x}, {1.5 - 1.5 * 2/3 * \x});
				\end{scope}
			\end{scope}

			\draw[black, fill=white] (0, 0) circle (1.5pt);
			\draw[black, fill=white] (0.75, 0) circle (1.5pt);
			\draw[red, fill=white] (1.5, 0) circle (1.5pt);
		\end{scope}


		\begin{scope}[shift = {(0, {0.2 + 2 * (1 + 2/3)})}]
			\draw[black] (0.75, 0) -- ++(0.75, 0);
			\node[anchor=east] at (0.75,0) {$k=2$};

			\begin{scope}[shift = {(0.75, 0)}, scale = 0.25]
				\begin{scope}[shift = {(0, 0)}]
					\draw[blue] (0, 0) plot[domain=0:1.5] ({\x}, {1.5 * 2/3 * \x});
					\draw[red] (0, 0) plot[domain=0:1.5] ({\x}, {1.5 * 1/4 * sqrt(6) * (pow((\x * 4/3 - 1),2) - 1)});
				\end{scope}

				\begin{scope}[shift = {(1.5, 0)}]
					\draw[blue] (0, 0) plot[domain=0:1.5] ({\x}, {1.5 - 1.5 * 2/3 * \x});
					\draw[red] (0, 0) plot[domain=0:1.5] ({\x}, {1.5 * 1/4 * sqrt(6) * (pow((\x * 4/3 - 1),2) - 1)});
				\end{scope}
			\end{scope}

			\draw[red, fill=white, inner sep = 6] (0.75, 0) circle (1.5pt);
			\draw[black, fill=white, inner sep = 6] ({0.75 * 1.5}, 0) circle (1.5pt);
			\draw[red, fill=white, inner sep = 6] (1.5, 0) circle (1.5pt);
		\end{scope}

	\end{scope}

\end{tikzpicture} & 
		\tikzsetnextfilename{scheme_multi_level_hp_2dR3}
\begin{tikzpicture}[3d view={-30}{30}, line join=round, font = \footnotesize]

	\node at (0,0,-.5) {};

	\fill[lightgray] (0,0,0) -- (3,0,0) -- (3,3,0) -- (0,3,0) -- cycle;
	\fill[lightred] (1,2,0) -- (1,1,0) -- (2,1,0) -- (2,0,0) -- (3,0,0) -- (3,2,0) -- cycle;

	\draw[black] (0,0,0) -- (2,0,0) -- (2,1,0) -- (1,1,0) -- (1,2,0) -- (3,2,0) -- (3,3,0) -- (0,3,0) -- cycle
							 (1,0,0) -- (1,1,0)
							 (0,1,0) -- (1,1,0)
							 (0,2,0) -- (1,2,0)
							 (1,2,0) -- (1,3,0)
							 (2,2,0) -- (2,3,0);
	\draw[red, densely dashed]
		(2,2,0) -- (2,1,0)
		(2,1,0) -- (3,1,0)
		(2,0,0) -- (3,0,0) -- (3,2,0);

	\draw[black, fill=white]
    (0,0,0) circle (1.5pt)
    (0,1,0) circle (1.5pt)
    (0,2,0) circle (1.5pt)
    (0,3,0) circle (1.5pt)
    (1,0,0) circle (1.5pt)
    (1,1,0) circle (1.5pt)
    (1,2,0) circle (1.5pt)
    (1,3,0) circle (1.5pt)
    (2,0,0) circle (1.5pt)
    (2,1,0) circle (1.5pt)
    (2,2,0) circle (1.5pt)
    (2,3,0) circle (1.5pt)
    (3,2,0) circle (1.5pt)
    (3,3,0) circle (1.5pt);

	\draw[red, fill=white]
    (3,0,0) circle (1.5pt)
    (3,1,0) circle (1.5pt);

	\fill[gray, opacity=.3, shift={(.22,.35,0)}]
		(1,2.2,0) -- (1,1,0) -- (2,1,0) -- (2,0,0) -- (3-.22,0,0) -- (3-.22,2.2,0) -- cycle;

	\draw[densely dotted, black]
		(1,2,0) -- ++(0,0,1.5)
		(1,1,0) -- ++(0,0,1.5)
		(2,1,0) -- ++(0,0,1.5)
		(2,0,0) -- ++(0,0,1.5)
		(3,0,0) -- ++(0,0,1.5)
		(3,2,0) -- ++(0,0,1.5);

	\begin{scope}[shift={(0,0,1.5)}]
		\fill[lightgray] (1,2,0) -- (1,1,0) -- (2,1,0) -- (2,0,0) -- (3,0,0) -- (3,2,0) -- cycle;
		\fill[lightred] (1.5,1,0) -- (2,1,0) --	(2,0,0) -- (2.5,0,0) -- (2.5,1.5,0) -- (2,1.5,0) --	(2,2,0) -- (1.5,2,0) -- cycle;

		\draw[black] (3,2,0) -- (3,0,0) -- (2.5,0,0)
								 (1,1.5,0) -- (1.5,1.5,0)
								 (1.5,1,0) -- (1.5,2,0)
								 (2,1.5,0) -- (2,2,0)
								 (2.5,0.5,0) -- (3,0.5,0)
								 (2.5,1,0) -- (3,1,0)
								 (2.5,1.5,0) -- (2.5,2,0)
								 (2.5,1.5,0) -- (3,1.5,0)
								 (2,1.5,0) -- (2.5,1.5,0) -- (2.5,0,0);
		\draw[red, densely dashed] (3,2,0) -- (1,2,0) -- (1,1,0) -- (2,1,0) -- (2,0,0) -- (2.5,0,0)
															 (1.5,1.5,0) -- (2,1.5,0)
															 (2,1,0) -- (2.5,1,0)
															 (2,0.5,0) -- (2.5,0.5,0)
															 (2,1.5,0) -- (2,1,0);

		\draw[black, fill=white]
			(1.5,1.5,0) circle (1.5pt)
			(2,1.5,0) circle (1.5pt)
			(2.5,0,0) circle (1.5pt)
			(2.5,1.5,0) circle (1.5pt)
			(2.5,1,0) circle (1.5pt)
			(2.5,0.5,0) circle (1.5pt)
			(3,0,0) circle (1.5pt)
			(3,0.5,0) circle (1.5pt)
			(3,1,0) circle (1.5pt)
			(3,1.5,0) circle (1.5pt);

		\draw[red, fill=white]
			(3,2,0) circle (1.5pt)
			(2,2,0) circle (1.5pt)
			(1,2,0) circle (1.5pt)
			(1,1,0) circle (1.5pt)
			(2,1,0) circle (1.5pt)
			(2,0,0) circle (1.5pt)
			(1,1.5,0) circle (1.5pt)
			(1.5,2,0) circle (1.5pt)
			(2.5,2,0) circle (1.5pt)
			(2,0.5,0) circle (1.5pt)
			(1.5,1,0) circle (1.5pt);

		\fill[gray, opacity=.3, shift={(.2,.33)}] (1.5,1,0) -- (2,1,0) --	(2,0,0) -- (2.6,0,0) -- (2.6,1.6,0) -- (2.1,1.6,0) --	(2.1,2-.33,0) -- (1.5,2-.33,0) -- cycle;

		\draw[densely dotted, black]
			(1.5,1,0) -- ++(0,0,1.5)
			(2,0,0) -- ++(0,0,1.5)
			(2.5,0,0) -- ++(0,0,1.5)
			(2.5,1.5,0) -- ++(0,0,1.5)
			(2,2,0) -- ++(0,0,1.5)
			(1.5,2,0) -- ++(0,0,1.5);

	\end{scope}

	\begin{scope}[shift={(0,0,3)}]
		\fill[lightgray] (1.5,1,0) -- (2,1,0) --	(2,0,0) -- (2.5,0,0) -- (2.5,1.5,0) -- (2,1.5,0) --	(2,2,0) -- (1.5,2,0) -- cycle;
		\draw[red, densely dashed] (1.5,2,0) -- (1.5,1,0) -- (2,1,0) --	(2,0,0);
		\draw[red, densely dashed] (2.5,0,0) -- (2.5,1.5,0) -- (2,1.5,0) --	(2,2,0);
		\draw[red, densely dashed] (1.5,2,0) -- ++(.5,0,0);
		\draw[black]
			(1.5,1.75,0) -- ++(.5,0,0)
			(1.5,1.5,0) -- ++(.5,0,0)
			(1.5,1.25,0) -- ++(1,0,0)
			(2,1,0) -- ++(.5,0,0)
			(2,0.75,0) -- ++(.5,0,0)
			(2,0.5,0) -- ++(.5,0,0)
			(2,0.25,0) -- ++(.5,0,0)
			(2,0.0,0) -- ++(.5,0,0)
			(1.75,2,0) -- (1.75,1,0)
			(2,1,0) -- (2,1.5,0)
			(2.25,1.5,0) -- (2.25,0,0);

		\draw[red, fill=white]
			(1.5,2,0) circle (1.5pt)
			(1.5,1.75,0) circle (1.5pt)
			(1.5,1.5,0) circle (1.5pt)
			(1.5,1.25,0) circle (1.5pt)
			(1.5,1,0) circle (1.5pt)
			(1.75,1,0) circle (1.5pt)
			(2,1,0) circle (1.5pt)
			(2,0.75,0) circle (1.5pt)
			(2,0.5,0) circle (1.5pt)
			(2,0.25,0) circle (1.5pt)
			(2,0.0,0) circle (1.5pt);

		\draw[red, fill=white]
			(2,2,0) circle (1.5pt)
			(2,1.75,0) circle (1.5pt)
			(2,1.5,0) circle (1.5pt)
			(2.25,1.5,0) circle (1.5pt)
			(2.5,1.5,0) circle (1.5pt)
			(2.5,1.25,0) circle (1.5pt)
			(2.5,1.25,0) circle (1.5pt)
			(2.5,1,0) circle (1.5pt)
			(2.5,0.75,0) circle (1.5pt)
			(2.5,0.5,0) circle (1.5pt)
			(2.5,0.25,0) circle (1.5pt)
			(2.5,0.0,0) circle (1.5pt)
			(1.75,2,0) circle (1.5pt);

		\draw[black, fill=white]
			(1.75,1.75,0) circle (1.5pt)
			(1.75,1.5,0) circle (1.5pt)
			(1.75,1.25,0) circle (1.5pt)
			(2,1.25,0) circle (1.5pt)
			(2.25,1.25,0) circle (1.5pt)
			(2.25,1,0) circle (1.5pt)
			(2.25,0.75,0) circle (1.5pt)
			(2.25,0.5,0) circle (1.5pt)
			(2.25,0.25,0) circle (1.5pt)
			(2.25,0,0) circle (1.5pt);

	\end{scope}




\end{tikzpicture} &
		\tikzsetnextfilename{scheme_multi_level_hp_3dR3}
\begin{tikzpicture}[3d view={-30}{30}, line join=round, font = \footnotesize]

	\fill[lightgray] (0,0,0) -- (3,0,0) -- (3,3,0) -- (0,3,0) -- cycle
							(0,3,0) -- (0,3,-1) -- (0,0,-1) -- (2,0,-1) -- (2,0,0) -- (0,0,0) -- cycle;

	\fill[lightred] (1,2,0) -- (1,1,0) -- (2,1,0) -- (2,0,0) -- (3,0,0) -- (3,2,0) -- cycle
									(2,0,-1) -- (3,0,-1) -- (3,0,0) -- (2,0,0) -- cycle;

	\draw[black] (0,0,0) -- (2,0,0) -- (2,1,0) -- (1,1,0) -- (1,2,0) -- (3,2,0) -- (3,3,0) -- (0,3,0) -- cycle
							 (1,0,0) -- (1,1,0)
							 (0,1,0) -- (1,1,0)
							 (0,2,0) -- (1,2,0)
							 (1,2,0) -- (1,3,0)
							 (2,2,0) -- (2,3,0);
	\draw[black] (0,3,0) -- (0,3,-1) -- (0,0,-1) -- (2,0,-1) -- (2,0,0)
							 (0,2,0) -- (0,2,-1)
							 (0,1,0) -- (0,1,-1)
							 (0,0,0) -- (0,0,-1)
							 (1,0,0) -- (1,0,-1);
	\draw[red, densely dashed]
		(2,2,0) -- (2,1,0)
		(2,1,0) -- (3,1,0)
		(2,0,0) -- (3,0,0) -- (3,2,0)
		(2,0,-1) -- (3,0,-1) -- (3,0,0);

	\foreach \x in {0, 1, 2, 3}
		\foreach \y in {0, 1, 2, 3}
		{
			\draw[black, fill=white] (\x,\y,0) circle (1.5pt);
		}

	\draw[black, fill=white]
    (0,0,0) circle (1.5pt)
    (0,1,0) circle (1.5pt)
    (0,2,0) circle (1.5pt)
    (0,3,0) circle (1.5pt)
    (1,0,0) circle (1.5pt)
    (1,1,0) circle (1.5pt)
    (1,2,0) circle (1.5pt)
    (1,3,0) circle (1.5pt)
    (2,0,0) circle (1.5pt)
    (2,1,0) circle (1.5pt)
    (2,2,0) circle (1.5pt)
    (2,3,0) circle (1.5pt)
    (3,2,0) circle (1.5pt)
    (3,3,0) circle (1.5pt);

	\draw[red, fill=white]
    (3,0,0) circle (1.5pt)
    (3,1,0) circle (1.5pt);

	\draw[red, fill=white] (3,0,-1) circle (1.5pt);
	\draw[black, fill=white] (2,0,-1) circle (1.5pt);
	\draw[black, fill=white] (1,0,-1) circle (1.5pt);
	\draw[black, fill=white] (0,0,-1) circle (1.5pt);
	\draw[black, fill=white] (0,1,-1) circle (1.5pt);
	\draw[black, fill=white] (0,2,-1) circle (1.5pt);
	\draw[black, fill=white] (0,3,-1) circle (1.5pt);

	\fill[gray, opacity=.3, shift={(.22,.35,0)}]
		(1,2.2,0) -- (1,1,0) -- (2,1,0) -- (2,0,0) -- (3-.22,0,0) -- (3-.22,2.2,0) -- cycle;

	\draw[densely dotted, black]
		(1,2,0) -- ++(0,0,0.8)
		(1,1,0) -- ++(0,0,0.8)
		(2,1,0) -- ++(0,0,0.8)
		(2,0,0) -- ++(0,0,0.8)
		(3,0,0) -- ++(0,0,0.8)
		(3,2,0) -- ++(0,0,0.8);

	\begin{scope}[shift={(0,0,1.8)}]

		\fill[lightgray] (1,2,0) -- (1,1,0) -- (2,1,0) -- (2,0,0) -- (3,0,0) -- (3,2,0) -- cycle;
		\fill[lightred] (1.5,1,0) -- (2,1,0) --	(2,0,0) -- (2.5,0,0) -- (2.5,1.5,0) -- (2,1.5,0) --	(2,2,0) -- (1.5,2,0) -- cycle
										(1,2,0) -- (1,2,-1) -- (1,1,-1) -- (2,1,-1) -- (2,0,-1) -- (3,0,-1) -- (3,0,0) -- (2,0,0) -- (2,1,0) -- (1,1,0) -- cycle;
		\fill[lightgray] (2,0,-1) -- (3,0,-1) -- (3,0,0) -- (2.5,0,0) -- (2.5,0,-.5) -- (2,0,-.5) --  cycle;

		\draw[black] (3,2,0) -- (3,0,0) -- (2.5,0,0)
								 (1,1.5,0) -- (1.5,1.5,0)
								 (1.5,1,0) -- (1.5,2,0)
								 (2,1.5,0) -- (2,2,0)
								 (2.5,0.5,0) -- (3,0.5,0)
								 (2.5,1,0) -- (3,1,0)
								 (2.5,1.5,0) -- (2.5,2,0)
								 (2.5,1.5,0) -- (3,1.5,0)
								 (2,1.5,0) -- (2.5,1.5,0) -- (2.5,0,0)
								 (2.5,0,-1) -- ++(0,0,1)
								 (2,0,-.5) --  (3,0,-.5)
								 (2,0,-1) -- (3,0,-1) -- (3,0,0);
		\draw[red, densely dashed] (3,2,0) -- (1,2,0) -- (1,1,0) -- (2,1,0) -- (2,0,0) -- (2.5,0,0)
															 (1.5,1.5,0) -- (2,1.5,0)
															 (2,1,0) -- (2.5,1,0)
															 (2,0.5,0) -- (2.5,0.5,0)
															 (2,1.5,0) -- (2,1,0)
															 (1,2,0) -- (1,2,-1) -- (1,1,-1) -- (2,1,-1) -- (2,0,-1)
															 (1,2,-.5) -- (1,1,-.5) -- (2,1,-.5) -- (2,0,-.5)
															 (1,1.5,-1) -- ++(0,0,1)
															 (1.5,1,-1) -- ++(0,0,1)
															 (2,.5,-1) -- ++(0,0,1)
															 (1,1,-1) -- (1,1,0)
															 (2,1,-1) -- (2,1,0)
															 (2,0,-1) -- (2,0,0);

		\draw[black, fill=white]
			(2.5,0,0) circle (1.5pt)
			(2.5,0,-1) circle (1.5pt)
			(3,0,-.5) circle (1.5pt)
			(2.5,0,-.5) circle (1.5pt)
			(1.5,1.5,0) circle (1.5pt)
			(2,1.5,0) circle (1.5pt)
			(2.5,1.5,0) circle (1.5pt)
			(2.5,1,0) circle (1.5pt)
			(2.5,0.5,0) circle (1.5pt)
			(3,0.5,0) circle (1.5pt)
			(3,0,0) circle (1.5pt)
			(3,0,-1) circle (1.5pt)
			(3,1,0) circle (1.5pt)
			(3,1.5,0) circle (1.5pt);

		\draw[red, fill=white]
			(2.5,2,0) circle (1.5pt)
			(3,2,0) circle (1.5pt)
			(2,2,0) circle (1.5pt)
			(1,2,0) circle (1.5pt)
			(1,1,0) circle (1.5pt)
			(2,1,0) circle (1.5pt)
			(2,0,0) circle (1.5pt)
			(1,1.5,0) circle (1.5pt)
			(1.5,2,0) circle (1.5pt)
			(2,0.5,0) circle (1.5pt)
			(1.5,1,0) circle (1.5pt)
			(1,2,-.5) circle (1.5pt)
			(1,1.5,-.5) circle (1.5pt)
			(1,1,-.5) circle (1.5pt)
			(1.5,1,-.5) circle (1.5pt)
			(2,1,-.5) circle (1.5pt)
			(2,.5,-.5) circle (1.5pt)
			(2,0,-.5) circle (1.5pt)
			(1,2,-1) circle (1.5pt)
			(1,1.5,-1) circle (1.5pt)
			(1,1,-1) circle (1.5pt)
			(1.5,1,-1) circle (1.5pt)
			(2,1,-1) circle (1.5pt)
			(2,.5,-1) circle (1.5pt)
			(2,0,-1) circle (1.5pt);

		\fill[gray, opacity=.3, shift={(.2,.33)}] (1.5,1,0) -- (2,1,0) --	(2,0,0) -- (2.6,0,0) -- (2.6,1.6,0) -- (2.1,1.6,0) --	(2.1,2-.33,0) -- (1.5,2-.33,0) -- cycle;

		\draw[densely dotted, black]
			(1.5,1,0) -- ++(0,0,.8)
			(2,0,0) -- ++(0,0,.8)
			(2.5,0,0) -- ++(0,0,.8)
			(2.5,1.5,0) -- ++(0,0,.8)
			(2,2,0) -- ++(0,0,.8)
			(1.5,2,0) -- ++(0,0,.8);
	\end{scope}

	\begin{scope}[shift={(0,0,3.1)}]
		\fill[lightgray] (1.5,1,0) -- (2,1,0) --	(2,0,0) -- (2.5,0,0) -- (2.5,1.5,0) -- (2,1.5,0) --	(2,2,0) -- (1.5,2,0) -- cycle
								(2,0,0) -- (2.5,0,0) -- (2.5,0,-.5) -- (2,0,-.5) -- cycle;
		\fill[lightred] (1.5,2,0) -- (1.5,2,-.5) -- (1.5,1,-.5) -- (2,1,-.5) --	(2,0,-.5) -- (2,0,0) -- (2,1,0) -- (1.5,1,0) -- cycle;

		\draw[red, densely dashed]
			(1.5,2,0) -- (1.5,1,0) -- (2,1,0) --	(2,0,0)
			(1.5,2,0) -- (1.5,2,-.5) -- (1.5,1,-.5) -- (2,1,-.5) --	(2,0,-.5) -- (2.5,0,-.5) -- (2.5,0,0)
			(1.5,2,-.25) -- (1.5,1,-.25) -- (2,1,-.25) --	(2,0,-.25)
			(2.5,0,0) -- (2.5,1.5,0) -- (2,1.5,0) --	(2,2,0);
		\draw[red, densely dashed]
			(1.5,2,0) -- ++(.5,0,0)
			(1.5,1.75,0) -- ++(0,0,-0.5)
			(1.5,1.5,0) -- ++(0,0,-0.5)
			(1.5,1.25,0) -- ++(0,0,-0.5)
			(1.5,1,0) -- ++(0,0,-0.5)
			(1.75,1,0) -- ++(0,0,-0.5)
			(2,1,0) -- ++(0,0,-0.5)
			(2,0.75,0) -- ++(0,0,-0.5)
			(2,0.5,0) -- ++(0,0,-0.5)
			(2,0.25,0) -- ++(0,0,-0.5)
			(2,0.0,0) -- ++(0,0,-0.5);

		\draw[black]
			(2,0,-.25) -- ++(.5,0,0)
			(1.5,1.75,0) -- ++(.5,0,0)
			(1.5,1.5,0) -- ++(.5,0,0)
			(1.5,1.25,0) -- ++(1,0,0)
			(2,1,0) -- ++(.5,0,0)
			(2,0.75,0) -- ++(.5,0,0)
			(2,0.5,0) -- ++(.5,0,0)
			(2,0.25,0) -- ++(.5,0,0)
			(2,0.0,0) -- ++(.5,0,0)
			(1.75,2,0) -- (1.75,1,0)
			(2,1,0) -- (2,1.5,0)
			(2.25,1.5,0) -- (2.25,0,0)
			(2.25,0,0) -- ++(0,0,-0.5);

		\draw[red, fill=white]
			(1.5,2,0) circle (1.5pt)
			(1.5,1.75,0) circle (1.5pt)
			(1.5,1.5,0) circle (1.5pt)
			(1.5,1.25,0) circle (1.5pt)
			(1.5,1,0) circle (1.5pt)
			(1.75,1,0) circle (1.5pt)
			(2,1,0) circle (1.5pt)
			(2,0.75,0) circle (1.5pt)
			(2,0.5,0) circle (1.5pt)
			(2,0.25,0) circle (1.5pt)
			(2,0.0,0) circle (1.5pt);
		\draw[red, fill=white, shift={(0,0,-.25)}]
			(1.5,2,0) circle (1.5pt)
			(1.5,1.75,0) circle (1.5pt)
			(1.5,1.5,0) circle (1.5pt)
			(1.5,1.25,0) circle (1.5pt)
			(1.5,1,0) circle (1.5pt)
			(1.75,1,0) circle (1.5pt)
			(2,1,0) circle (1.5pt)
			(2,0.75,0) circle (1.5pt)
			(2,0.5,0) circle (1.5pt)
			(2,0.25,0) circle (1.5pt)
			(2,0.0,0) circle (1.5pt)
			(2.5,0,0) circle (1.5pt);
		\draw[red, fill=white, shift={(0,0,-.5)}]
			(1.5,2,0) circle (1.5pt)
			(1.5,1.75,0) circle (1.5pt)
			(1.5,1.5,0) circle (1.5pt)
			(1.5,1.25,0) circle (1.5pt)
			(1.5,1,0) circle (1.5pt)
			(1.75,1,0) circle (1.5pt)
			(2,1,0) circle (1.5pt)
			(2,0.75,0) circle (1.5pt)
			(2,0.5,0) circle (1.5pt)
			(2,0.25,0) circle (1.5pt)
			(2,0.0,0) circle (1.5pt)
			(2.25,0,0) circle (1.5pt)
			(2.5,0,0) circle (1.5pt);

		\draw[red, fill=white]
			(1.75,2,0) circle (1.5pt)
			(2,2,0) circle (1.5pt)
			(2,1.75,0) circle (1.5pt)
			(2,1.5,0) circle (1.5pt)
			(2.25,1.5,0) circle (1.5pt)
			(2.5,1.5,0) circle (1.5pt)
			(2.5,1.25,0) circle (1.5pt)
			(2.5,1.25,0) circle (1.5pt)
			(2.5,1,0) circle (1.5pt)
			(2.5,0.75,0) circle (1.5pt)
			(2.5,0.5,0) circle (1.5pt)
			(2.5,0.25,0) circle (1.5pt)
			(2.5,0.0,0) circle (1.5pt);

		\draw[black, fill=white]
			(1.75,1.75,0) circle (1.5pt)
			(1.75,1.5,0) circle (1.5pt)
			(1.75,1.25,0) circle (1.5pt)
			(2,1.25,0) circle (1.5pt)
			(2.25,1.25,0) circle (1.5pt)
			(2.25,1,0) circle (1.5pt)
			(2.25,0.75,0) circle (1.5pt)
			(2.25,0.5,0) circle (1.5pt)
			(2.25,0.25,0) circle (1.5pt)
			(2.25,0,-.25) circle (1.5pt)
			(2.25,0,0) circle (1.5pt);

	\end{scope}



\end{tikzpicture} \\
		$1$D & $2$D & $3$D \\
	\end{tabular}
	\caption{Principle of the multi-level \(hp\) FCM in one, two and three dimensions. The graphic denoted 1D depicts the shape functions applied on each refinement level. The graphics labeled 2D and 3D show the individual topological components to which shape functions are assigned. No shape functions are assigned to the components in red to ensure linear independence from the function space of the underlying level and overall \(\mathcal{C}^0\)-continuity, all other components share usual shape functions.}
	\label{multilevelhp}
\end{figure}

\noindent The implemented formalism for the generation of overlay meshes leads to an element-wise introduction of a suitable set of shape functions as summarized by the following steps:
\begin{itemize}
\renewcommand\labelitemi{$\filledsquare$}
	\item Each element to be refined is bisected isotropically resulting in \(2^d\) overlays for an element of dimension \(d\) on the following mesh level.
	\item All topological components such as nodes, edges and faces which define interior boundaries of the considered mesh level are marked by a flag. The marked components are not equipped with shape functions in the following solution process to ensure global \(\mathcal{C}^0\)-continuity.
	\item All topological components which are refined through overlaying meshes are disabled, if their overlays are not disabled due to continuity constraints. This way linear independence between the shape functions of all levels is ensured.
	\item Whereas the original multi-level \(hp\)-method (\emph{Zander}, \cite{Zander:15.1}) keeps nodes active on the lowest possible level, our implementation follows (\emph{Kopp}, \cite{Kopp:2022}) by keeping them active on the highest possible level. This eliminates the need for any special treatment of nodes and allows them to be handled in the same way as all other topological components.
\end{itemize}

For a detailed description of the refinement principles in the framework of the Finite Cell Method we refer the interested reader to \cite{Zander:17.1, Zander:16.1, Zander:17.2, Schillinger:11.4} or in a different context of isogeometric analysis to \cite{Schillinger:12.1}.

\subsection{Refinement strategy}
\label{sec:refinementstrategy}
The need for $p$- or $h$-refinement is generally not apparent \emph{a priori} and depends strongly on the characteristics of the problem under consideration. Consequently, \emph{a posteriori} error indicators remain a key ingredient of automated adaptive refinement procedures, as they identify regions with insufficient approximation quality. For cut-cell discretizations, residual-based error estimation has been investigated by \emph{Di Stolfo et al.} \cite{DiStolfo:2021}, who proposed a reliable residual-based estimator for the Finite Cell Method. Complementary to this work, \emph{Sun et al.} \cite{Sun:2020} developed an implicit \emph{a posteriori} error estimation strategy for cut finite elements, demonstrating its suitability for driving adaptive refinement.

The standard adaptive procedure for linear analyses, comprising (i) an initial coarse-mesh solution, (ii) the determination of element-wise a posteriori error indicators, and (iii) a refined mesh solution, can be transferred to the non-linear problems considered in this work. 
For history-dependent material models, however, refinement or coarsening between
individual load steps would require a careful transfer of the current deformation state and all history-dependent variables to the new discretization~\cite{Sartorti:2023}. To avoid such intermediate state transfers, the adaptive procedure employed in the present work is performed in successive refinement cycles, each comprising the complete loading history on a fixed
discretization. After the final load step, the error indicator is evaluated and the mesh is refined accordingly. The complete analysis is then restarted from the initial unloaded state on the refined discretization. Hence, all history-dependent variables are initialized and evolved directly on the respective discretization, and no transfer of internal variables between successive meshes is required. This procedure is repeated for the prescribed number of adaptive refinement cycles.

The error indicator used herein quantifies the need for refinement with a scalar \(\eta(\Omega_e)\) for every element domain \(\Omega_e\), similar to the approach proposed by \emph{Kelly et al.} \cite{Kelly:1983}:
\begin{equationarray}{rclcl}
	\label{Kelly}\eta(\Omega_e) &=& \left[\frac{1}{|\partial \hat{\Omega}_e|}\Bigg\{\sum_{i=1}^{n_i}\sum_{j=1}^{n_j}\Big(\int_{\partial \hat{\Omega}_e} \Big\llbracket \frac{\partial u_i}{\partial x_j} \Big\rrbracket^2 d(\partial\hat{\Omega}_e) \Big)  \Bigg\}\right]^{\raisebox{1.3ex}{$\frac{1}{2}$}}
\end{equationarray}
in which \(n_i\) denotes the number of components of the primary field and \(n_j\) the geometric dimension. The domain \(\partial\hat{\Omega}_e = \partial\Omega_e\cap\Omega_\mathit{phys} \) denotes the cell boundary and \(\llbracket \cdot \rrbracket\) denotes the jump operator characterizing the displacement-gradient jump across adjacent element boundaries. 
Note that only the portion of the cell boundary located inside the physical domain \(\Omega_\mathit{phys}\) contributes to the boundary integral. The corresponding boundary integration is performed using the same quadrature strategy as employed for the associated cut cell, i.e., adaptive space-tree quadrature or NNMF, respectively.

In general, we assume that equation \eqref{Kelly} returns large values where the solution is not particularly smooth and requires further local \(h\)-refinement. 
The multi-level \(hp\)-formulation generally allows the polynomial order to be specified independently on each overlay level. In the adaptive strategy employed herein, however, the polynomial orders are prescribed a priori, while adaptivity is restricted to hierarchical \(h\)-refinement through the addition of overlay levels.

To eliminate discretization bias across heterogeneous cell sizes and successive refinement levels, the error estimator $\eta(\Omega_e)$ was normalized with respect to the cell boundary. For the one-dimensional example in Section~\ref{sec:oneDWMaterialJump}, refinement is restricted to the cell exhibiting the highest value of $\eta(\Omega_e)$. For the examples in Sections~\ref{sec:twoDContinuousMaterial} and~\ref{sec:plateWH}, we employed a fixed-fraction marking strategy, restricting refinement to the \(30\%\) of cells exhibiting the highest values of $\eta(\Omega_e)$. For the porous plate example in Section~\ref{sec:porousPlate}, a D\"orfler marking strategy~\cite{Doerfler1996} with a threshold of $\theta=0.4$ was employed instead. The respective marking parameters were selected empirically to provide a suitable balance between solution accuracy and computational effort.\\

\subsection{One-dimensional model with material kink}
\label{sec:oneDWMaterialJump}
In the following, the proposed refinement strategy is tested using a one- and a two-dimensional benchmark problem that simulates the weakly discontinuous material behavior at the elasto-plastic transition front in plastic deformation analysis.

The governing differential equation of the considered one-dimensional elasticity problem reads:
\begin{equationarray}{rclclclclc}
	\frac{d}{dx} \left[EA\frac{du}{dx} \right] &=& 0 \quad \mbox{ with }\quad \ & E(x) &=&
		\begin{cases}
			E_0+\frac{\Delta E}{a}(a-x) \quad \ & 0\leq x<a\\
			E_0\quad \ & a\leq x\leq b
		\end{cases} \nonumber \\[6mm]
		&&&\forall x\in \Omega &=&[0,a)\cup[a,b] \mbox{ and } 0<a<b 
\end{equationarray}
where \(\Delta E\) controls the linear variation of Young's modulus in \(x\in[0,a)\), while \(E=E_0\) is prescribed in \(x\in[a,b]\). The displacement solution \(u(x)\) should satisfy the following boundary conditions:
\begin{equationarray}{rcl}
	u(x=0) &=& 0\\[2mm]
	\left. E\frac{du}{dx}\right|_{x=b} &=& t_0\;
\end{equationarray}
which leads to the closed-form reference solution provided in Appendix A. The parameters that determine the geometry, the material and the load properties were selected as $a = \tfrac{2}{3}$, $b= 1$, $t_0=1$, $A=1$, $E_0 = 10$ and $\Delta E=90$.

\begin{figure}[h!]
	\centering
	\tikzsetnextfilename{solution_c0_material_benchmarmk_1d}
\begin{tikzpicture}[spy using outlines={black, rectangle, size=2cm, magnification=10, connect spies}]

	\begin{axis} [
		shift = {(0, 0)},
		font={\footnotesize},
		axis lines = box,
		title = {$k = 0$, $p = 20$ with $21$ DoFs},
		xticklabel = \empty,
		ylabel = Displacement $u$,
		scaled y ticks = false,
		width = 8cm,
		height = 6cm,
		max space between ticks = 80,
		minor x tick num = 0,
		minor y tick num = 0,
		cycle list name = color_list_mark,
		legend pos = north west,
		legend style = { nodes = { scale = 0.75, transform shape }, },
		legend cell align = { left },
		domain = 0:1,
		ymin = -0.007,
		ymax =  0.058,
	]

		\addplot[ lightgray, line width = 0.05mm ] coordinates { (0, -0.007) (0, 0.058) };
		\addplot[ lightgray, line width = 0.05mm ] coordinates { (1, -0.007) (1, 0.058) };

		\addplot[ black, domain = 0:0.66666 ] { 0.66666 * ( ln( 0.66666 * 100 ) - ln( 0.66666 * 100 - 90 * x ) ) / 90 };
		\addplot[ black, domain = 0.66666:1 ] { 0.66666 * ( ln( 0.66666 * 100 ) - ln( 0.66666 * 10 ) ) / 90 + ( x - 0.66666 ) / 10 };

		\addplot[ blue ] table [x = x, y = u, col sep = comma] {solution_c0_material_benchmark_1d_pfem_p20R2.csv};

		\legend{
		};

	\end{axis}

	\begin{axis} [
		shift = {(0, -4.5cm)},
		font={\footnotesize},
		axis lines = box,
		xlabel = $x$,
		ylabel = Strain $\varepsilon$,
		scaled y ticks = false,
		width = 8cm,
		height = 6cm,
		max space between ticks = 80,
		minor x tick num = 0,
		minor y tick num = 0,
		cycle list name = color_list_mark,
		legend pos = north west,
		legend style = { nodes = { scale = 0.75, transform shape }, },
		legend cell align = { left },
		domain = 0:1,
		ymin = 0,
		ymax = 0.11,
	]

		\addplot[ lightgray, line width = 0.05mm ] coordinates { (0, 0) (0, 0.11) };
		\addplot[ lightgray, line width = 0.05mm ] coordinates { (1, 0) (1, 0.11) };

		\addplot[ black, domain = 0:0.66666 ] { 0.66666 / ( 0.66666 * 100 - 90 * x ) };
		\addplot[ black, domain = 0.66666:1 ] { 0.1 };

		\addplot[ blue ] table [x = x, y = eps, col sep = comma] {solution_c0_material_benchmark_1d_pfem_p20R2.csv};

		\legend{
		};

	\end{axis}

	\begin{axis} [
		shift = {(6.5cm, 0)},
		font={\footnotesize},
		axis lines = box,
		title = {$k = 14$, $p \in [2, 15]$ with $929$ DoFs},
		xticklabel = \empty,
		yticklabel = \empty,
		scaled y ticks = false,
		width = 8cm,
		height = 6cm,
		max space between ticks = 80,
		minor x tick num = 0,
		minor y tick num = 0,
		cycle list name = color_list_mark,
		legend pos = north west,
		legend style = { nodes = { scale = 0.75, transform shape }, },
		legend cell align = { left },
		domain = 0:1,
		ymin = -0.007,
		ymax =  0.058,
	]

		\addplot[ lightgray, line width = 0.05mm ] coordinates { (0, -0.007) (0, 0.058) };

		\addplot[ lightgray, line width = 0.05mm ] coordinates { ( 0.0, -1 ) ( 0.0, 1 ) };
		\addplot[ lightgray, line width = 0.05mm ] coordinates { ( 0.5, -1 ) ( 0.5, 1 ) };
		\addplot[ lightgray, line width = 0.05mm ] coordinates { ( 0.625, -1 ) ( 0.625, 1 ) };
		\addplot[ lightgray, line width = 0.05mm ] coordinates { ( 0.65625, -1 ) ( 0.65625, 1 ) };
		\addplot[ lightgray, line width = 0.05mm ] coordinates { ( 0.6640625, -1 ) ( 0.6640625, 1 ) };
		\addplot[ lightgray, line width = 0.05mm ] coordinates { ( 0.666015625, -1 ) ( 0.666015625, 1 ) };
		\addplot[ lightgray, line width = 0.05mm ] coordinates { ( 0.666259765625, -1 ) ( 0.666259765625, 1 ) };
		\addplot[ lightgray, line width = 0.05mm ] coordinates { ( 0.6663818359375, -1 ) ( 0.6663818359375, 1 ) };
		\addplot[ lightgray, line width = 0.05mm ] coordinates { ( 0.66650390625, -1 ) ( 0.66650390625, 1 ) };
		\addplot[ lightgray, line width = 0.05mm ] coordinates { ( 0.66656494140625, -1 ) ( 0.66656494140625, 1 ) };
		\addplot[ lightgray, line width = 0.05mm ] coordinates { ( 0.6666259765625, -1 ) ( 0.6666259765625, 1 ) };
		\addplot[ lightgray, line width = 0.05mm ] coordinates { ( 0.66668701171875, -1 ) ( 0.66668701171875, 1 ) };
		\addplot[ lightgray, line width = 0.05mm ] coordinates { ( 0.666748046875, -1 ) ( 0.666748046875, 1 ) };
		\addplot[ lightgray, line width = 0.05mm ] coordinates { ( 0.66680908203125, -1 ) ( 0.66680908203125, 1 ) };
		\addplot[ lightgray, line width = 0.05mm ] coordinates { ( 0.6668701171875, -1 ) ( 0.6668701171875, 1 ) };
		\addplot[ lightgray, line width = 0.05mm ] coordinates { ( 0.6669921875, -1 ) ( 0.6669921875, 1 ) };
		\addplot[ lightgray, line width = 0.05mm ] coordinates { ( 0.66796875, -1 ) ( 0.66796875, 1 ) };
		\addplot[ lightgray, line width = 0.05mm ] coordinates { ( 0.671875, -1 ) ( 0.671875, 1 ) };
		\addplot[ lightgray, line width = 0.05mm ] coordinates { ( 0.6875, -1 ) ( 0.6875, 1 ) };
		\addplot[ lightgray, line width = 0.05mm ] coordinates { ( 0.75, -1 ) ( 0.75, 1 ) };

		\addplot[ lightgray, line width = 0.05mm ] coordinates { (1, -0.007) (1, 0.058) };

		\addplot[ black, domain = 0:0.66666 ] { 0.66666 * ( ln( 0.66666 * 100 ) - ln( 0.66666 * 100 - 90 * x ) ) / 90 };
		\addplot[ black, domain = 0.66666:1 ] { 0.66666 * ( ln( 0.66666 * 100 ) - ln( 0.66666 * 10 ) ) / 90 + ( x - 0.66666 ) / 10 };

		\addplot[ orange ] table [x = x, y = u, col sep = comma] {solution_c0_material_benchmark_1d_kelly_k14R2.csv};

		\legend{
		};

	\end{axis}

	\begin{axis} [
		shift = {(6.5cm, -4.5cm)},
		font={\footnotesize},
		axis lines = box,
		xlabel = $x$,
		yticklabels = \empty,
		width = 8cm,
		height = 6cm,
		max space between ticks = 80,
		minor x tick num = 0,
		minor y tick num = 0,
		cycle list name = color_list_mark,
		legend pos = north west,
		legend style = { nodes = { scale = 0.75, transform shape }, },
		legend cell align = { left },
		domain = 0:1,
		ymin = 0.0,
		ymax = 0.11,
	]

		\addplot[ lightgray, line width = 0.05mm ] coordinates { (0, 0) (0, 0.11) };

		\addplot[ lightgray, line width = 0.05mm ] coordinates { ( 0.0, -1 ) ( 0.0, 1 ) };
		\addplot[ lightgray, line width = 0.05mm ] coordinates { ( 0.5, -1 ) ( 0.5, 1 ) };
		\addplot[ lightgray, line width = 0.05mm ] coordinates { ( 0.625, -1 ) ( 0.625, 1 ) };
		\addplot[ lightgray, line width = 0.05mm ] coordinates { ( 0.65625, -1 ) ( 0.65625, 1 ) };
		\addplot[ lightgray, line width = 0.05mm ] coordinates { ( 0.6640625, -1 ) ( 0.6640625, 1 ) };
		\addplot[ lightgray, line width = 0.05mm ] coordinates { ( 0.666015625, -1 ) ( 0.666015625, 1 ) };
		\addplot[ lightgray, line width = 0.05mm ] coordinates { ( 0.666259765625, -1 ) ( 0.666259765625, 1 ) };
		\addplot[ lightgray, line width = 0.05mm ] coordinates { ( 0.6663818359375, -1 ) ( 0.6663818359375, 1 ) };
		\addplot[ lightgray, line width = 0.05mm ] coordinates { ( 0.66650390625, -1 ) ( 0.66650390625, 1 ) };
		\addplot[ lightgray, line width = 0.05mm ] coordinates { ( 0.66656494140625, -1 ) ( 0.66656494140625, 1 ) };
		\addplot[ lightgray, line width = 0.05mm ] coordinates { ( 0.6666259765625, -1 ) ( 0.6666259765625, 1 ) };
		\addplot[ lightgray, line width = 0.05mm ] coordinates { ( 0.66668701171875, -1 ) ( 0.66668701171875, 1 ) };
		\addplot[ lightgray, line width = 0.05mm ] coordinates { ( 0.666748046875, -1 ) ( 0.666748046875, 1 ) };
		\addplot[ lightgray, line width = 0.05mm ] coordinates { ( 0.66680908203125, -1 ) ( 0.66680908203125, 1 ) };
		\addplot[ lightgray, line width = 0.05mm ] coordinates { ( 0.6668701171875, -1 ) ( 0.6668701171875, 1 ) };
		\addplot[ lightgray, line width = 0.05mm ] coordinates { ( 0.6669921875, -1 ) ( 0.6669921875, 1 ) };
		\addplot[ lightgray, line width = 0.05mm ] coordinates { ( 0.66796875, -1 ) ( 0.66796875, 1 ) };
		\addplot[ lightgray, line width = 0.05mm ] coordinates { ( 0.671875, -1 ) ( 0.671875, 1 ) };
		\addplot[ lightgray, line width = 0.05mm ] coordinates { ( 0.6875, -1 ) ( 0.6875, 1 ) };
		\addplot[ lightgray, line width = 0.05mm ] coordinates { ( 0.75, -1 ) ( 0.75, 1 ) };

		\addplot[ lightgray, line width = 0.05mm ] coordinates { (1, 0) (1, 0.11) };

		\addplot[ black, domain = 0:0.66666 ] { 0.66666 / ( 0.66666 * 100 - 90 * x ) };
		\addplot[ black, domain = 0.66666:1 ] { 0.1 };

		\addplot[ orange ] table [x = x, y = eps, col sep = comma] {solution_c0_material_benchmark_1d_kelly_k14R2.csv};

		\legend{
		};

	\end{axis}

	\spy[height = 2.0cm, width = 3cm] on (10.6, -0.5) in node[fill = white] at (8.35, -1.5);

\end{tikzpicture}
	\caption{Solution of the displacement \(u(x)\) and the strain \(\varepsilon(x)\) obtained from different meshes. The figures on the left show the solution of a single high-order element with polynomial degree \(p = 20\). The figures on the right show an iteratively refined mesh using the Kelly indicator. The refined overlay meshes are indicated by the gray vertical lines. The exact solution for \(u(x)\) and \(\varepsilon(x)\) is shown in black in each case.}
	\label{OneDSolution}
\end{figure}

\paragraph{Solution assessment} 
The solution of the problem is depicted in Figure~\ref{OneDSolution} in terms of the displacement \(u(x)\) and strain \(\varepsilon(x)\). The left side shows the solution obtained with a single high-order element of polynomial degree \(p=20\), whereas the right side corresponds to a discretization with a base mesh of \(p=15\) and \(k=14\) local refinement steps concentrated at the location \(a\) of the material transition. The polynomial degree is reduced by one with each successive overlay level.

The benchmark was selected to emulate the localized reduction in regularity that is commonly encountered at evolving elasto-plastic fronts. As expected, the globally smooth \(\mathcal{C}^{\infty}\)-continuous approximation is unable to accurately reproduce the \(\mathcal{C}^{0}\)-continuous strain field and exhibits pronounced oscillations in the vicinity of the material kink. In contrast, the locally refined solution accurately resolves the discontinuity and reproduces the exact solution without noticeable overshooting or oscillatory artifacts. The example demonstrates that the proposed refinement strategy is capable of efficiently resolving localized non-smooth solution features that are representative of those encountered in non-linear material analyses.

\paragraph{Convergence behavior}
To quantify the accuracy of the approximation, the relative error in the energy norm is evaluated for a sequence of successive refinement steps:
\begin{equationarray}{rclclclcl}
	\|e\|_E &=& \sqrt{\frac{|a(u_h,u_h)-a(u,u)|}{|a(u,u)|}} \cdot 100\% \;
\end{equationarray}
where \(a(u_h,u_h)\) denotes the computed energy norm of the approximate solution \(u_h\) and \(a(u,u)\) the corresponding value of the exact solution \(u(x)\).

As a baseline, we use a series of uniform \(p\)-refinements, shown in Fig. \ref{OneDConvergence}. The series is compared against two local refinement strategies: for the results shown in orange (filled dots), we applied one additional refinement layer in each cycle at the location \(a\) of the jump of the Young's modulus. For the second refinement strategy shown in green (filled triangles), we used the refinement indicator of Eq. \eqref{Kelly} and refined the element which yielded the value of largest magnitude. The convergence results of the local refinements confirm to a high degree the superiority gained compared to the global \(p\)-refinement. With comparable effort, we obtain a solution that is several orders of magnitude more accurate than in the non-refined case whereby for this simple test case both, manual and an error indicator based refinement, exhibit an identical behavior. 
\begin{figure}[h!]
	\centering
	\tikzsetnextfilename{convergence_c0_material_benchmark_1d}
\begin{tikzpicture}
	\begin{axis} [
		font={\footnotesize},
		axis lines = box,
		xlabel = Number of unconstrained DoFs $N$,
		ylabel = Relative error in the energy norm $||e||_E$ in \%,
		width = 0.75\textwidth,
		height = 9cm,
		max space between ticks = 80,
		minor x tick num = 0,
		minor y tick num = 0,
		cycle list name = color_list_mark,
		legend pos = south west,
		legend style = { nodes = { scale = 0.75, transform shape }, },
		legend cell align = { left },
		xmode = log,
		ymode = log,
		grid = major,
		grid style = {densely dashed, line width = 0.1pt},
		]

		\addplot table [x = n, y = e, col sep = comma] {error_c0_material_benchmark_1d_pfemR2.csv};
		\addplot table [x = n, y = e, col sep = comma] {error_c0_material_benchmark_1d_mlhp_manuelR2.csv};
		\addplot table [x = n, y = e, col sep = comma] {error_c0_material_benchmark_1d_mlhp_kellyR2.csv};

		\legend{
			Global \textit{p}-refinements,
			Manual local \textit{hp}-refinements,
			Automated local \textit{hp}-refinements,
		};

	\end{axis}
\end{tikzpicture}
	\caption{Convergence behavior of the total potential energy norm with respect to the number of unknown degrees of freedom for different refinement strategies.}
	\label{OneDConvergence}
\end{figure}

\subsection{Two-dimensional model with \(\mathcal{C}^0\)-continuous material}
\label{sec:twoDContinuousMaterial}
Next, we extend the one-dimensional benchmark problem to two dimensions to investigate the performance of the proposed approach along a material jump interface and along element boundaries in the relevant region. The solution domain has an outer radius \(b\) and a material interface at \(r=a\). The outer circumference of the domain is subjected to uniform stresses normal to the boundary.
\begin{figure}[h!]
	\centering
	\definecolor{body}{RGB}{220, 220, 220}
\tikzsetnextfilename{c0materialBenchmark2D}
\begin{tikzpicture}

	\draw[black, densely dotted] (4, 0) arc(0:-12:4);
	\draw[black, densely dotted] (0, 4) arc(90:102:4);
	\draw[black, fill = body] (0, 0) -- ++(4, 0) arc(0:90:4) -- cycle;
	\draw[black, densely dotted] (1.5, 0) arc(0:90:1.5);
	\draw[black, densely dotted] (0, 0) rectangle ++(2.5, 2.5);

	\draw[black] (0.15, -0.2) -- ++(3.7, 0);
	\draw[black] (-0.2, 0.15) -- ++(0, 3.7);
	\foreach \x in {0.25, 0.75, ..., 3.75}
	{
		\draw[black, fill = lightgray!20!white] (\x, -0.1) circle (0.1);
		\draw[black, fill = lightgray!20!white] (-0.1, \x) circle (0.1);
	}

	\foreach \phi in {-10, 0, ..., 100}
	{
		\draw[black, -latex] (\phi:4.1) -- (\phi:4.6);
	}

	\draw[black, -latex] (0, 0) -- node[below, rotate = 60, pos = 0.9] {$a$} (60:1.5);
	\draw[black, -latex] (0, 0) -- node[below, rotate = 80, pos = 0.95] {$b$} (80:4);
	\draw[black, latex-latex] (-0.6, 0) -- node[left] {$c$} ++(0, 2.5);
	\draw[black] (-0.7, 2.5) -- ++(0.4, 0);
	\draw[black] (-0.7, 0) -- ++(0.6, 0);

	\draw[black, -latex] (0, 0) -- node[below, rotate = 30, pos = 0.9] {$r$} (30:2.5);
	\draw[black, -latex] (2, 0) arc(0:30:2);
	\node[right, rotate = 10] at (10:2) {$\varphi$};

	\begin{axis} [
		shift = {(6cm, 0)},
		font={\footnotesize},
		axis lines = box,
		xlabel = $r$,
		title = {Displacement},
		xticklabels = {0, $\frac{2}{3}$, $2$},
		xtick = {0, 0.666666, 2},
		width = 5cm,
		height = 6cm,
		max space between ticks = 80,
		minor x tick num = 0,
		minor y tick num = 0,
		cycle list name = color_list,
		legend pos = south east,
		legend style = { nodes = { scale = 0.75, transform shape }, },
		legend cell align = { left },
		grid = major,
		grid style = {densely dashed, line width = 0.1pt}
	]

		\addplot[blue] table [x = x, y = u, col sep = comma] {solution_c0_material_benchmark_2d_analytic_u_innerR2.csv};
		\addplot[blue] table [x = x, y = u, col sep = comma] {solution_c0_material_benchmark_2d_analytic_u_outerR2.csv};

		\legend{
			$u_r$
		};

	\end{axis}

	\begin{axis} [
		shift = {(11cm, 0)},
		font={\footnotesize},
		axis lines = box,
		xlabel = $r$,
		title = {Strains},
		xticklabels = {0, $\frac{2}{3}$, $2$},
		xtick = {0, 0.666666, 2},
		width = 5cm,
		height = 6cm,
		max space between ticks = 80,
		minor x tick num = 0,
		minor y tick num = 0,
		cycle list name = color_list,
		legend pos = south east,
		legend style = { nodes = { scale = 0.75, transform shape }, },
		legend cell align = { left },
		grid = major,
		grid style = {densely dashed, line width = 0.1pt}
	]

		\addplot[blue] table [x = x, y = err, col sep = comma] {solution_c0_material_benchmark_2d_analytic_err_innerR2.csv};
		\addplot[orange] table [x = x, y = epp, col sep = comma] {solution_c0_material_benchmark_2d_analytic_epp_innerR2.csv};
		\addplot[blue] table [x = x, y = err, col sep = comma] {solution_c0_material_benchmark_2d_analytic_err_outerR2.csv};
		\addplot[orange] table [x = x, y = epp, col sep = comma] {solution_c0_material_benchmark_2d_analytic_epp_outerR2.csv};

		\legend{
			$\varepsilon_{rr}$,
			$\varepsilon_{\varphi\varphi}$
		};

	\end{axis}

\end{tikzpicture} 
	\caption{Problem setup of a circular plate with \(\mathcal{C}^0\)-continuous material (left). Due to symmetry of the problem, only one quarter of the total domain is considered. The exact solution of the displacement field \(u_r(r)\) (mid) and the state of strain \(\boldsymbol{\varepsilon}(r)\) (right) is shown along the radial coordinate \(r\).}
	\label{TwoDBenchmarkSolution}
\end{figure}

The modulus of elasticity has a linearly variable value for \(r \in [0,a]\) and a constant value in the remaining domain (cf. Fig. \ref{TwoDBenchmarkSolution}):
\begin{equationarray}{rcl}
	E(r) &=& 
	\begin{cases}
		E_0+\frac{\Delta E}{a}(a-r) \quad \ & r\in [0,a]\\
		E_0 & r\in (a,b]\;.
	\end{cases} 
\end{equationarray}

\noindent With a Poisson ratio \(\nu=0, \;\forall \mathbf{x}\in \Omega\) the problem reduces to the linear ordinary differential equation
\begin{equationarray}{rclclclclc}
	\frac{\partial}{\partial r}\biggl[E \frac{\partial u_r}{\partial r}\biggr]+\frac{E}{r}\biggl[\frac{\partial u_r}{\partial r}-\frac{u_r}{r}\biggr] &=& 0
\end{equationarray}
with boundary conditions \(u_r(0)\,=\,0\) and \(\sigma_{rr}\bigl|_{r=b}\,=\,E_0 (C_1 - C_2 / b^2)\). The radial stress is expressed in terms of \(C_1 = (a \,t_0 + E_0 \,\hat{u}_r(a)) / (2 a \,E_0)\) and \(C_2 = (-t_0 \,a^3 + a^2 \,E_0 \,\hat{u}_r(a)) / (2 a\, E_0)\) where \(t_0\) denotes the radial stress at the material interface \(r = a\) and \(\hat{u}_r\) the displacement field in interval \([0,a]\). The analytical solution of the problem and the expression for \(\hat{u}_r\) are provided in Appendix B.

Exploiting symmetry, the computational domain is restricted to the square of side length \(c\) shown in Fig.~\ref{TwoDBenchmarkSolution}. Symmetry conditions are imposed along the two symmetry boundaries, while the exact displacement
solution is prescribed along the remaining two edges. 
For all following computations, the geometry and material parameters are chosen as \(a=\tfrac{2}{3}\), \(b=3\), \(c=2\), \(E_0=10\), \(\Delta E = 90\), \(\nu=0\) and \(t_0=1\). With the aforementioned choice of parameter values the analytic value of the total potential energy evaluates to \(a(\mathbf{u},\mathbf{u}) \approx  0.3398437370278343 \).

\begin{figure}[h!]
	\centering
	\input{solution_c0_material_benchmark_2dR3.tex}
	\caption{Pointwise displacement gradient norm for different refinement models: uniform and unrefined model with \(p=7\) (left), manually refined with \(k=6\) overlay meshes (mid), automatically refined with \(k=6\) overlay meshes (right).}
	\label{TwoDBenchmarkSolutionApprox}
\end{figure}

\paragraph{Solution assessment} 
The magnitude of the displacement gradient \(\|\nabla\mathbf{u}\|_2\) obtained from three different refinement strategies is depicted in Fig.~\ref{TwoDBenchmarkSolutionApprox}. A \((5\times 5)\) base mesh with polynomial degree \(p=7\) covering the square domain \([0,c]\times[0,c]\) was employed for all computations. Similar to the one-dimensional benchmark, the solution obtained by uniform \(p\)-refinement exhibits pronounced oscillations in the vicinity of the material interface \(r=a\), indicating insufficient local resolution. A substantial improvement is achieved by manually refining all elements intersected by the interface in six successive refinement cycles, resulting in an accurate representation of the material kink. For the automatically refined model, the error indicator of Eq.~\eqref{Kelly} was evaluated after each solution step and the \(30\%\) of elements with the largest indicator values were refined. In addition, up to \(20\%\) of the elements with the smallest indicator values were considered for de-refinement, provided that all \(2^d\) child elements of the corresponding parent element satisfied the criterion. The resulting discretization concentrates refinement in the vicinity of the material interface and achieves a solution quality comparable to the manually refined model, while requiring no prior knowledge of the location of the material transition.

\begin{figure}[h!]
	\centering
	\tikzsetnextfilename{c0materialBenchmark2DConvergence}
\begin{tikzpicture}
	\begin{axis} [
		font={\footnotesize},
		axis lines = box,
		xlabel = Number of unconstrained DoFs $N$,
		ylabel = Relative error in the energy norm $||e||_E$ in \%,
		width = 0.75\textwidth,
		height = 9cm,
		max space between ticks = 80,
		minor x tick num = 0,
		minor y tick num = 0,
		cycle list name = color_list_mark,
		legend pos = south west,
		legend style = { nodes = { scale = 0.75, transform shape }, },
		legend cell align = { left },
		xmode = log,
		ymode = log,
		grid = major,
		grid style = {densely dashed, line width = 0.1pt}
		]

		\addplot table [x = n, y = e, col sep = comma] {error_c0_material_benchmark_2d_pfemR2.csv};
		\addplot table [x = n, y = e, col sep = comma] {error_c0_material_benchmark_2d_mlhp_manuelR2.csv};
		\addplot table [x = n, y = e, col sep = comma] {error_c0_material_benchmark_2d_mlhp_kellyR2.csv};

		\legend{
			Global \textit{p}-refinements,
			Manual local \textit{hp}-refinements,
			Automated local \textit{hp}-refinements,
		};

	\end{axis}
\end{tikzpicture}
	\caption{Convergence behavior of the total potential energy norm with respect to the number of unknown degrees of freedom for different refinement strategies.}
	\label{TwoDBenchmarkConvergence}
\end{figure}

\paragraph{Convergence behavior} A comparison of the convergence rates of all three models is depicted in Fig. \ref{TwoDBenchmarkConvergence}. 
The results are comparable to the one-dimensional problem shown before. The global \(p\)-refinements converged only with a low rate and exhibited strong oscillations. With manual local refinements, on the other hand, the model converged noticeably faster. However, after eight cycles of refinement the rate stagnates and the curve levels off with a final error of \(3.75 \cdot 10^{-3}\%\). 
This loss of convergence can largely be attributed to the empirically motivated refinement strategy, whose sub-cell distribution is not optimally aligned with the local solution characteristics. 
In contrast, the automated refinements using the error indicator converged at a similar error level with only slightly more degrees-of-freedom and kept the overall convergence rate for continuing refinement. 

\section{Thermo-viscoplasticity}
\label{sec:Thermo-viscoplasticity}
The focus of this publication is on multi-level \(hp\)-FCM analyses of thermo-viscoplastic problems of complex analysis domains. In the following, we briefly present the implemented material model including the fundamental relations of elasto-viscoplasticity and linear thermo-elasticity, followed by the specification of the model parameters which are characteristic to rate-independent and rate-dependent plasticity. For an in-depth treatment of the thermo-viscoplastic material model, we refer the interested reader to relevant literature \cite{Oppermann:2022,deSouzaNeto:2011}.

\subsection{Model specification}
\label{ModelSpecification}
We consider a solid body which occupies a bounded Lipschitz domain \(\Omega_{}\subset\mathbb{R}^d, d\leq 3\) and time domain \(\mathcal{T}\subset\mathbb{R}_{+}\) with the following state variables: 
(i) the temperature field \(T(\mathbf{x},t): \Omega_{}\times \mathcal{T}\rightarrow \mathbb{R}\) and (ii) the displacement field \(\mathbf{u}(\mathbf{x},t): \Omega_{}\times \mathcal{T}\rightarrow \mathbb{R}^3\) from which the small strain tensor is derived as \(\boldsymbol{\varepsilon}=\tfrac{1}{2}\left(\nabla\mathbf{u}+(\nabla\mathbf{u})^T\right)\). The total strain tensor is split into an elastic, a viscoplastic and a thermal strain contribution, according to the parts that contribute to the Helmholtz free energy potential, which forms the basis for the formulation of a constitutive model consistent with thermodynamic principles:
\begin{equationarray}{lclcl}
	\boldsymbol{\varepsilon} &=& \boldsymbol{\varepsilon}^{e}+\boldsymbol{\varepsilon}^{vp}+\boldsymbol{\varepsilon}^{th}  \,.
\end{equationarray}
\noindent The constitutive equations are governed by the general relation:
\begin{equationarray}{lclcl}
	\label{consti}\boldsymbol{\sigma} &=& \mathbf{C}\,:\,\left(\boldsymbol{\varepsilon}-\boldsymbol{\varepsilon}^{vp}-\boldsymbol{\varepsilon}^{th}\right)
\end{equationarray}
with the Cauchy stress tensor \(\boldsymbol{\sigma}\) and a fourth order elastic material tensor \(\mathbf{C}\).

\subsection{Linear thermo-elasticity}\label{linearThermoElasticity}
The thermal part of Eq. \eqref{consti} depends on the change in temperature assuming isotropic and homogeneous thermal expansion 
\begin{equationarray}{lclcl}
	\label{epsth}\boldsymbol{\varepsilon}^{th} &=& \gamma\,(T-T_0)\,\mathbf{I}
\end{equationarray}
\noindent with the coefficient of thermal expansion \(\gamma=\gamma(T)\), where \(T\) denotes the temperature at $\mathbf{x}$, \(T_0\) is a reference temperature and \(\mathbf{I}\) is the second-order identity tensor, providing a purely volumetric strain contribution. The thermal expansion is temperature-dependent and considered as a piecewise linear fit to the experimental data from \cite{Spittel:2009.2}, cf. Figs. \ref{ConstitutivePropertiesI} and \ref{ConstitutivePropertiesII}, Appendix C. 
Moreover, a deformation-induced temperature change is neglected in this model due to the assumption of small strains, which leads to a one-sided coupling of the thermal and elastic problem and enables a partitioned solution scheme as shown in \cite{Zander:12.1}.

In the framework of viscoplasticity, the governing system of equations is solved in incremental steps with a global Newton-Raphson scheme to enforce non-linear equilibrium. The partitioned solution approach allows for an initial temperature solution in each incremental step that is independent of the viscoplastic solution scheme. The thermal strain \(\boldsymbol{\varepsilon}^{th}\) is then directly used in Eq. \eqref{consti} in the local scope of the material equations. The temperatures of Eq. \eqref{epsth} result from the solution of the stationary heat flow problem in its weak form as introduced in Eq. \eqref{WFHeat}.

While the assumption of a one-way coupling between the thermal and mechanical fields is appropriate for problems characterized by moderate deformation rates and small strains, it neglects heat generation due to plastic dissipation and thermo-mechanical feedback effects. In processes involving high strain rates or significant plastic work, this simplification may lead to an underestimation of local temperature increases and, consequently, affect the predicted material response. Nevertheless, the adopted formulation provides a computationally efficient framework and represents a reasonable approximation for a wide range of quasi-static and moderately coupled thermo-mechanical problems. The extension to fully coupled thermo-viscoplastic formulations, including deformation-induced heat generation, constitutes a potential direction for future work.

\subsection{Viscoplasticity model}
Our model uses the classical $J_2$ flow theory applying the \emph{von-Mises} yield function
\begin{equationarray}{rclcl}
	\Phi(\boldsymbol{\sigma},\bar{\varepsilon}^p) &=& \sigma_v(\boldsymbol{\sigma})-\sigma_y(\bar{\varepsilon}^p,T) &&
	\begin{cases}
		\;<0 & \mbox{elastic deformation}\\
		\;\geq 0 & \mbox{viscoplastic yielding}
	\end{cases}\\[6mm]
	\sigma_v(\boldsymbol{\sigma}) &=& \sqrt{\frac{3}{2}}\,\|\boldsymbol{\sigma}^{dev}\| &&\mbox{with } \boldsymbol{\sigma}^{dev} \,=\, \mathbb{P}:\boldsymbol{\sigma}
\end{equationarray}
where \(\sigma_v(\boldsymbol{\sigma})\) denotes the \emph{von-Mises} equivalent stress considering the deviatoric stresses \(\boldsymbol{\sigma}^{dev}\) obtained from projection with the fourth order deviatoric projection tensor \(\mathbb{P}\) \cite{Simo:1998}. The yield stress can harden for each temperature from the initial to the maximum saturated level, cf. Fig. \ref{ConstitutivePropertiesI}, Appendix C.

\noindent The model describes hardening by the evolution of the yield stress 
\begin{equationarray}{rcl}
	\label{strainhardeningA}\sigma_y(\bar{\varepsilon}^p,T) &=& \sigma_{y_0}(T)-\Phi_{\chi}(T)\,\chi(\bar{\varepsilon}^p)
\end{equationarray}
starting with the initial temperature-dependent yield stress \(\sigma_{y_0}(T)\) and evolving with a contribution which combines strain hardening \(\chi(\bar{\varepsilon}^p)\) and thermal softening \(\Phi_{\chi}(T)\).
The latter is defined by the arctangent ansatz
\begin{equationarray}{rclclclcl}
	\Phi_\chi(T) &=& 1-\frac{2\omega_\chi}{\pi-2\varphi_\chi}
\left[\tan^{-1}\!\left(\chi_\chi (T-T_r)+\tan(\varphi_\chi)\right)-\varphi_\chi\right].
\end{equationarray}
The corresponding material parameters are listed in Table~\ref{vpParameters} in Appendix D. Further details of the constitutive model and the associated parameter identification can be found in \cite{Oppermann:2022}. 

\noindent The strain hardening law is given by the exponential relation
\begin{equationarray}{rcl}
	\label{strainhardeningB}\chi(\bar{\varepsilon}^p) &=& \Delta\sigma_{y_{\infty 0}}\left[\exp\biggl(-\frac{\bar{\varepsilon}^p}{\bar{\varepsilon}_0^p}\biggr)-1 \right]
\end{equationarray}
where \(\Delta\sigma_{y_{\infty 0}}\) is the isothermal saturation yield stress increment denoting the increase in yield stress due to plastic deformation and \(\bar{\varepsilon}_0^p\) denotes the coefficient of saturation as limit point of appreciable hardening. 
The hardening modulus is obtained by the slope of the yield stress as:
\begin{equationarray}{rclcl}
	H(\bar{\varepsilon}^p) &=& \frac{\partial \sigma_y}{\partial \bar{\varepsilon}^p} &=& \Phi_{\chi}(T) \frac{\Delta\sigma_{y_{\infty 0}}}{\bar{\varepsilon}_0^p} \exp\biggl(-\frac{\bar{\varepsilon}^p}{\bar{\varepsilon}_0^p} \biggr)\;.
\end{equationarray}

\noindent The rate-dependent plasticity model used here to describe the viscous reaction follows the \emph{Perzyna}-type constitutive law \cite{Perzyna:1971}:
\begin{equationarray}{rcl}
	\dot{\gamma}(\boldsymbol{\sigma},\bar{\varepsilon}^p) &=& 
	\begin{cases}
		\frac{1}{\mu}\left[\frac{\sigma_v(\boldsymbol{\sigma})}{\sigma_y(\bar{\varepsilon}^p)}-1 \right]^{\frac{1}{m}}   & \mbox{if }\;\; \Phi(\boldsymbol{\sigma},\bar{\varepsilon}^p)\geq 0\\
		0 & \mbox{if }\;\; \Phi(\boldsymbol{\sigma},\bar{\varepsilon}^p)< 0 \hfill \\
	\end{cases}
\end{equationarray}
where \(\mu>0\) is the relaxation time in seconds representing the viscosity and \(m>0\) is the
non-dimensional viscoplastic exponent representing the rate sensitivity for the considered material.

The material parameters employed in the present study are taken from literature and experimental data for representative steels and are not subject to calibration within this work. In general, the parameters of the viscoplastic model such as the initial yield stress, hardening modulus, saturation parameters, and rate sensitivity, are identified from standard mechanical tests, including uniaxial tension and creep experiments at different temperatures. The focus of this contribution lies on the numerical treatment of the governing equations rather than on parameter identification.

The Perzyna-type flow rule introduced above is employed for the thermo-viscoplastic analysis of the porous plate in Section \ref{sec:porousPlate}. For the benchmark problem considered in Section \ref{sec:plateWH}, the modified flow rule proposed by Perić is adopted to ensure consistency with the reference solution and is introduced there.

\subsection{Numerical solution framework}
Both the thermal and the mechanical sub-problems of the thermo-viscoplastic analysis are solved in the framework of the multi-level \(hp\)-Finite Cell Method. To this end, the analysis domain is discretized with finite cells on a Cartesian grid and locally refined by hierarchical overlay meshes as needed, cf. sections \ref{sec:ImmersedLocalRefinement} and \ref{sec:refinement}. The solution of the thermal problem follows the stationary heat problem described in Section \ref{sec:FCM}.  

The governing equations of the viscoplastic problem are solved incrementally in time using an implicit backward Euler scheme. With the displacement field at time \(t_{n+1}\) expressed as
\begin{equationarray}{rcl}
    \mathbf{u}_{(n+1)} &=&  \mathbf{u}_{(n)} + \Delta\mathbf{u}_{(n+1)}  \qquad  \text{in } [t_n,t_{n+1}],
\end{equationarray}
the global equilibrium condition at \(t_{n+1}\) is given by
\begin{equationarray}{rclcl}
    \mathbf{R}_{(n+1)} &=& \mathbf{f}^{\mathrm{ext}}_{(n+1)}-\mathbf{f}^{\mathrm{int}}_{(n+1)} &=& \mathbf{0},
\label{eq:equilibrium}
\end{equationarray}
where the internal and external force vectors are defined as
\begin{equationarray}{rclclc}
    \mathbf{f}^{\mathrm{int}}_{(n+1)} &=&  \mathop{\A}\limits_e
    \int_{\Omega_e} \boldsymbol{\varepsilon}(\mathbf{v}) : \alpha\, \boldsymbol{\sigma}_{(n+1)} \,\mathrm{d}\Omega
\label{eq:internal_force} \qquad\mathrm{ and}\\[3mm]

    \mathbf{f}^{\mathrm{ext}}_{(n+1)} &=& \mathop{\A}\limits_e 
    \int_{\Omega_e} \mathbf{v}\cdot   \alpha\,\bar{\mathbf{p}}_{(n+1)} \,\mathrm{d}\Omega  +
    \mathop{\A}\limits_e  \int_{\Gamma_{N,e}}  \mathbf{v}\cdot \alpha\,\bar{\mathbf{t}}_{(n+1)}
    \,\mathrm{d}\Gamma .
\label{eq:external_force}
\end{equationarray}
and $\mathbf{v}$ denotes an admissible test function.
Here, \(\bar{\mathbf{p}}\) and \(\bar{\mathbf{t}}\) denote prescribed body forces and boundary tractions, respectively. Due to the non-linear constitutive response, the equilibrium equations are solved by a Newton--Raphson iteration. At iteration \(i\), linearization about the current displacement iterate \(\mathbf{u}_{(n+1)}^{(i)}\) yields
\begin{equationarray}{rclcl}
    \mathop{\A}\limits_e  \int_{\Omega_e}  \boldsymbol{\varepsilon}(\mathbf{v})  :
    \alpha\,   \mathbb{C}_{(n+1)}^{\mathrm{vp},(i)}
    :  \boldsymbol{\varepsilon}  \left(  \delta\mathbf{u}_{(n+1)}^{(i+1)} \right)
    \,\mathrm{d}\Omega  &=&  \mathbf{R}_{(n+1)}^{(i)},
\label{eq:newton}
\end{equationarray}
with
\begin{equationarray}{rclcl}
    \mathbf{R}_{(n+1)}^{(i)}  &=&
    \mathbf{f}^{\mathrm{ext}}_{(n+1)}  - \mathbf{f}^{\mathrm{int},(i)}_{(n+1)},
\qquad  \mathbf{f}^{\mathrm{int},(i)}_{(n+1)} &=&
    \mathop{\A}\limits_e   \int_{\Omega_e}   \boldsymbol{\varepsilon}(\mathbf{v})
    :   \alpha\,  \boldsymbol{\sigma}_{(n+1)}^{(i)}    \,\mathrm{d}\Omega 
\label{eq:residual}
\end{equationarray}
where \(\A\) denotes the assembly operator of the discretization and \(\alpha\) the indicator function introduced in Section~\ref{sec:FCM}.
The displacement field is updated according to \(\mathbf{u}_{(n+1)}^{(i+1)} = \mathbf{u}_{(n+1)}^{(i)} + \delta\mathbf{u}_{(n+1)}^{(i+1)} \) until the prescribed convergence criterion is satisfied. Here, \(\delta\mathbf{u}_{(n+1)}^{(i+1)}\) denotes the Newton correction, whereas \(\Delta\mathbf{u}_{(n+1)}=\mathbf{u}_{(n+1)}-\mathbf{u}_{(n)}\) denotes the total displacement increment over the time step. 

The algorithmic consistent viscoplastic tangent operator is defined as
\begin{equationarray}{rcl}
	\mathbb{C}_{(n+1)}^{vp,(i)}&=& \frac{\partial\boldsymbol{\sigma}_{(n+1)}^{(i)}}{\partial\boldsymbol{\varepsilon}_{(n+1)}^{(i)}}
\end{equationarray}
and recomputed in each iteration cycle $i$ using a classical \emph{elastic predictor-return-mapping} algorithm \cite{Simo:1998,deSouzaNeto:2011}.

\noindent Essential boundary conditions are prescribed as
\begin{equationarray}{rclclcl}
    u_j=\bar{u}_j,
    \qquad
    j\in\{1,2,3\},
    \qquad
    \forall\,\mathbf{x}\in\Gamma_D ,
\end{equationarray}
where \(\Gamma_D\) and \(\Gamma_N\) denote the Dirichlet and Neumann boundaries, respectively, with
\(\partial\Omega=\Gamma_D\cup\Gamma_N\) and \(\Gamma_D\cap\Gamma_N=\emptyset\).

 In this context, it should be noted that the calculation results are very sensitive to the full consideration of the strain history data stored in the quadrature points and used to calculate the governing equations. The overall effort for storage, updating and usage of the strain history data drops dramatically by the use of the proposed \emph{Non-Negative Moment Fitting} compared to the \emph{Adaptive Space-Tree} quadrature method as demonstrated with the examples of the following section.
A detailed analysis of the convergence behavior of the non-linear solution procedure is deferred to the numerical examples section, where iteration counts and robustness are evaluated for representative problems.

\section{Numerical examples}
\label{sec:examples}

In the following section, the proposed formulation is assessed by means of two numerical examples. The classical plate-with-a-hole benchmark introduced by \emph{Simo} and \emph{Taylor} \cite{Simo:1986} is first employed to verify the formulation, assess the stability of the return-mapping algorithm in the presence of plasticity initiation within the fictitious domain, and subsequently extend the analysis to viscoplastic material behavior. A porous plate is then considered to demonstrate the capabilities of the proposed error-indicator-driven hierarchical refinement strategy for thermo-viscoplastic analyses of geometrically complex domains.

\subsection{Plate with circular hole}
\label{sec:plateWH}

\begin{figure}[hb!!]
	\centering
		\begin{tikzpicture}
			\node (A) {\tikzsetnextfilename{model_perforated_plate}
\begin{tikzpicture}[isometric view, scale = 0.75]

	\fill[red, opacity = 0.2] (0, 9, 0) -- ++(5, 0, 0) -- ++(0, 0, 0.5) -- ++(-5, 0, 0) -- cycle;
	\draw[red, thick] (0, 9, 0) -- ++(5, 0, 0) -- ++(0, 0, 0.5) -- ++(-5, 0, 0) -- cycle;

	\node[transform shape, above] at (2.5, 9, 0.6) {$\overline{u}$};

	\foreach \x in {0, 1, ..., 5}
		{
			\foreach \z in {0, 0.5}
				{
					\draw[red, thick, -stealth] (\x, 8, \z) -- ++(0, 1, 0);
				}
		}

	\draw[black] (5, 8, 0) -- ++(0, 0, 0.5);
	\begin{scope}[canvas is xy plane at z = 0]
		\fill[blue, opacity = 0.2] (2.5, 0) -- ++(2.5, 0) -- ++(0, 8) -- ++(-5, 0) -- ++(0, -5.5) arc(90:0:2.5) -- cycle;
		\draw[blue, pattern color = blue, pattern = north west lines] (2.5, 0) -- ++(2.5, 0) -- ++(0, 8) -- ++(-5, 0) -- ++(0, -5.5) arc(90:0:2.5) -- cycle;
	\end{scope}

	\begin{scope}[canvas is xy plane at z = 0.25]
		\draw[
			black,
			top color = lightgray,
			bottom color = lightgray,
			middle color = lightgray!0!white,
			shading angle = 15
		] (2.5, 0) -- ++(2.5, 0) -- ++(0, 8) -- ++(-5, 0) -- ++(0, -5.5) arc(90:0:2.5) -- cycle;
	\end{scope}

	\begin{scope}[canvas is xy plane at z = 0.5]
		\draw[black] (2.5, 0) -- ++(2.5, 0) -- ++(0, 8) -- ++(-5, 0) -- ++(0, -5.5) arc(90:0:2.5) -- cycle;
    \node[transform shape] at (2.5, 4.5) {
      \begin{minipage}{5cm}
        \Large
        \begin{equation*}
          \begin{aligned}
            h &= 36\,\mathrm{mm} \\
            w &= 20\,\mathrm{mm} \\
            t &= 1\,\mathrm{mm} \\
            R &= 5\,\mathrm{mm}
          \end{aligned}
        \end{equation*}
      \end{minipage}
    };
	\end{scope}

	\draw[black] (2.5, 0, 0) -- ++(0, 0, 0.5);
	\draw[black] (0, 2.5, 0) -- ++(0, 0, 0.5);
	\draw[black] (5, 0, 0) -- ++(0, 0, 0.5);
	\draw[black] (0, 8, 0) -- ++(0, 0, 0.5);

	\draw[black, densely dashed] (0, 0, 0) -- ++(2.5, 0, 0);
	\draw[black, densely dashed] (0, 0, 0) -- ++(0, 2.5, 0);
	\draw[black, densely dashed] (0, 0, 0.5) -- ++(2.5, 0, 0);
	\draw[black, densely dashed] (0, 0, 0.5) -- ++(0, 2.5, 0);

	\begin{scope}[canvas is xz plane at y = 0]
		\fill[blue, opacity = 0.2] (2.5, 0) rectangle ++(2.5, 0.5);
		\draw[blue, pattern color = blue, pattern = north west lines] (2.5, 0) rectangle ++(2.5, 0.5);
	\end{scope}

	\begin{scope}[canvas is yz plane at x = 0]
		\fill[blue, opacity = 0.2] (2.5, 0) rectangle ++(5.5, 0.5);
		\draw[blue, pattern color = blue, pattern = north east lines] (2.5, 0) rectangle ++(5.5, 0.5);
	\end{scope}

	\draw[black, -stealth] (0, 0, 0) -- node[transform shape, pos = 1.0, below right] {$x_1$} ++(1.5, 0, 0);
	\draw[black, -stealth] (0, 0, 0) -- node[transform shape, pos = 1.0, below left] {$x_2$} ++(0, 1.5, 0);
	\draw[black, -stealth] (0, 0, 0) -- node[transform shape, pos = 1.0, below left] {$x_3$} ++(0, 0, 1.5);

	\draw[black, stealth-stealth] (0, -1.5, 0) -- node[transform shape, below right] {$w/2$} ++(5, 0, 0);
	\draw[black, stealth-stealth] (-1.5, 0, 0) -- node[transform shape, below left] {$h/2$} ++(0, 8, 0);
	\draw[black, stealth-stealth] (5, -1.5, 0) -- node[transform shape, right, pos = 0.2] {$t$} ++(0, 0, 0.5);

	\draw[black] (0, -0.5, 0) -- ++(0, -1.5, 0);
	\draw[black] (5, -0.5, 0) -- ++(0, -1.5, 0);
	\draw[black] (5, -0.5, 0.5) -- ++(0, -1.5, 0);

	\draw[black] (-0.5, 0, 0) -- ++(-1.5, 0, 0);
	\draw[black] (-0.5, 8, 0) -- ++(-1.5, 0, 0);

	\draw[black, -stealth] (0, 0, 0) -- node[transform shape, above left, shift = {(0.2, 0, 0)}] {$R$} ++({2.5 * cos(20)}, {2.5 * sin(20)}, 0);

\end{tikzpicture}};
			\node[right, shift = {(10mm, 0mm)}] at (A.east) {
				\begin{minipage}{5cm}
					\footnotesize%
					\textbf{Material properties (\emph{von-Mises})} \\[-2mm]
					\begin{equation*}
						\begin{aligned}
							E &= 70\,\mathrm{GPa} \\
							\nu &= 0.2 \\
						\end{aligned}
					\end{equation*}
					\indent linear hardening \\[-2mm]
					\begin{equation*}
						\sigma_y(\overline{\varepsilon}^p) = 0.243 + 0.2 \overline{\varepsilon}^p\,\mathrm{[GPa]}
					\end{equation*}

					\vspace*{2mm}
					\textbf{Boundary conditions} \\
					symmetry boundary conditions \\[-2mm]
					\begin{equation*}
						\begin{aligned}
							u_1 & = 0 && \mathrm{along}\ (0,x_2,x_3) \\
							u_2 & = 0 && \mathrm{along}\ (x_1,0,x_3) \\
							u_3 & = 0 && \mathrm{along}\ (x_1,x_2,0) \\
						\end{aligned}
					\end{equation*}
					prescribed displacement \\[-2mm]
					\begin{equation*}
						u_2  = \overline{u} \quad \mathrm{along}\ (x_1,h/2,x_3) \\
					\end{equation*}
					
				\end{minipage}
			};
		\end{tikzpicture}
		\caption{Plate with a circular hole benchmark \cite{Simo:1986,deSouzaNeto:2011}: geometry, material properties and boundary conditions. }
		\label{PerforatedPlate}
\end{figure}

The geometry, boundary conditions and the applied material model of the first example are depicted in Fig. \ref{PerforatedPlate}. Due to the symmetry of the problem only a quarter of the total model was considered in the analysis. The linear hardening model applied an initial yield stress \(\sigma_0\,=\, 243\,\mathrm{MPa}\), a hardening modulus \(H\,=\,200\,\mathrm{MPa}\) and the remaining material parameters correspond to aluminum. A total loading \(\bar{u}\,=\,0.6\,\mathrm{mm}\) was applied in $43$ equal displacement increments.

\begin{figure}[htb!!]
	\centering
		\begin{tikzpicture}
			\node (m3) at (8.0cm, 0) {\includegraphics[width = 0.45\textwidth]{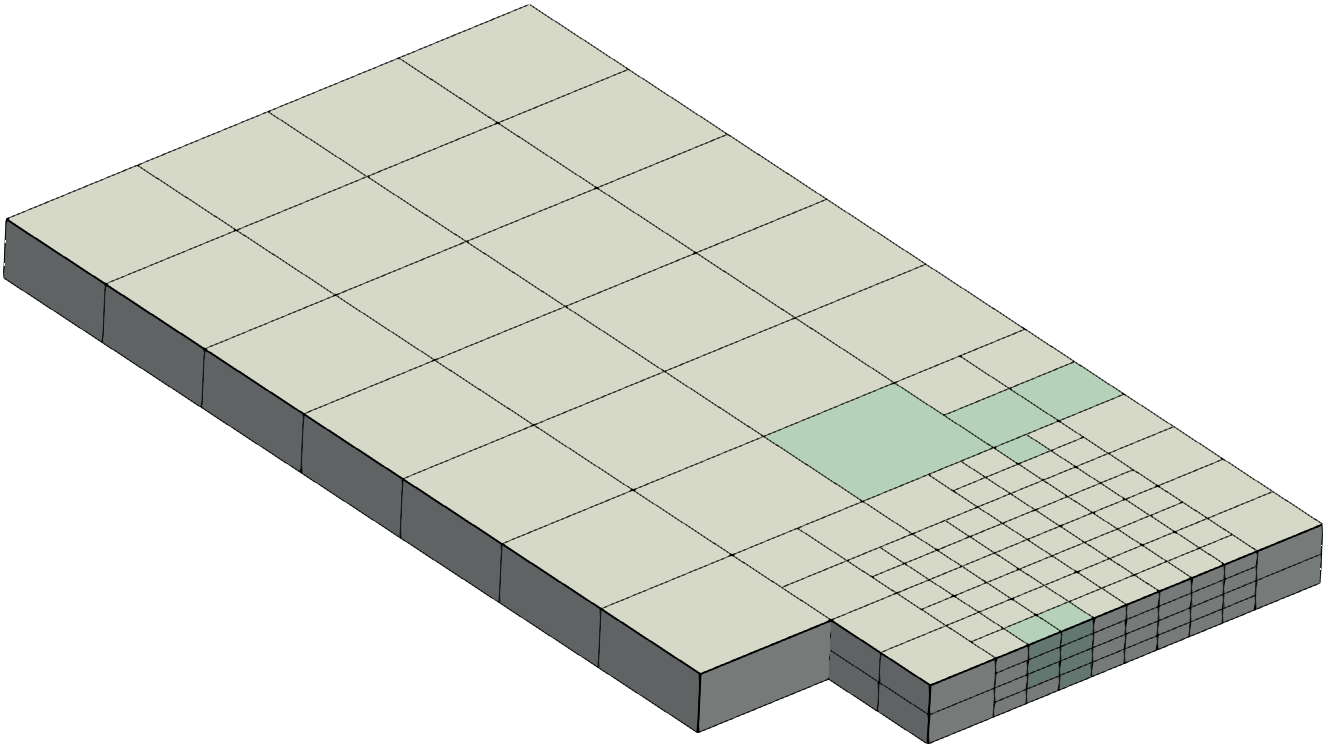}};
			\node (m2) at (4.0cm, 0) {\includegraphics[width = 0.45\textwidth]{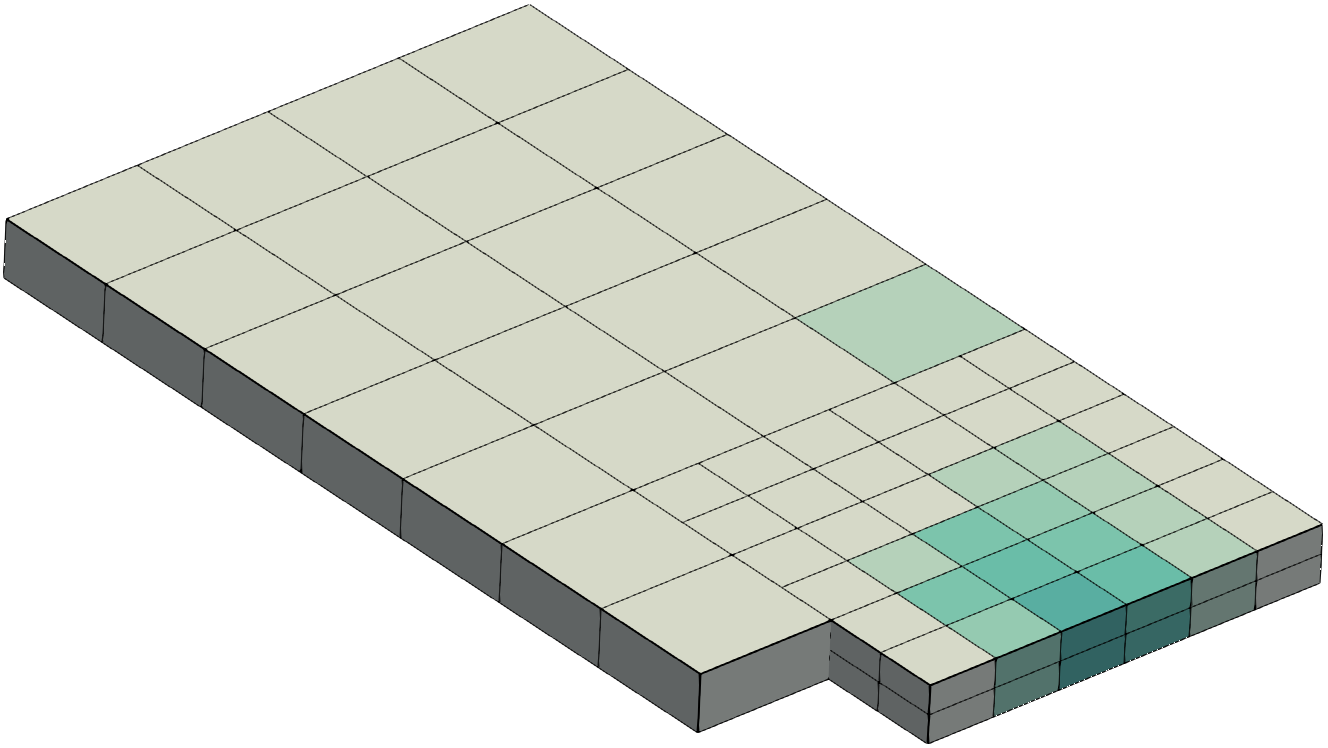}};
			\node (m1) at (0.0cm, 0) {\includegraphics[width = 0.45\textwidth]{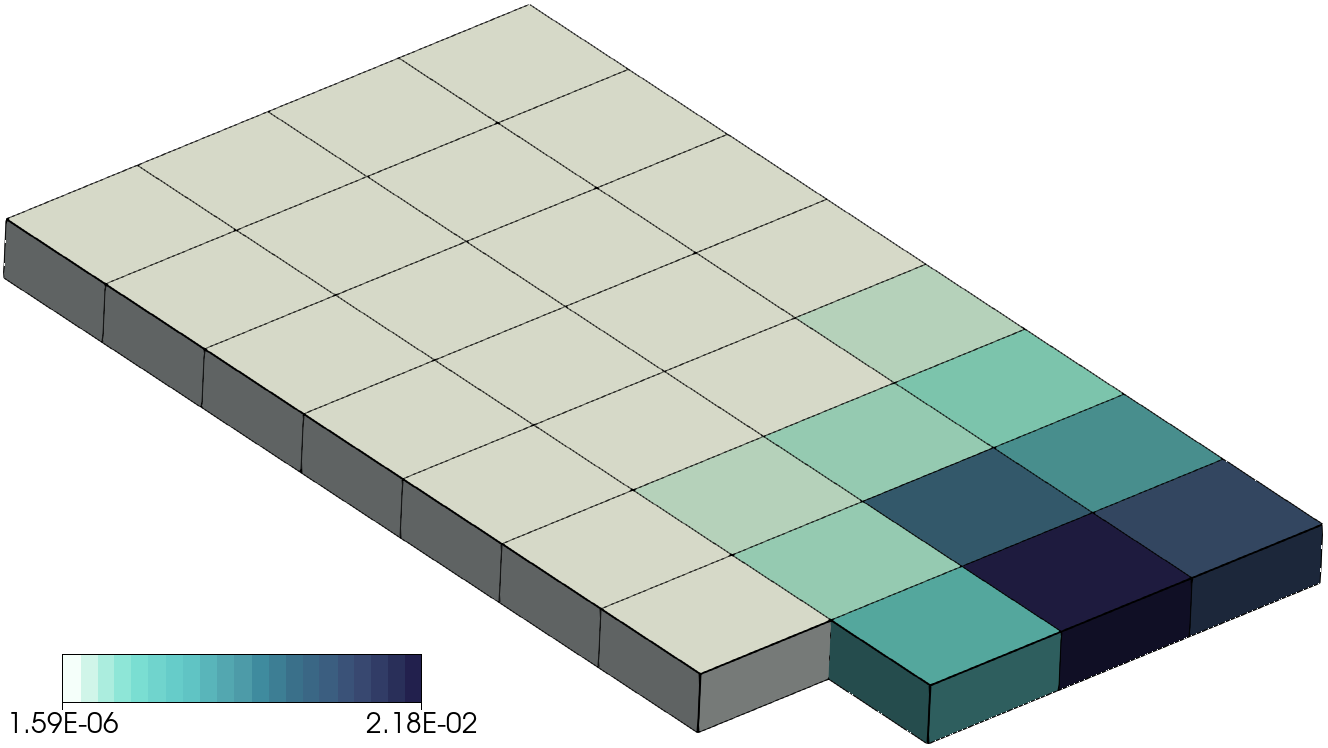}};

			\node[below] at (m1.south) {\footnotesize FCM base mesh $4 \times 8 \times 1$};
			\node[below] at (m2.south) {\footnotesize refinement $k = 1$};
			\node[below] at (m3.south) {\footnotesize refinement $k = 2$};
		\end{tikzpicture}
		\caption{Plate with a circular hole: gradient jump used in the error indicator of Eq. \eqref{Kelly} for the base model and two successive overlay refinements.}
		\label{pwh_GradientJump}
\end{figure}

\paragraph{Overlay refinement}
The initial mesh was constructed from \((4\times 8 \times 1)\) finite cells with a polynomial degree \(p=\{3,3,2\}\). After the first cycle, the mesh was refined in two successive cycles using the gradient-jump indicator (Eq.~\eqref{Kelly}), cf. Fig.~\ref{pwh_GradientJump}. The gradient jump between neighboring elements, also depicted in Fig.~\ref{pwh_GradientJump}, serves as a predictor for the evolution of plastic deformation and provides the basis for local refinement decisions according to Eq.~\eqref{Kelly}. The figure indicates an isotropic refinement, although the mesh could have been kept constant over the entire thickness.

\begin{figure}[h!]
	\centering
		\begin{tikzpicture}
			\node (m1) at (0.0cm, 0.0cm) {\includegraphics[width = 0.35\textwidth]{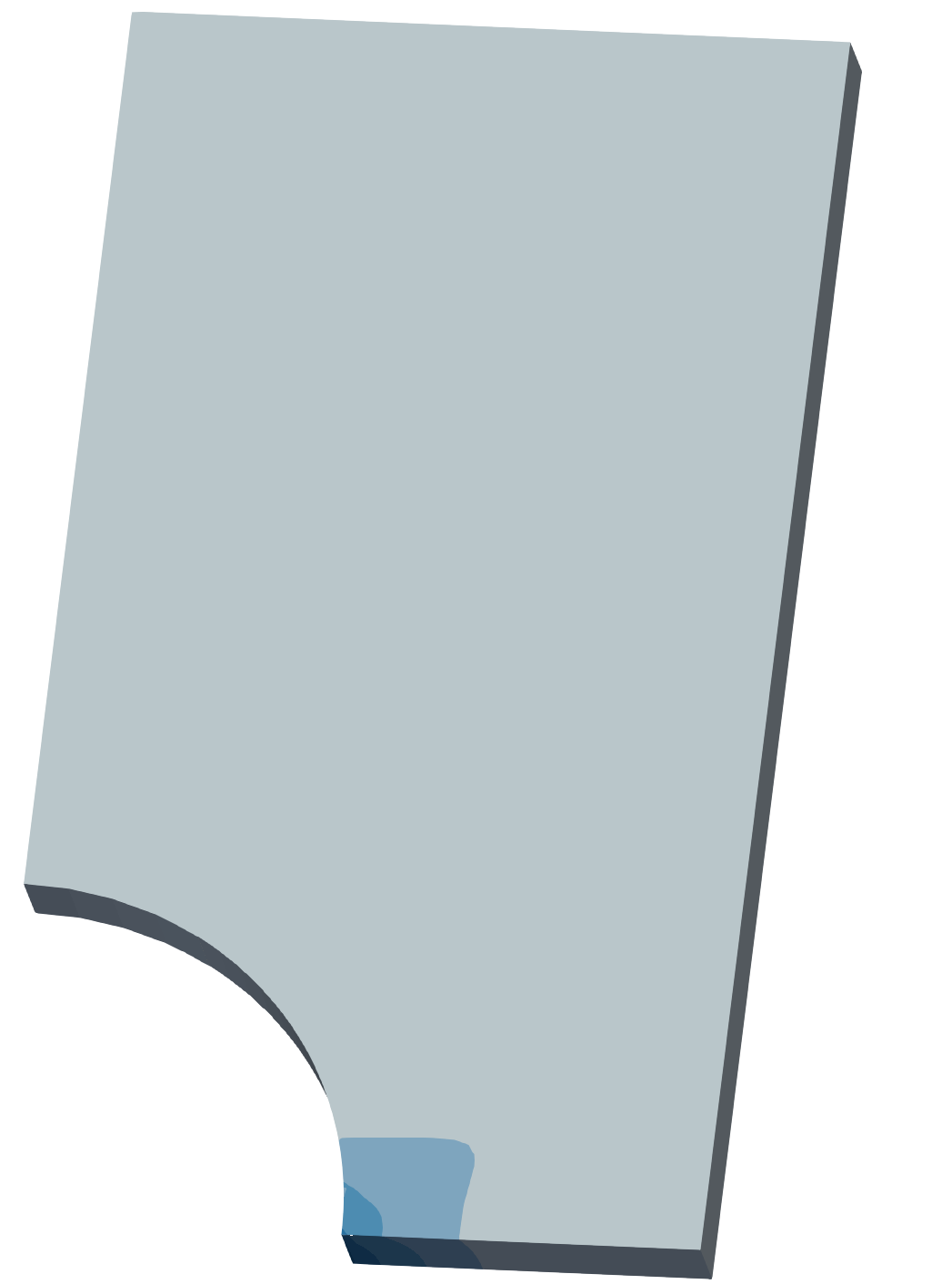}};
			\node (m2) at (2.5cm, 0.5cm) {\includegraphics[width = 0.35\textwidth]{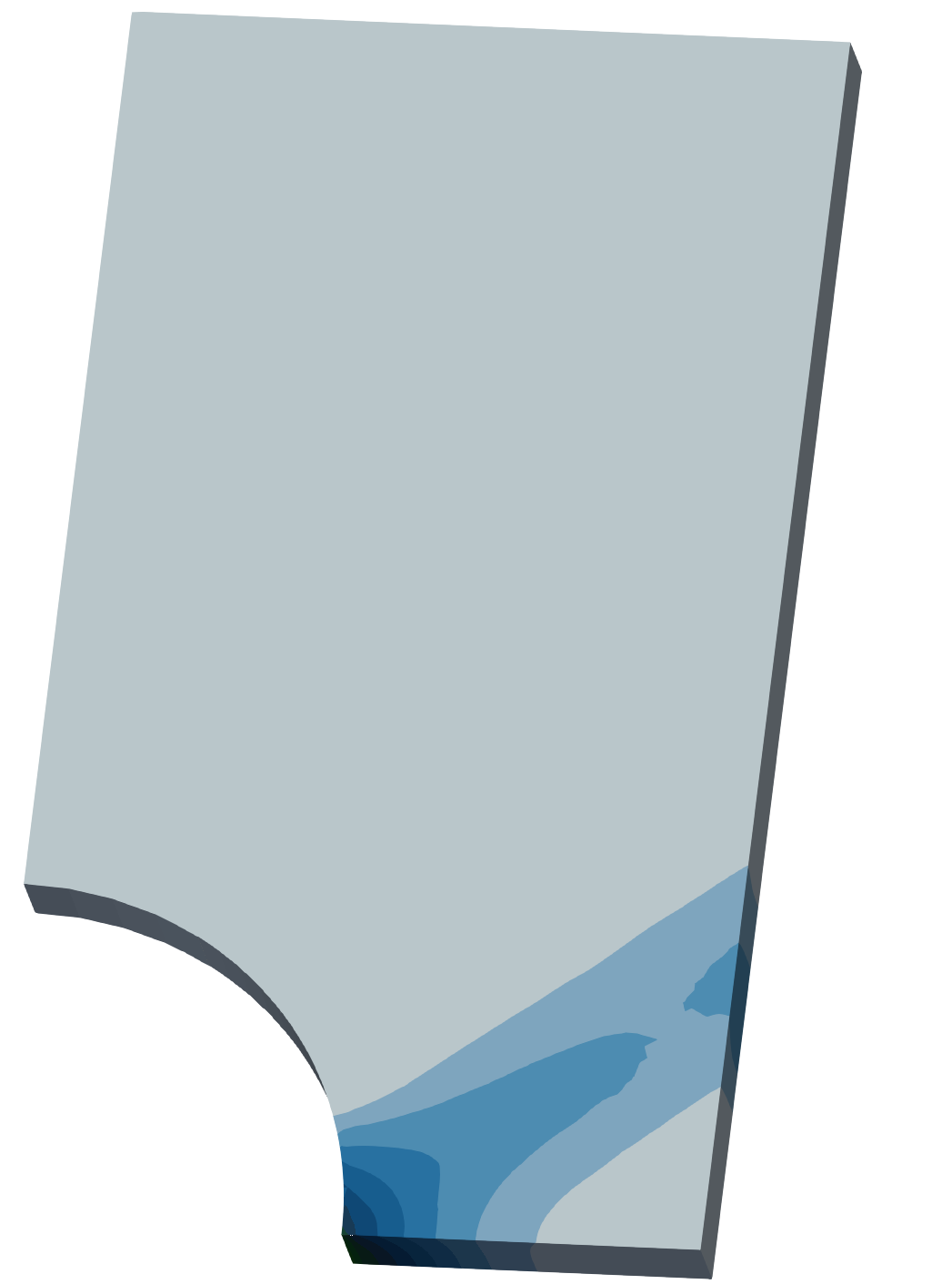}};
			\node (m3) at (5.0cm, 1.0cm) {\includegraphics[width = 0.35\textwidth]{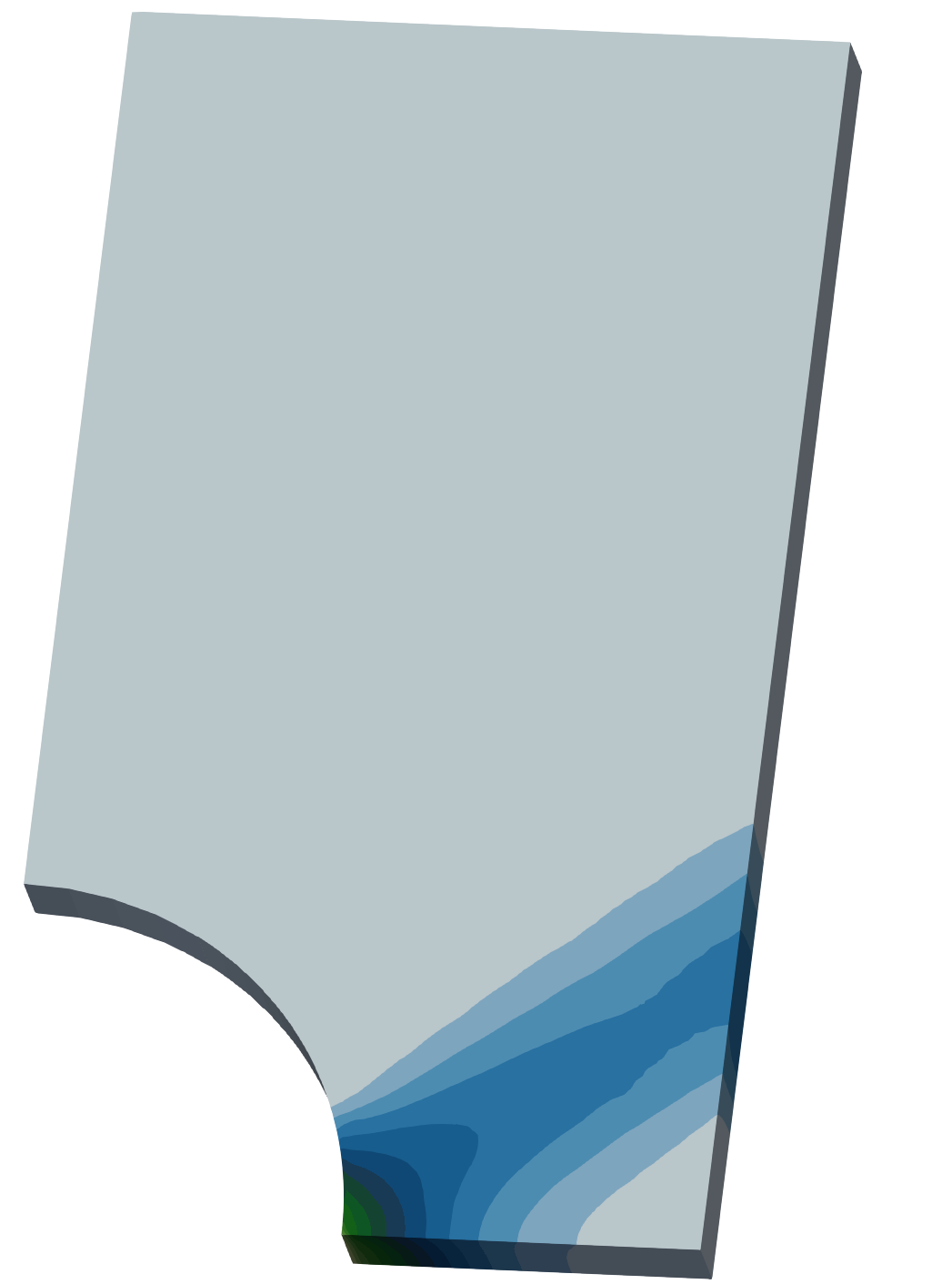}};
			\node (m4) at (7.5cm, 1.5cm) {\includegraphics[width = 0.35\textwidth]{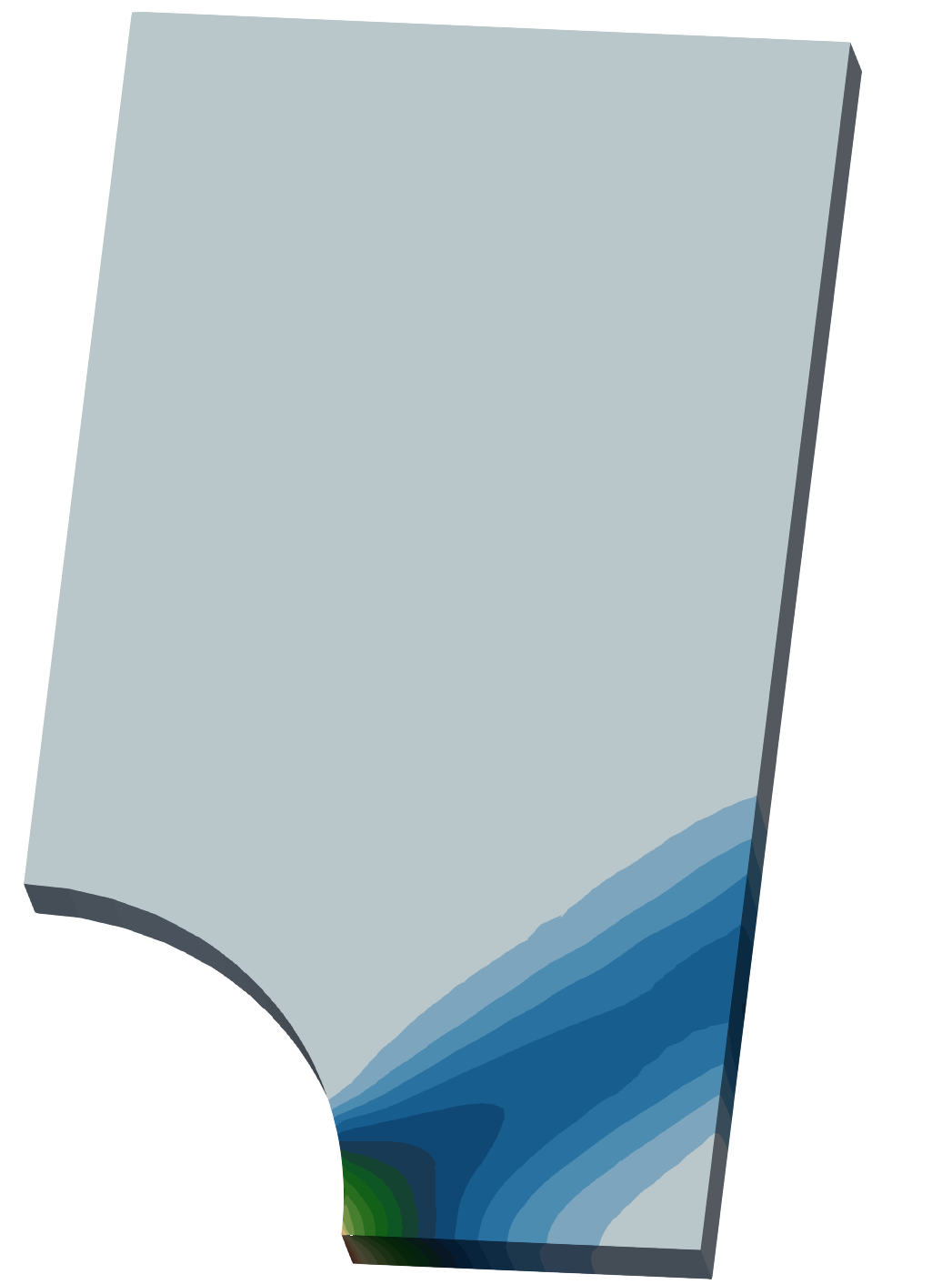}};
			\node (m5) at (10.0cm, 2.0cm) {\includegraphics[width = 0.35\textwidth]{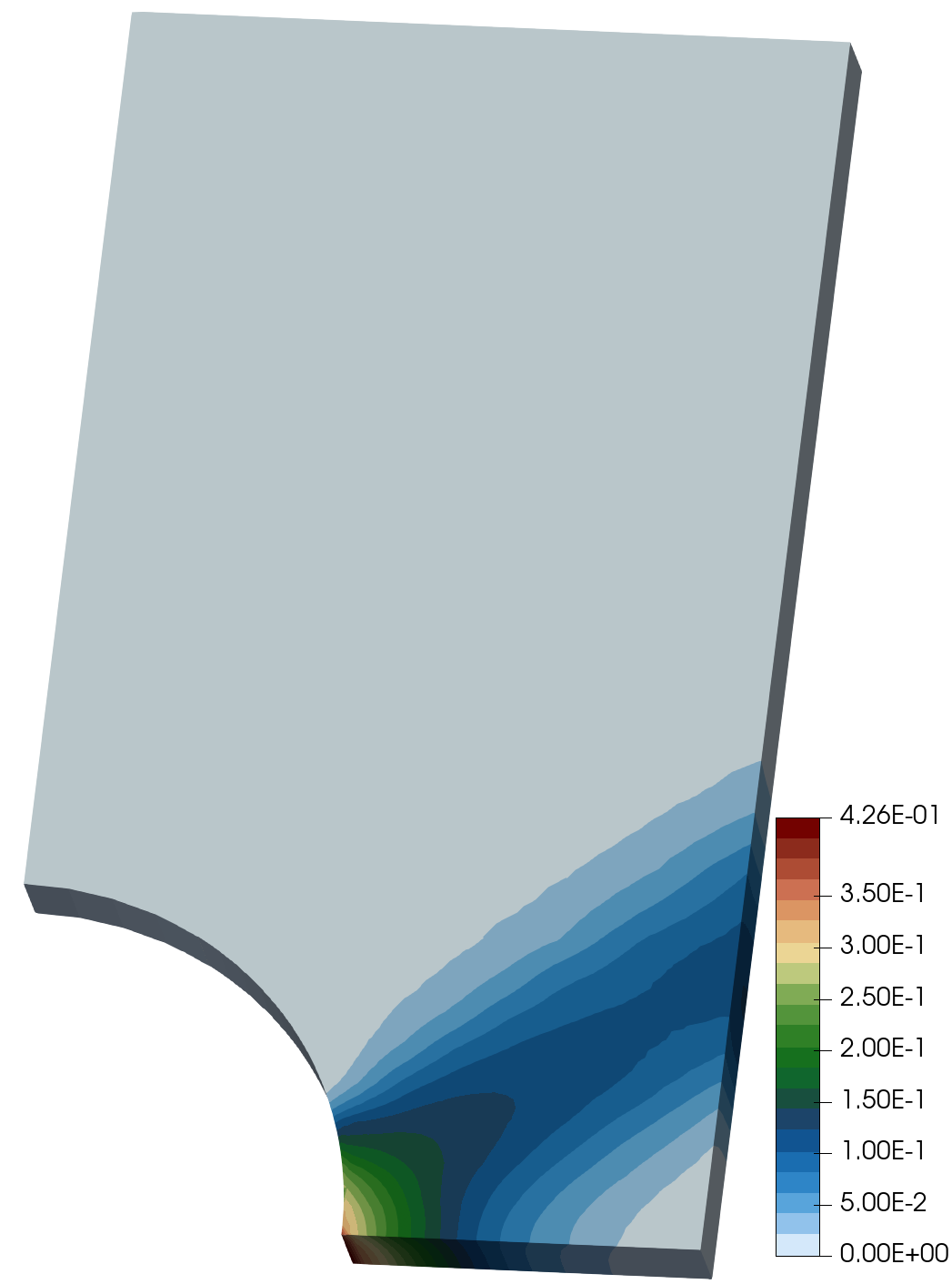}};

			\node[below] at (m1.south) {$\overline{u} = 0.1$};
			\node[below] at (m2.south) {$\overline{u} = 0.2$};
			\node[below] at (m3.south) {$\overline{u} = 0.3$};
			\node[below] at (m4.south) {$\overline{u} = 0.4$};
			\node[below] at (m5.south) {$\overline{u} = 0.5$};
		\end{tikzpicture}
		\caption{Plate with a circular hole: development of the accumulated plastic strain \(\bar{\varepsilon}^p\) with increasing prescribed displacement $\overline{u}\,\mathrm{[mm]}$.}
		\label{PPAccumulatedStrainProgression}
\end{figure}

The refinement strategy does not modify the discretization between individual load steps. Preliminary investigations showed that refinement based on indicators evaluated at intermediate load steps resulted in only negligible differences compared with the refinement pattern obtained from the final load step. \color{mygreenD} This indicates that, for the present example, the regions of localized plastic deformation developed during the loading history are sufficiently represented by the final load state. \color{black} Therefore, each adaptive cycle is performed on a fixed discretization over the complete loading history, and the error indicator is evaluated only after the final load step. The mesh is refined accordingly, and the complete analysis is subsequently restarted from the initial unloaded state on the refined discretization. This procedure avoids the transfer of history-dependent variables between successive meshes while retaining an adequate resolution of the evolving plastic zone, see Fig.~\ref{PPAccumulatedStrainProgression}.

\begin{figure}[h!]
	\centering
		\begin{tikzpicture}
			\node (m3) at (7.0cm, 0) {\includegraphics[width = 0.45\textwidth]{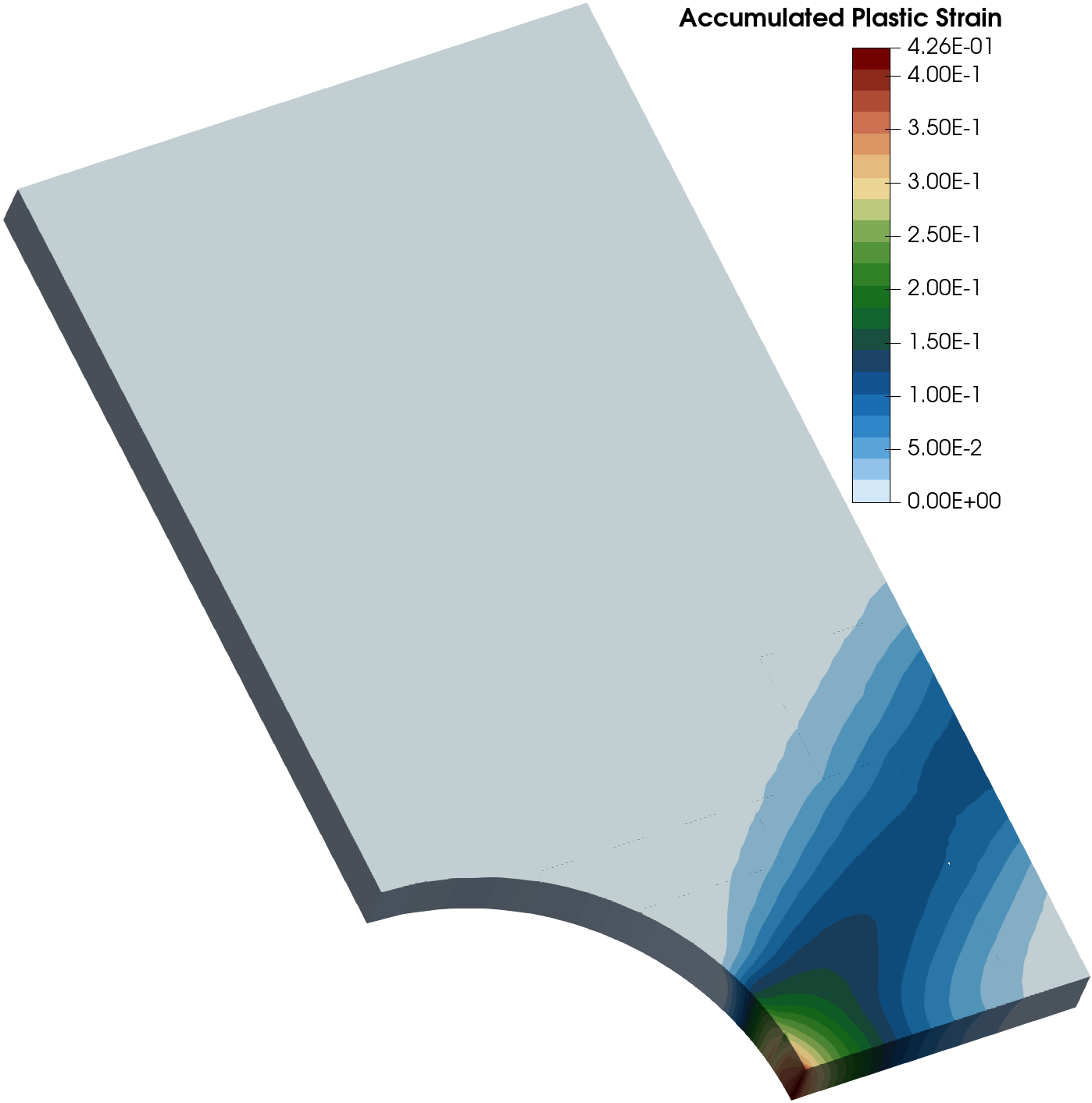}};
			\node (m1) at (0.0cm, 0) {\includegraphics[width = 0.45\textwidth]{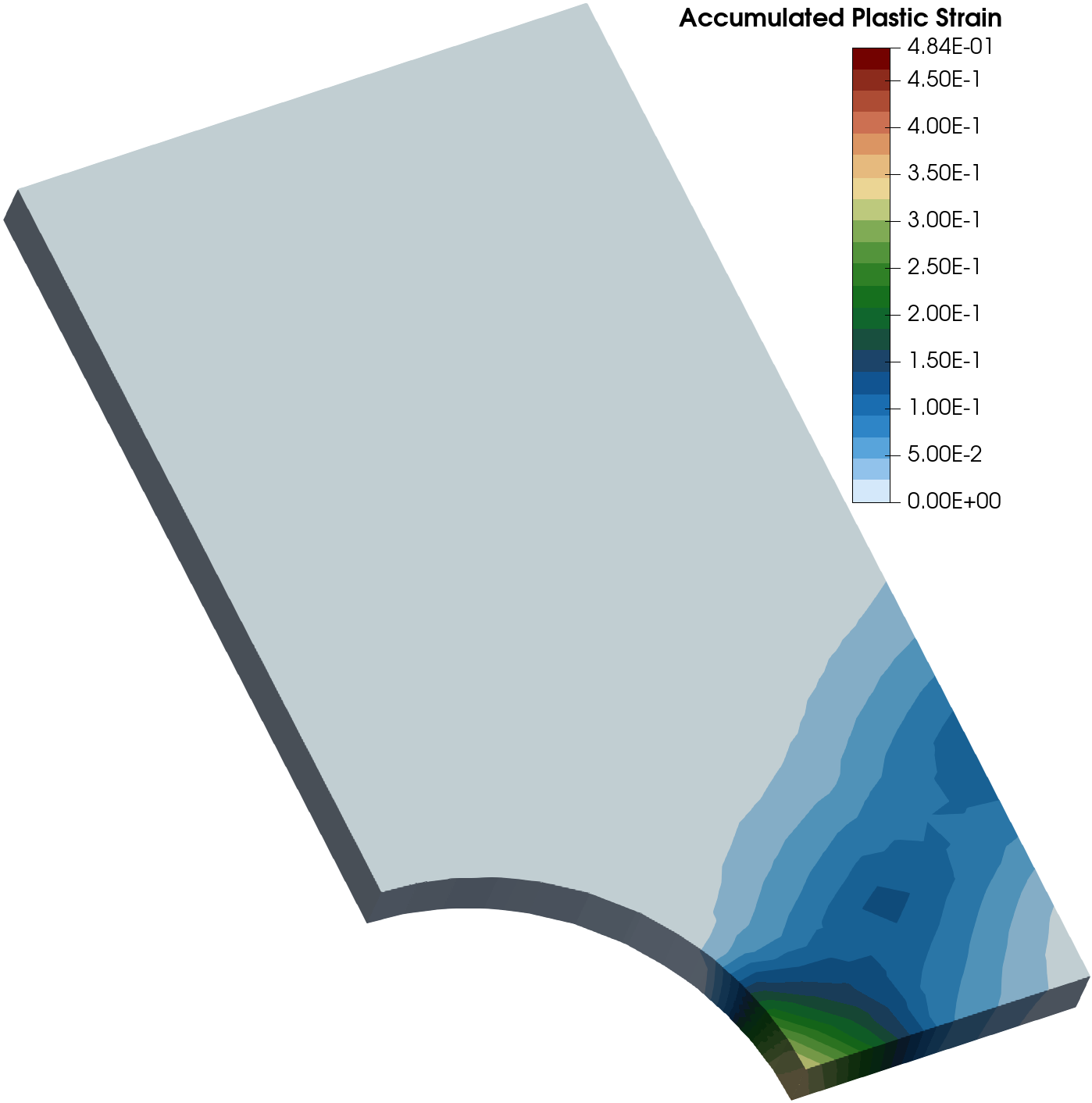}};

			\node[below] at (m1.south) {\footnotesize unrefined};
			\node[below] at (m3.south) {\footnotesize refinement $k = 2$};
		\end{tikzpicture}
		\caption{Plate with a circular hole: comparison of the accumulated plastic strain \(\bar{\varepsilon}^p\) of the last load step between the unrefined and the refined model according to the $k=2$ mesh of Fig. \ref{pwh_GradientJump}.}
		\label{AccumStrainComparison}
\end{figure}

A comparison of the accumulated plastic strain of the unrefined and refined model provides a good impression of the reliability and accuracy of the predicted plastic zone as well as the corresponding maximum values, which differ between \(12-14\%\), cf. Fig.~\ref{AccumStrainComparison}.

In Figure \ref{PPAccumulatedStrainProgression}, the progression of the propagating plastic deformation is shown for various loading steps. The results agree well with those shown in \cite{Simo:1986}. 
\begin{figure}[h!]
	\centering
		\begin{tikzpicture}[spy using outlines = {densely dashed, rectangle, magnification=7, connect spies}]

	\begin{axis} [
		font={\footnotesize},
		axis lines = box,
		xlabel = Edge displacement $\bar{u}$ in  mm,
		ylabel = Edge reaction force in  kN,
		width = .9\textwidth,
		height = 7.5cm,
		xtick distance = 0.1,
		ytick distance = 0.5,
		grid = major,
		grid style = {densely dashed, line width = 0.1pt},
		cycle list name = color_list_mark,
		legend pos = south east,
		legend style = {nodes = {scale = 0.65, transform shape}},
		legend cell align={left},
		enlarge x limits=0.05,
		enlarge y limits=0.05,
	]
		\addplot+[color = black, mark = none, cTab6] table [
						x = u,
						y = f,
						col sep = comma] {reaction_perforated_plate_referenceR2.csv};
		\addlegendentryexpanded{Reference FEM CST plane stress};

		\addplot+[color = black, mark = x, mark size = 1.7pt] table [x = u,
										y = f,
										col sep = comma] {reaction_perforated_plate_femR2.csv};
		\addlegendentryexpanded{FEM $6\times10\times2$ hexahedra};

		\pgfplotsset{cycle list shift = -2}

		\addplot+[mark = square*, densely dashed] table [x = u,
				y = f,
				col sep = comma] {reaction_perforated_plate_fcm_ast_k0R2.csv};
		\addlegendentryexpanded{FCM AST cycle $k = 0$};

		\addplot+[mark = square*, densely dashed] table [x = u,
				y = f,
				col sep = comma] {reaction_perforated_plate_fcm_ast_k1R2.csv};
		\addlegendentryexpanded{FCM AST cycle $k = 1$};

		\addplot+[mark = square*, densely dashed] table [x = u,
				y = f,
				col sep = comma] {reaction_perforated_plate_fcm_ast_k2R2.csv};
		\addlegendentryexpanded{FCM AST cycle $k = 2$};

		\addplot+[mark = square*, densely dashed] table [x = u,
				y = f,
				col sep = comma] {reaction_perforated_plate_fcm_ast_k3R2.csv};
		\addlegendentryexpanded{FCM AST cycle $k = 3$};

		\pgfplotsset{cycle list shift = -6}

		\addplot+[mark = *] table [x = u,
										y = f,
										col sep = comma] {reaction_perforated_plate_fcm_nnmf_k0R2.csv};
		\addlegendentryexpanded{FCM NNMF cycle $k = 0$};

		\addplot+[mark = *] table [x = u,
										y = f,
										col sep = comma] {reaction_perforated_plate_fcm_nnmf_k1R2.csv};
		\addlegendentryexpanded{FCM NNMF cycle $k = 1$};

		\addplot+[mark = *] table [x = u,
										y = f,
										col sep = comma] {reaction_perforated_plate_fcm_nnmf_k2R2.csv};
		\addlegendentryexpanded{FCM NNMF cycle $k = 2$};

		\addplot+[mark = *] table [x = u,
										y = f,
										col sep = comma] {reaction_perforated_plate_fcm_nnmf_k3R2.csv};
		\addlegendentryexpanded{FCM NNMF cycle $k = 3$};

		\coordinate (spypoint1) at (axis cs:0.588, 2.76);
		\coordinate (magnifyglass1) at (axis cs:0.5, 2.15);

		\coordinate (spypoint2) at (axis cs:0.07, 2.448);
		\coordinate (magnifyglass2) at (axis cs:0.19, 1.75);

	\end{axis}

	\spy[width = 4.8cm, height = 1.5cm] on (spypoint1) in node[fill = white] at (magnifyglass1);
	\spy[width = 5cm, height = 2cm] on (spypoint2) in node[fill = white] at (magnifyglass2);

	\begin{axis} [
		yshift = -3.8cm,
		font={\footnotesize},
		axis lines = box,
		xlabel = Load step,
		ylabel = Global iterations,
		width = .9\textwidth,
		height = 4cm,
		max space between ticks = 20,
		xtick distance = 5,
		grid = major,
		grid style = {densely dashed, line width = 0.1pt},
		cycle list name = color_list_mark,
		legend pos = outer north east,
		legend style = {nodes = {scale = 0.65, transform shape}},
		legend cell align={left},
		enlarge x limits=0.05,
	]

		\addplot+[black, mark = x, mark size = 1.7pt] table [x expr=\coordindex,
				y = i,
				col sep = comma] {reaction_perforated_plate_femR2.csv};
		\addlegendentryexpanded{FEM};

		\pgfplotsset{cycle list shift = -1}

		\addplot+[mark = square*, densely dashed] table [x expr=\coordindex,
				y = i,
				col sep = comma] {reaction_perforated_plate_fcm_ast_k0R2.csv};
		\addlegendentryexpanded{FCM AST cycle $k = 0$};

		\addplot+[mark = square*, densely dashed] table [x expr=\coordindex,
				y = i,
				col sep = comma] {reaction_perforated_plate_fcm_ast_k1R2.csv};
		\addlegendentryexpanded{FCM AST cycle $k = 1$};

		\addplot+[mark = square*, densely dashed] table [x expr=\coordindex,
				y = i,
				col sep = comma] {reaction_perforated_plate_fcm_ast_k2R2.csv};
		\addlegendentryexpanded{FCM AST cycle $k = 2$};

		\addplot+[mark = square*, densely dashed] table [x expr=\coordindex,
				y = i,
				col sep = comma] {reaction_perforated_plate_fcm_ast_k3R2.csv};
		\addlegendentryexpanded{FCM AST cycle $k = 3$};

		\pgfplotsset{cycle list shift = -5}

		\addplot+[mark = *] table [x expr=\coordindex,
				y = i,
				col sep = comma] {reaction_perforated_plate_fcm_nnmf_k0R2.csv};
		\addlegendentryexpanded{FCM NNMF cycle $k = 0$};

		\addplot+[mark = *] table [x expr=\coordindex,
				y = i,
				col sep = comma] {reaction_perforated_plate_fcm_nnmf_k1R2.csv};
		\addlegendentryexpanded{FCM NNMF cycle $k = 1$};

		\addplot+[mark = *] table [x expr=\coordindex,
				y = i,
				col sep = comma] {reaction_perforated_plate_fcm_nnmf_k2R2.csv};
		\addlegendentryexpanded{FCM NNMF cycle $k = 2$};

		\addplot+[mark = *] table [x expr=\coordindex,
				y = i,
				col sep = comma] {reaction_perforated_plate_fcm_nnmf_k3R2.csv};
		\addlegendentryexpanded{FCM NNMF cycle $k = 3$};

		\legend{};

	\end{axis}

	\begin{scope}[shift = {(6.8, .4)}, isometric view, scale = 0.3, line join=round]
		\fill[cSecond, opacity = 0.2] (0, 9, 0) -- ++(5, 0, 0) -- ++(0, 0, 0.5) -- ++(-5, 0, 0) -- cycle;
		\draw[cSecond] (0, 9, 0) -- ++(5, 0, 0) -- ++(0, 0, 0.5) -- ++(-5, 0, 0) -- cycle;

		\node[above, font = \footnotesize, fill = white, fill opacity = 0.8, text opacity = 1, inner sep = 1, outer sep = 2] at (2.5, 9, 0.6) {$\overline{u}$};

		\foreach \x in {0, 1, ..., 5}
			{
				\foreach \z in {0, 0.5}
					{
						\draw[cSecond, -stealth] (\x, 8, \z) -- ++(0, 1, 0);
					}
			}

		\draw[black] (5, 8, 0) -- ++(0, 0, 0.5);
		\begin{scope}[canvas is xy plane at z = 0]
			\draw[black] (2.5, 0) -- ++(2.5, 0) -- ++(0, 8) -- ++(-5, 0) -- ++(0, -5.5) arc(90:0:2.5) -- cycle;
		\end{scope}

		\begin{scope}[canvas is xy plane at z = 0.25]
			\draw[
				black,
				top color = bulk,
				bottom color = bulk,
				middle color = white,
				shading angle = 15
			] (2.5, 0) -- node[pos = 0.5, anchor = south east, inner sep = 1.2, cFirst, opacity = .65] {} ++(2.5, 0) -- ++(0, 8) -- ++(-5, 0) -- node[pos = 0.5, anchor = south west, inner sep = 2, cFirst, opacity = .65] {} ++(0, -5.5) arc(90:0:2.5) -- cycle;
		\end{scope}

		\begin{scope}[canvas is xy plane at z = 0.5]
			\draw[black] (2.5, 0) -- ++(2.5, 0) -- ++(0, 8) -- ++(-5, 0) -- ++(0, -5.5) arc(90:0:2.5) -- cycle;
		\end{scope}

		\draw[black] (2.5, 0, 0) -- ++(0, 0, 0.5);
		\draw[black] (0, 2.5, 0) -- ++(0, 0, 0.5);
		\draw[black] (5, 0, 0) -- ++(0, 0, 0.5);
		\draw[black] (0, 8, 0) -- ++(0, 0, 0.5);

		\begin{scope}[canvas is xz plane at y = 0]
			\fill[cFirst, opacity = .0] (2.5, 0) rectangle (5, 0.5);
		\end{scope}
		\begin{scope}[canvas is yz plane at x = 0]
			\fill[cFirst, opacity = .0] (8, 0) rectangle (2.5, 0.5);
		\end{scope}

		\draw[black, densely dashed] (0, 0, 0) -- ++(2.5, 0, 0);
		\draw[black, densely dashed] (0, 0, 0) -- ++(0, 2.5, 0);
		\draw[black, densely dashed] (0, 0, 0.5) -- ++(2.5, 0, 0);
		\draw[black, densely dashed] (0, 0, 0.5) -- ++(0, 2.5, 0);

	\end{scope}

\end{tikzpicture}
		\caption{Plate with a circular hole: convergence of the reaction force along the loaded edge.}
		\label{PPIntegrationConvergence}
\end{figure}

\vspace{3mm}
\paragraph{Non-negative moment fitting}
The \emph{Non-Negative Moment Fitting (NNMF)} from Section \ref{sec:NNMFquadrature} was applied and compared to the established \emph{Adaptive Space-Tree quadrature (AST)}. The convergence of the computation is depicted in Fig. \ref{PPIntegrationConvergence} (top) considering three successive refinement steps. In terms of solution quality, there is practically no difference between the two integration methods, which demonstrates the equivalent accuracy of the NNMF approach compared to the established standard. Both results are in good agreement with those reported in \emph{de Souza Neto et al.} \cite{deSouzaNeto:2011}. The numerical stability of the computation is illustrated in Fig. \ref{PPIntegrationConvergence} (bottom) with the number of iterations per load step to ensure a reliable convergence level. Although the NNMF approach occasionally requires slightly more iterations than the computationally more expensive AST approach, both schemes converge with a comparable iteration count, requiring roughly 20\% fewer iterations than the FEM reference solution. The significantly higher numerical effort per step of the \(hp\)-refined results (red curve) offsets this saving, yielding a substantially more accurate solution in terms of accuracy and resolution of the physical response. It is worth mentioning that the unrefined FCM model ($p=\{3,3,2\}$) with results comparable to the FEM reference model requires \(42\%\) fewer iterations.

Table \ref{tab:quadpoints} reveals the number of quadrature points used in the different integration schemes, AST and NNMF, respectively. The reduction in the number of quadrature points is substantial for all considered discretizations and amounts to approximately \(75\%\) for the more highly refined models. The computational benefit of NNMF for a geometrically more complex three-dimensional problem is investigated separately in Section~\ref{sec:porousPlate}, where its efficiency is assessed directly in terms of execution time. Overall, the analyses considered in this work demonstrate substantial reductions in numerical integration effort, making NNMF computationally superior to AST for the non-linear finite cell analyses considered herein.
\begin{table}
\begin{spacing}{1.2}
\hrule
\centering
\begin{tabular}{crrrrrrrrr}
	
	\multicolumn{4}{c}{\textbf{Total}} & \multicolumn{3}{r}{\textbf{Physical domain}} & \multicolumn{3}{c}{\hspace{4mm}\textbf{Fictitious domain} }  \\ 
	Cycle $i$ & \multicolumn{1}{c}{AST} & \multicolumn{1}{c}{NNMF} & \multicolumn{1}{c}{Reduction} & \hspace{2mm} \ & \multicolumn{1}{c}{AST} & \multicolumn{1}{r}{NNMF} & \hspace{2mm} \ & \multicolumn{1}{c}{AST} & \multicolumn{1}{r}{NNMF} \\\hline
	0 & 20\,412 &  2\,468 & 87.91\% && 20\,142 &  2\,198 &&    270 & 270\\
	1 & 38\,031 &  6\,039 & 84.12\% && 36\,477 &  4\,485 && 1\,444 & 1\,554 \\
	2 & 58\,544 & 14\,389 & 75.42\% && 54\,432 & 10\,051 && 4\,112 & 4\,338 \\
	3 & 58\,936 & 14\,781 & 74.92\% && 54\,432 & 10\,051 && 4\,504 & 4\,730
\end{tabular}
\hrule
\end{spacing}
\caption{Plate with a circular hole: comparison of the number of quadrature points for AST and NNMF quadrature.}
\label{tab:quadpoints}
\end{table}

\paragraph{Viscoplasticity}
We consider the benchmark for viscoplastic analysis based on the material model introduced in Section \ref{sec:Thermo-viscoplasticity}, but employ the flow rule proposed by Peri\'{c} \cite{Peric:1993} 
to ensure consistency with the reference results of \cite{deSouzaNeto:2011}: 
\begin{equationarray}{rcl}
	\dot{\gamma}(\boldsymbol{\sigma},\bar{\varepsilon}^p) &=& 
	\begin{cases}
		\frac{1}{\mu}\left[\biggl(\frac{\sigma_v(\boldsymbol{\sigma})}{\sigma_y(\bar{\varepsilon}^p)}\biggr)^{\frac{1}{m}}-1\right]   & \mbox{if }\;\; \Phi(\boldsymbol{\sigma},\bar{\varepsilon}^p)\geq 0\\
		0 & \mbox{if }\;\; \Phi(\boldsymbol{\sigma},\bar{\varepsilon}^p)< 0 \,.\hfill \\
	\end{cases}
\end{equationarray}

\noindent The relaxation time was held constant at \(\mu = 500 \,\mathrm{s}\), while viscoplastic exponents of \(m=0.1\) and \(m=1.0\) were considered for varying deformation rates \(\dot{\bar{u}}\). The remaining material parameters were chosen identical to the elasto-plastic model.
\begin{figure}[h!]
	\centering
		\begin{tabular}{@{} c @{} c @{}}
			\pgfplotsset{
    legend image with text/.style={
        legend image code/.code={%
            \node[anchor=west, align = flush right, text width = 1cm] at (0.0cm,0cm) {#1};
        }
    },
}

\begin{tikzpicture}[spy using outlines = {densely dashed, rectangle, magnification=7, connect spies}]

	\begin{axis} [
		font={\footnotesize},
		axis lines = box,
		xlabel = Edge displacement $\bar{u}$ in  mm,
		ylabel = Edge reaction force in  kN,
		width = 0.5\textwidth,
		height = 12cm,
		xtick distance = 0.1,
		ytick distance = 2,
		extra y ticks = {-1, 1, 3, ..., 15},
		extra y tick labels = {},
		grid = major,
		grid style = {densely dashed, line width = 0.1pt},
		legend style = {at={(1,1.02)}, anchor=south east, nodes = {scale = 0.55, transform shape},%
										/tikz/every even column/.append style={column sep = 3.5mm}},
		transpose legend,
		legend cell align={left},
		legend columns = 5,
		mark size = 1.25pt,
		enlarge x limits=0.05,
		enlarge y limits=0.05,
	]

	  \addlegendimage{legend image with text = {~}}
		\addlegendentryexpanded{Reference}

		\addplot+[cTab3, mark = none, densely dashed, thick] table [
						x = u4,
						y = 0_555_E-02,
						col sep = comma] {reaction_perforated_plate_rate1.0_referenceR2.csv};
		\addlegendentryexpanded{$\dot{\bar{u}} = 10^{-1}\,\sfrac{\mathrm{mm}}{\mathrm{s}}$};

		\addplot+[cTab2, mark = none, densely dashed, thick] table [
						x = u3,
						y = 0_555_E-03,
						col sep = comma] {reaction_perforated_plate_rate1.0_referenceR2.csv};
		\addlegendentryexpanded{$\dot{\bar{u}} = 10^{-2}\,\sfrac{\mathrm{mm}}{\mathrm{s}}$};

		\addplot+[cTab1, mark = none, densely dashed, thick] table [
						x = u2,
						y = 0_555_E-04,
						col sep = comma] {reaction_perforated_plate_rate1.0_referenceR2.csv};
		\addlegendentryexpanded{$\dot{\bar{u}} = 10^{-3}\,\sfrac{\mathrm{mm}}{\mathrm{s}}$};

		\addplot+[mark = none, black, thick] table [
						x = u,
						y = f,
						col sep = comma] {reaction_perforated_plate_rateindependent_referenceR2.csv}; 
		\addlegendentryexpanded{$\dot{\bar{u}} \rightarrow 0\,\sfrac{\mathrm{mm}}{\mathrm{s}}$};

		\addlegendimage{legend image with text = {$\dot{\bar{u}} =$}}
		\addlegendentryexpanded{$\mspace{-2mu}10^{-3}\,\sfrac{\mathrm{mm}}{\mathrm{s}}$}

		\addplot+[solid, cTab20c4, mark = triangle*, mark options = {solid, fill = cTab20c4!50!white}] table [x = u,
										y = f,
										col sep = comma] {reaction_perforated_plate_rate1.0_mu0.55E-4_k0R2.csv};
		\addlegendentryexpanded{cycle $k = 0$};

		\addplot+[solid, cTab20c3, mark = square*, mark options = {solid, fill = cTab20c3!50!white}] table [x = u,
										y = f,
										col sep = comma] {reaction_perforated_plate_rate1.0_mu0.55E-4_k1R2.csv};
		\addlegendentryexpanded{cycle $k = 1$};

		\addplot+[solid, cTab20c2, mark = diamond*, mark options = {solid, fill = cTab20c2!50!white}] table [x = u,
										y = f,
										col sep = comma] {reaction_perforated_plate_rate1.0_mu0.55E-4_k2R2.csv};
		\addlegendentryexpanded{cycle $k = 2$};

		\addplot+[solid, cTab20c1, mark = *,  mark options = {solid, fill = cTab20c1!50!white}] table [x = u,
										y = f,
										col sep = comma] {reaction_perforated_plate_rate1.0_mu0.55E-4_k3R2.csv};
		\addlegendentryexpanded{cycle $k = 3$};

		\addlegendimage{legend image with text = {$\dot{\bar{u}} = $}}
		\addlegendentryexpanded{$\mspace{-2mu}10^{-2}\,\sfrac{\mathrm{mm}}{\mathrm{s}}$}

		\addplot+[solid, cTab20c8, mark = triangle*, mark options = {solid, fill = cTab20c8!50!white}] table [x = u,
										y = f,
										col sep = comma] {reaction_perforated_plate_rate1.0_mu0.55E-3_k0R2.csv};
		\addlegendentryexpanded{cycle $k = 0$};

		\addplot+[solid, cTab20c7, mark = square*, mark options = {solid, fill = cTab20c7!50!white}] table [x = u,
										y = f,
										col sep = comma] {reaction_perforated_plate_rate1.0_mu0.55E-3_k1R2.csv};
		\addlegendentryexpanded{cycle $k = 1$};

		\addplot+[solid, cTab20c6, mark = diamond*, mark options = {solid, fill = cTab20c6!50!white}] table [x = u,
										y = f,
										col sep = comma] {reaction_perforated_plate_rate1.0_mu0.55E-3_k2R2.csv};
		\addlegendentryexpanded{cycle $k = 2$};

		\addplot+[solid, cTab20c5, mark = *, mark options = {solid, fill = cTab20c5!50!white}] table [x = u,
										y = f,
										col sep = comma] {reaction_perforated_plate_rate1.0_mu0.55E-3_k3R2.csv};
		\addlegendentryexpanded{cycle $k = 3$};

		\addlegendimage{legend image with text = {\hfill $\dot{\bar{u}} = $}}
		\addlegendentryexpanded{$\mspace{-2mu}10^{-1}\,\sfrac{\mathrm{mm}}{\mathrm{s}}$}

		\addplot+[solid, cTab20c12, mark = triangle*, mark options = {solid, fill = cTab20c12!50!white}] table [x = u,
										y = f,
										col sep = comma] {reaction_perforated_plate_rate1.0_mu0.55E-2_k0R2.csv};
		\addlegendentryexpanded{cycle $k = 0$};

		\addplot+[solid, cTab20c11, mark = square*, mark options = {solid, fill = cTab20c11!50!white}] table [x = u,
										y = f,
										col sep = comma] {reaction_perforated_plate_rate1.0_mu0.55E-2_k1R2.csv};
		\addlegendentryexpanded{cycle $k = 1$};

		\addplot+[solid, cTab20c10, mark = diamond*, mark options = {solid, fill = cTab20c10!50!white}] table [x = u,
										y = f,
										col sep = comma] {reaction_perforated_plate_rate1.0_mu0.55E-2_k2R2.csv};
		\addlegendentryexpanded{cycle $k = 2$};

		\addplot+[solid, cTab20c9, mark = *, mark options = {solid, fill = cTab20c9!50!white}] table [x = u,
										y = f,
										col sep = comma] {reaction_perforated_plate_rate1.0_mu0.55E-2_k3R2.csv};
		\addlegendentryexpanded{cycle $k = 3$};

		\coordinate (spypoint1) at (axis cs:0.582, 12.42);
		\coordinate (magnifyglass1) at (axis cs:0.443, 8.65);

		\coordinate (spypoint2) at (axis cs:0.582, 4.28);
		\coordinate (magnifyglass2) at (axis cs:0.443, 6);

		\coordinate (spypoint3) at (axis cs:0.582, 2.99);
		\coordinate (magnifyglass3) at (axis cs:0.443, 1);

		\coordinate (l1) at (axis cs:0.2, 2.2);
		\coordinate (l2) at (axis cs:0.22, 3.2);
		\coordinate (l3) at (axis cs:0.43, 3.8);
		\coordinate (l4) at (axis cs:0.4, 12.5);
	\end{axis}

	\spy[width = 3cm, height = 1.5cm, magnification=5.5] on (spypoint1) in node[fill = white] at (magnifyglass1);
	\spy[width = 3cm, height = 1.5cm, magnification=5.5] on (spypoint2) in node[fill = white] at (magnifyglass2);
	\spy[width = 3cm, height = 1.5cm, magnification=5.5] on (spypoint3) in node[fill = white] at (magnifyglass3);

	\node[font = \scriptsize] at (l1) {(i)~$\dot{\overline{u}} = 0\,\mathrm{mm}/\mathrm{s}$};
	\node[font = \scriptsize] at (l2) {(ii)~$\dot{\overline{u}} = 10^{-3}\,\mathrm{mm}/\mathrm{s}$};
	\node[font = \scriptsize] at (l3) {(iii)~$\dot{\overline{u}} = 10^{-2}\,\mathrm{mm}/\mathrm{s}$};
	\node[font = \scriptsize] at (l4) {(iv)~$\dot{\overline{u}} = 10^{-1}\,\mathrm{mm}/\mathrm{s}$};

	\node[text width = 1.6cm] at (1.2, 9.75) {\footnotesize$m = 1.0$\\$\mu = 500\,\mathrm{s}$};

	\begin{scope}[shift = {(1.5, 7.5)}, isometric view, scale = 0.2, line join=round]
		\fill[cSecond, opacity = 0.2] (0, 9, 0) -- ++(5, 0, 0) -- ++(0, 0, 0.5) -- ++(-5, 0, 0) -- cycle;
		\draw[cSecond] (0, 9, 0) -- ++(5, 0, 0) -- ++(0, 0, 0.5) -- ++(-5, 0, 0) -- cycle;

		\node[above, font = \footnotesize, fill = white, fill opacity = 0.8, text opacity = 1, inner sep = 1, outer sep = 2] at (2.5, 9, 0.6) {$\overline{u}$};

		\foreach \x in {0, 1, ..., 5}
			{
				\foreach \z in {0, 0.5}
					{
						\draw[cSecond, -stealth] (\x, 8, \z) -- ++(0, 1, 0);
					}
			}

		\draw[black] (5, 8, 0) -- ++(0, 0, 0.5);
		\begin{scope}[canvas is xy plane at z = 0]
			\draw[black] (2.5, 0) -- ++(2.5, 0) -- ++(0, 8) -- ++(-5, 0) -- ++(0, -5.5) arc(90:0:2.5) -- cycle;
		\end{scope}

		\begin{scope}[canvas is xy plane at z = 0.25]
			\draw[
				black,
				top color = bulk,
				bottom color = bulk,
				middle color = white,
				shading angle = 15
			] (2.5, 0) -- node[pos = 0.5, anchor = south east, inner sep = 1.2, cFirst, opacity = .65] {} ++(2.5, 0) -- ++(0, 8) -- ++(-5, 0) -- node[pos = 0.5, anchor = south west, inner sep = 2, cFirst, opacity = .65] {} ++(0, -5.5) arc(90:0:2.5) -- cycle;
		\end{scope}

		\begin{scope}[canvas is xy plane at z = 0.5]
			\draw[black] (2.5, 0) -- ++(2.5, 0) -- ++(0, 8) -- ++(-5, 0) -- ++(0, -5.5) arc(90:0:2.5) -- cycle;
		\end{scope}

		\draw[black] (2.5, 0, 0) -- ++(0, 0, 0.5);
		\draw[black] (0, 2.5, 0) -- ++(0, 0, 0.5);
		\draw[black] (5, 0, 0) -- ++(0, 0, 0.5);
		\draw[black] (0, 8, 0) -- ++(0, 0, 0.5);

		\begin{scope}[canvas is xz plane at y = 0]
			\fill[cFirst, opacity = .0] (2.5, 0) rectangle (5, 0.5);
		\end{scope}
		\begin{scope}[canvas is yz plane at x = 0]
			\fill[cFirst, opacity = .0] (8, 0) rectangle (2.5, 0.5);
		\end{scope}

		\draw[black, densely dashed] (0, 0, 0) -- ++(2.5, 0, 0);
		\draw[black, densely dashed] (0, 0, 0) -- ++(0, 2.5, 0);
		\draw[black, densely dashed] (0, 0, 0.5) -- ++(2.5, 0, 0);
		\draw[black, densely dashed] (0, 0, 0.5) -- ++(0, 2.5, 0);

	\end{scope}

\end{tikzpicture} &
			\pgfplotsset{
    legend image with text/.style={
        legend image code/.code={%
            \node[anchor=west, align = flush right, text width = 1cm] at (0.0cm,0cm) {#1};
        }
    },
}

\begin{tikzpicture}[spy using outlines = {densely dashed, rectangle, magnification=14, connect spies}]

	\begin{axis} [
		font={\footnotesize},
		axis lines = box,
		xlabel = Edge displacement $\bar{u}$ in  mm,
		ylabel = Edge reaction force in  kN,
		width = 0.5\textwidth,
		height = 12cm,
		xtick distance = 0.1,
		ytick distance = 1,
		grid = major,
		grid style = {densely dashed, line width = 0.1pt},
		legend style = {at={(1,1.02)}, anchor=south east, nodes = {scale = 0.55, transform shape},%
										/tikz/every even column/.append style={column sep = 3.5mm}},
		transpose legend,
		legend cell align={left},
		legend columns = 5,
		enlarge x limits=0.05,
		enlarge y limits=0.05,
		mark size = 1.25pt,
		ymax =10.
	]

	  \addlegendimage{legend image with text = {~}}
		\addlegendentryexpanded{Reference}

		\addplot+[cTab3, mark = none, densely dashed, thick] table [
						x = u8,
						y = 0_555_E+02,
						col sep = comma] {reaction_perforated_plate_rate0.01_referenceR2.csv};
		\addlegendentryexpanded{$\dot{\bar{u}} = 10^3\,\sfrac{\mathrm{mm}}{\mathrm{s}}$};

		\addplot+[cTab2, mark = none, densely dashed, thick] table [
						x = u6,
						y = 0_555_E+00,
						col sep = comma] {reaction_perforated_plate_rate0.01_referenceR2.csv};
		\addlegendentryexpanded{$\dot{\bar{u}} = 10^{1}\,\sfrac{\mathrm{mm}}{\mathrm{s}}$};

		\addplot+[cTab1, mark = none, densely dashed, thick] table [
						x = u5,
						y = 0_555_E-01,
						col sep = comma] {reaction_perforated_plate_rate0.01_referenceR2.csv};
		\addlegendentryexpanded{$\dot{\bar{u}} = 10^{0}\,\sfrac{\mathrm{mm}}{\mathrm{s}}$};

		\addplot+[mark = none, black, thick] table [
						x = u,
						y = f,
						col sep = comma] {reaction_perforated_plate_rateindependent_referenceR2.csv};
		\addlegendentryexpanded{$\dot{\bar{u}} \rightarrow 0\,\sfrac{\mathrm{mm}}{\mathrm{s}}$};

		\addlegendimage{legend image with text = {$\dot{\bar{u}} =$}}
		\addlegendentryexpanded{$\mspace{-2mu}10^{0}\,\sfrac{\mathrm{mm}}{\mathrm{s}}$}

		\addplot+[solid, cTab20c3, mark = square*, mark options = {solid, fill = cTab20c4!50!white}] table [x = u,
										y = f,
										col sep = comma] {reaction_perforated_plate_rate0.01_mu0.55E-1_k0R2.csv};
		\addlegendentryexpanded{cycle $k = 0$};

		\addplot+[solid, cTab20c2, mark = diamond*, mark options = {solid, fill = cTab20c3!50!white}] table [x = u,
										y = f,
										col sep = comma] {reaction_perforated_plate_rate0.01_mu0.55E-1_k1R2.csv};
		\addlegendentryexpanded{cycle $k = 1$};

		\addplot+[solid, cTab20c1, mark = *, mark options = {solid, fill = cTab20c2!50!white}] table [x = u,
										y = f,
										col sep = comma] {reaction_perforated_plate_rate0.01_mu0.55E-1_k2R2.csv};
		\addlegendentryexpanded{cycle $k = 2$};

		\addlegendimage{legend image with text = {\,}}
		\addlegendentryexpanded{\,}

		\addlegendimage{legend image with text = {$\dot{\bar{u}} = $}}
		\addlegendentryexpanded{$\mspace{-2mu}10^1\,\sfrac{\mathrm{mm}}{\mathrm{s}}$}

		\addplot+[solid, cTab20c7, mark = square*, mark options = {solid, fill = cTab20c8!50!white}] table [x = u,
										y = f,
										col sep = comma] {reaction_perforated_plate_rate0.01_mu0.55E-0_k0R2.csv};
		\addlegendentryexpanded{cycle $k = 0$};

		\addplot+[solid, cTab20c6, mark = diamond*, mark options = {solid, fill = cTab20c7!50!white}] table [x = u,
										y = f,
										col sep = comma] {reaction_perforated_plate_rate0.01_mu0.55E-0_k1R2.csv};
		\addlegendentryexpanded{cycle $k = 1$};

		\addplot+[solid, cTab20c5, mark = *, mark options = {solid, fill = cTab20c6!50!white}] table [x = u,
										y = f,
										col sep = comma] {reaction_perforated_plate_rate0.01_mu0.55E-0_k2R2.csv};
		\addlegendentryexpanded{cycle $k = 2$};

		\addlegendimage{legend image with text = {\,}}
		\addlegendentryexpanded{\,}

		\addlegendimage{legend image with text = {\hfill $\dot{\bar{u}} = $}}
		\addlegendentryexpanded{$\mspace{-2mu}10^3\,\sfrac{\mathrm{mm}}{\mathrm{s}}$}

		\addplot+[solid, cTab20c11, mark = square*, mark options = {solid, fill = cTab20c12!50!white}] table [x = u,
										y = f,
										col sep = comma] {reaction_perforated_plate_rate0.01_mu0.55E+2_k0R2.csv};
		\addlegendentryexpanded{cycle $k = 0$};

		\addplot+[solid, cTab20c10, mark = diamond*, mark options = {solid, fill = cTab20c11!50!white}] table [x = u,
										y = f,
										col sep = comma] {reaction_perforated_plate_rate0.01_mu0.55E+2_k1R2.csv};
		\addlegendentryexpanded{cycle $k = 1$};

		\addplot+[solid, cTab20c9, mark = *, mark options = {solid, fill = cTab20c10!50!white}] table [x = u,
										y = f,
										col sep = comma] {reaction_perforated_plate_rate0.01_mu0.55E+2_k2R2.csv};
		\addlegendentryexpanded{cycle $k = 2$};

		\coordinate (spypoint1) at (axis cs:0.582, 8.175);
		\coordinate (magnifyglass1) at (axis cs:0.443, 9.5);

		\coordinate (spypoint2) at (axis cs:0.582, 5.19);
		\coordinate (magnifyglass2) at (axis cs:0.443, 6.5);

		\coordinate (spypoint3) at (axis cs:0.582, 4.15);
		\coordinate (magnifyglass3) at (axis cs:0.443, 0.5);

		\coordinate (l1) at (axis cs:0.2, 2.2);
		\coordinate (l2) at (axis cs:0.275, 3.6);
		\coordinate (l3) at (axis cs:0.40, 4.7);
		\coordinate (l4) at (axis cs:0.50, 7.7);

	\end{axis}

	\spy[width = 3cm, height = 1.5cm, magnification=5.5] on (spypoint1) in node[fill = white] at (magnifyglass1);
	\spy[width = 3cm, height = 1.5cm, magnification=5.5] on (spypoint2) in node[fill = white] at (magnifyglass2);
	\spy[width = 3cm, height = 1.5cm, magnification=5.5] on (spypoint3) in node[fill = white] at (magnifyglass3);

	\node[font = \scriptsize] at (l1) {(i)~$\dot{\overline{u}} = 0\,\mathrm{mm}/\mathrm{s}$};
	\node[font = \scriptsize] at (l2) {(ii)~$\dot{\overline{u}} = 10^{0}\,\mathrm{mm}/\mathrm{s}$};
	\node[font = \scriptsize] at (l3) {(iii)~$\dot{\overline{u}} = 10^{1}\,\mathrm{mm}/\mathrm{s}$};
	\node[font = \scriptsize] at (l4) {(iv)~$\dot{\overline{u}} = 10^{3}\,\mathrm{mm}/\mathrm{s}$};

	\node[text width = 1.6cm] at (1.2, 9.75) {\footnotesize$m = 0.1$\\$\mu = 500\,\mathrm{s}$};

	\begin{scope}[shift = {(1.5, 7.5)}, isometric view, scale = 0.2, line join=round]
		\fill[cSecond, opacity = 0.2] (0, 9, 0) -- ++(5, 0, 0) -- ++(0, 0, 0.5) -- ++(-5, 0, 0) -- cycle;
		\draw[cSecond] (0, 9, 0) -- ++(5, 0, 0) -- ++(0, 0, 0.5) -- ++(-5, 0, 0) -- cycle;

		\node[above, font = \footnotesize, fill = white, fill opacity = 0.8, text opacity = 1, inner sep = 1, outer sep = 2] at (2.5, 9, 0.6) {$\overline{u}$};

		\foreach \x in {0, 1, ..., 5}
			{
				\foreach \z in {0, 0.5}
					{
						\draw[cSecond, -stealth] (\x, 8, \z) -- ++(0, 1, 0);
					}
			}

		\draw[black] (5, 8, 0) -- ++(0, 0, 0.5);
		\begin{scope}[canvas is xy plane at z = 0]
			\draw[black] (2.5, 0) -- ++(2.5, 0) -- ++(0, 8) -- ++(-5, 0) -- ++(0, -5.5) arc(90:0:2.5) -- cycle;
		\end{scope}

		\begin{scope}[canvas is xy plane at z = 0.25]
			\draw[
				black,
				top color = bulk,
				bottom color = bulk,
				middle color = white,
				shading angle = 15
			] (2.5, 0) -- node[pos = 0.5, anchor = south east, inner sep = 1.2, cFirst, opacity = .65] {} ++(2.5, 0) -- ++(0, 8) -- ++(-5, 0) -- node[pos = 0.5, anchor = south west, inner sep = 2, cFirst, opacity = .65] {} ++(0, -5.5) arc(90:0:2.5) -- cycle;
		\end{scope}

		\begin{scope}[canvas is xy plane at z = 0.5]
			\draw[black] (2.5, 0) -- ++(2.5, 0) -- ++(0, 8) -- ++(-5, 0) -- ++(0, -5.5) arc(90:0:2.5) -- cycle;
		\end{scope}

		\draw[black] (2.5, 0, 0) -- ++(0, 0, 0.5);
		\draw[black] (0, 2.5, 0) -- ++(0, 0, 0.5);
		\draw[black] (5, 0, 0) -- ++(0, 0, 0.5);
		\draw[black] (0, 8, 0) -- ++(0, 0, 0.5);

		\begin{scope}[canvas is xz plane at y = 0]
			\fill[cFirst, opacity = .0] (2.5, 0) rectangle (5, 0.5);
		\end{scope}
		\begin{scope}[canvas is yz plane at x = 0]
			\fill[cFirst, opacity = .0] (8, 0) rectangle (2.5, 0.5);
		\end{scope}

		\draw[black, densely dashed] (0, 0, 0) -- ++(2.5, 0, 0);
		\draw[black, densely dashed] (0, 0, 0) -- ++(0, 2.5, 0);
		\draw[black, densely dashed] (0, 0, 0.5) -- ++(2.5, 0, 0);
		\draw[black, densely dashed] (0, 0, 0.5) -- ++(0, 2.5, 0);

	\end{scope}

\end{tikzpicture} \\
		\end{tabular}
		\caption{Plate with a circular hole: viscoplastic response for various strain rates. Dashed lines refer to FEM reference solutions from literature \cite{deSouzaNeto:2011}. All other lines refer to the numerical analysis results for which the curves of different refinement levels coincide and only the last refinement result indicated by the circular markers are identifiable. }
		\label{PPVP}
\end{figure}

In Fig. \ref{PPVP} the load-displacement curves of the analysis for a viscoplastic exponent \(m=1.0\) (left) and \(m=0.1\) (right) are depicted. The black solid line (i) refers to the rate-independent analysis result, cf. Fig. \ref{PPIntegrationConvergence}. The different coloring of the curves refers to different deformation rates. A reference FEM result from a plane stress analysis is shown for each curve as dashed line \cite{deSouzaNeto:2011}. All computation results include curves for the various refinement levels, of which only the most refined result is identifiable as the curve above the other results. The computational results are in very good agreement with the reference solutions but behave slightly weaker which can be attributed to the quadratic approximation space used in our models compared to the linear functions (CST element) used in the FEM reference model. Both models of Fig. \ref{PPVP} were loaded with a total prescribed displacement of \(\bar{u}=0.6\,\mathrm{mm}\) in $43$ equal load increments. In both cases, a similar iteration behavior was observed with an almost constant number of iterations per load step which set in after $5-8$ initial steps. In good agreement with the error analysis in \cite{deSouzaNeto:2011}, the number of iterations per step decreased with increasing deformation rate. The average numbers of iterations for the model with \(m=1.0\) were: (ii) 14 iterations/load step, (iii) 7 iterations/load step and (iv) 4 iterations/load step. The model with \(m=0.1\) experienced an iteration effort of: (ii) 10 iterations/load step, (iii) 8 iterations/load step and (iv) 7 iterations/load step.

\subsection{Porous plate}
\label{sec:porousPlate}
The performance of the proposed analysis method, including error-indicator-controlled overlay refinement and non-negative moment-fitted quadrature, is assessed for thermo-viscoplastic problems in terms of numerical stability, reliability, and computational effort. To this end, we consider a porous plate subjected to thermal loading and prescribed displacements, cf. Fig. \ref{PourousPlateModel}. A detailed description of the geometry is provided in Appendix~E. The plate was subjected to a total prescribed displacement of $\bar{u}_y=0.09\,\mathrm{mm}$, applied in $30$ equal load steps. A viscosity exponent of $m=1$, a relaxation time of $\mu=0.001\,\mathrm{s}$, and a constant deformation rate of $\dot{\bar{u}}_y=1\,\mathrm{mm/s}$ were used throughout all computations for this example. In addition to the results obtained with the proposed framework, we report a reference solution computed using the open-source library \emph{deal.II} \cite{Arndt:2021} with identical material and loading parameters. A highly resolved reference solution was obtained using \(21\, 948\) quadratic Lagrange elements with \(3\times3\times3\) Gauss quadrature.

\begin{figure}[h!]
	\centering 
	\begin{tikzpicture}
		\node (A) {\tikzsetnextfilename{model_porus_plate}
\begin{tikzpicture}[scale = 0.75, line join=round]

  \begin{scope}[isometric view]
		\fill[blue, opacity = 0.2] (0, 9, 0) -- ++(8, 0, 0) -- ++(0, 0, 0.5) -- ++(-8, 0, 0) -- cycle;
		\draw[blue] (0, 9, 0) -- ++(8, 0, 0) -- ++(0, 0, 0.5) -- ++(-8, 0, 0) -- cycle;

		\node[blue, below, font = \footnotesize] at (4, -1, 0.6) {$u_y=0$};
		\node[blue, below, font = \footnotesize] at (-1, 4, 0.6) {$u_x=0$};
		\node[blue, above, font = \footnotesize] at (4, 9.5, 0.6) {$u_y=\bar{u}_y$};

		\foreach \x in {0, 1, ..., 8}
			{
				\foreach \z in {0, 0.5}
					{
						\draw[blue, -latex, very thin] (\x, 8, \z) -- ++(0, 1, 0);
					}
			}

		\begin{scope}[canvas is yz plane at x = 0.0]
			\fill[
				top color = lightgray!40!white,
				bottom color = lightgray!80!black,
				shading angle = 90
			]
			 (0, 0) -- ++(0, 0.5) -- ++(8, 0) -- ++(0, -.5) -- cycle;
		\end{scope}

		\begin{scope}[canvas is xy plane at z = 0.5]
			\draw[
				top color = lightgray,
				bottom color = lightgray,
				middle color = white,
				shading angle = 15
			]
			 (0, 0) -- ++(8, 0) -- ++(0, 8) -- ++(-8, 0) -- cycle;
		\end{scope}

		\draw[] (0, 0, 0.5) -- (0, 8, 0.5) -- (0, 8, 0) -- (0, 0, 0) -- (8, 0, 0) -- (8, 0, 0.5) --  cycle;

		\begin{scope}[canvas is xy plane at z = 0.5]
			\begin{axis}[xshift = 2pt, yshift = 0pt, hide axis, width = 9.5cm, height = 9.5cm, xmin=0, xmax=10, ymin=0, ymax=10, enlarge x limits=0, enlarge y limits=0]
				\addplot[scatter, black!70, only marks, mark=*, mark options = {solid, fill = black!70},
						point meta=explicit symbolic,
						scatter/@pre marker code/.style={/tikz/mark size=\pgfplotspointmeta*0.9 cm},
						scatter/@post marker code/.style={}]
					table [x = x, y = y, meta = r, col sep = comma,
						x filter/.code={\pgfmathparse{\pgfmathresult * 1.25}},
						y filter/.code={\pgfmathparse{\pgfmathresult * 1.25}}]
					{model_porous_plate_holesR2.csv}; 
			\end{axis}
		\end{scope}

		\draw[] (0, 0, 0) -- ++(0, 0, 0.5);
		\draw[] (8, 0, 0) -- ++(0, 0, 0.5);
		\draw[] (0, 8, 0) -- ++(0, 0, 0.5);

		\begin{scope}[canvas is xz plane at y = 0]
			\fill[blue, opacity = .5] (0, 0) rectangle (8, 0.5);
			\draw[pattern color = blue, pattern = north west lines] (0, 0) rectangle (8, 0.5);
		\end{scope}

		\begin{scope}[canvas is yz plane at x = 0]
			\fill[blue, opacity = .5] (0, 0) rectangle (8, 0.5);
			\draw[pattern color = blue, pattern = north east lines] (0, 0) rectangle (8, 0.5);
		\end{scope}

		\begin{scope}[canvas is yz plane at x = 0]
			\draw[red, fill=red, fill opacity = 0.15] (0, 0.5) -- (0, 1.2) sin (4, 4.5) cos (8, 1.2) -- (8, 0.5);
			\node[red, above] at (2.5, 4.6) {$\bar{T}_l(y)$};
			\fill[black!70] (3.25, 0) rectangle ++(.27, .5);
		\end{scope}

		\begin{scope}[canvas is yz plane at x = 8]
			\draw[red, fill=red, fill opacity = 0.15] (0, 0.5) -- (0, 1.2) sin (4, 2.8) cos (8, 1.2) -- (8, 0.5);
			\node[red, anchor=south west, inner sep = 0] at (3, 2.6) {$\bar{T}_r(y)$};
		\end{scope}

		\draw[-latex] (0,0,0) -- (1.5, 0, 0) node[below right, inner sep = 1] {$x$};
		\draw[-latex] (0,0,0) -- (0, 1.5, 0) node[below left, inner sep = 1] {$y$};
		\draw[-latex] (0,0,0) -- (0, 0, 1.5) node[right] {$z$};

  \end{scope}


\end{tikzpicture}};
		\node[right, shift = {(10mm, -6mm)}] at (A.east) {
			\begin{minipage}{5cm}
				\scriptsize%
				\textbf{Boundary conditions} \\[1mm]
				essential boundary conditions \\[-6mm]
				\begin{equationarray*}{rcl}
            			u_x &=& 0\quad\text{along}\quad(0,y,z)\\[1mm]
            			u_y &=& 0\quad\text{along}\quad(x,0,z)\\[1mm]
            			u_z &=& 0\quad\text{along}\quad(x,y,0)\\[-5mm]
				\end{equationarray*}
				prescribed displacement \\[-6mm]
				\begin{equationarray*}{rcl} 
					u_y & =& \overline{u}_y\quad\text{along}\quad(x,10,z)\\[-5mm]
				\end{equationarray*}
				prescribed temperatures \\[-6mm]
				\begin{equationarray*}{rcl}
						T &=& \overline{T}_l(y) \quad\text{along}\quad(0,y,z) \\[1mm]
						T &=& \overline{T}_r(y) \quad\text{along}\quad(10,y,z) \\[3mm]
 
					 \overline{T}_l(y) &=& 500 + 20\,\sin(\pi y/10) \, [^\circ\mathrm{C}]\\[1mm]
					 \overline{T}_r(y) &=& 500 + 5\,\sin(\pi y/10) \, [^\circ\mathrm{C}]\\
				\end{equationarray*}
				
			\end{minipage}
		};

		\node[above right] at (A.south west) {\footnotesize $10.0 \times 10.0 \times 0.75\ [\mathrm{mm}]$};
	\end{tikzpicture}
	\vspace{-6mm}
	\caption{Porous plate model -- geometry and boundary conditions. The plate is constrained according to the symmetry/support conditions indicated in the figure and subjected to a prescribed displacement \(\bar u_y\) at the upper boundary face.
	 Different sinusoidal temperature distributions \(\bar{T}_{\{l,r\}}(y)\) are prescribed along the left and right boundary faces. The base mesh consists of \((32\times32\times1)\) elements and is successively refined by multiple overlay levels.}	
	\label{PourousPlateModel}
\end{figure}

\begin{figure}[h!]
	\centering 
	\includegraphics[width=0.5\textwidth]{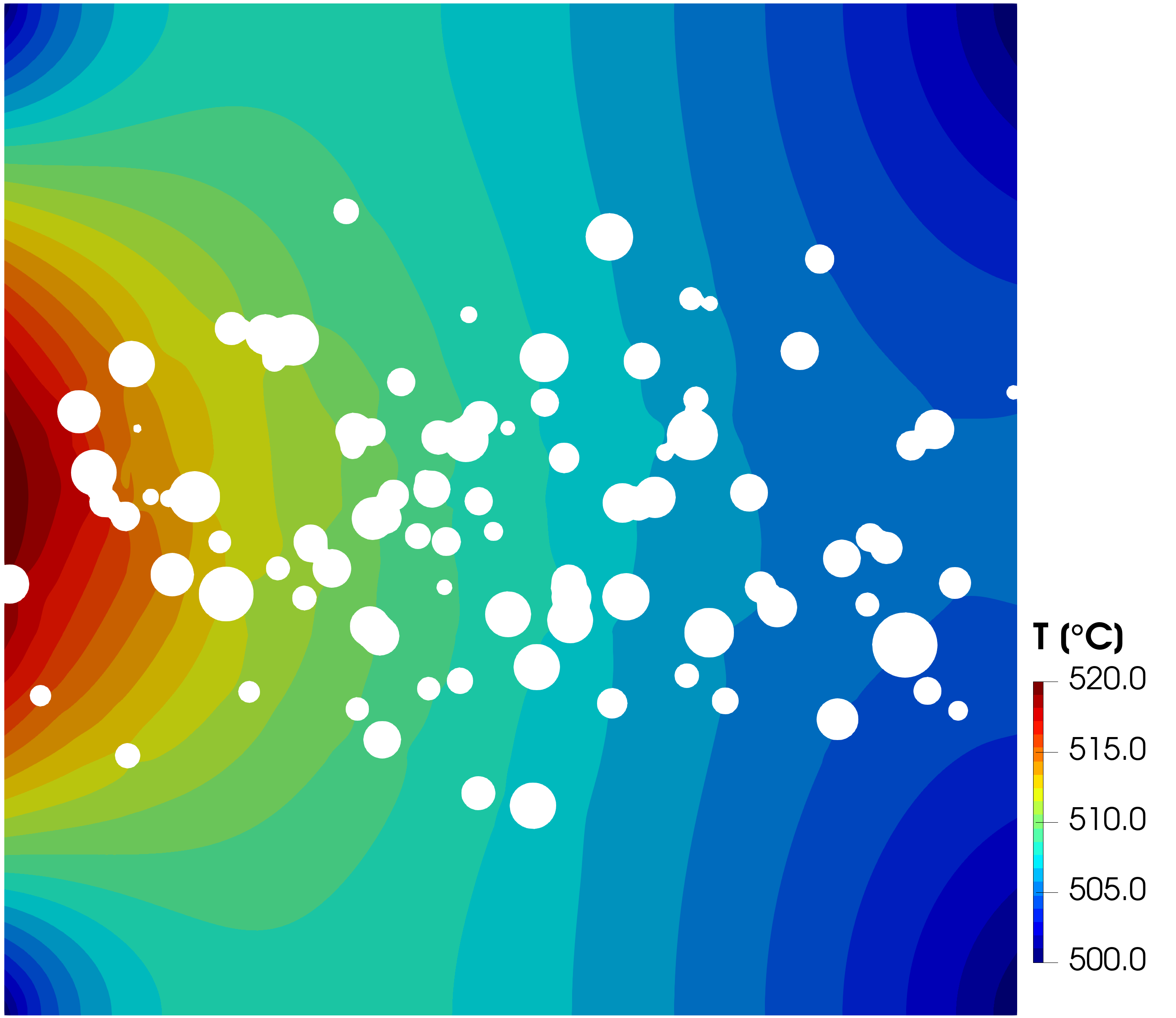}
	\caption{Porous plate: applied temperature profile.}
	\label{PorousPlateTemp}
\end{figure}
\paragraph{Thermo-viscoplasticity}
%
The thermal model of the example considers an isotropic thermal expansion corresponding to a sinusoidally changing temperature load along two opposite interfaces, but does not consider the feedback of heat from plastic deformation, assuming that this effect is significantly lower than the temperature load shown in Fig. \ref{PorousPlateTemp}. The thermal strain tensor is adapted to the temperature dependence of the thermal expansion coefficient as:
\begin{equationarray}{rcl}
	\boldsymbol{\varepsilon}^{th} &=& \color{mygreenD}\left[\gamma_m(T) [T-T_{\mathrm{ref}}]-\gamma_m(T_r) [T_r-T_{\mathrm{ref}}]\right]\,\mathbf{I}
\end{equationarray} 
where \(\gamma_m(T)\) is the mean thermal expansion coefficient and the temperatures \(T_{\mathrm{ref}}=273.15\,\mathrm{K}\) and \(T_r=293.15\,\mathrm{K}\) mark the origin and a reference temperature of the Celsius scale, respectively. Due to the discontinuity of the measured data at the austenitisation temperature, cf. Fig. \ref{ConstitutivePropertiesII}, Appendix C, a piecewise interpolation is used for the mean thermal expansion coefficient \(\gamma_m(T)\) using two linear polynomials and for the heat conductivity \(\kappa(T)\) using a fourth order polynomial and an exponential function, cf. Appendix D.

\paragraph{Non-negative moment fitting}
We compared the established \emph{Adaptive Space-Tree} (AST) approach with the \emph{Non-Negative Moment Fitting} (NNMF) method of Section~\ref{sec:NNMFquadrature}. For the AST method, four recursive bisection levels were employed in the \(x\)-\(y\) plane, whereas the through-thickness direction remained unrefined. For each sub-cell, \((p+1)^3\) Gauss--Legendre quadrature points were used to evaluate the integrals with sufficient accuracy, resulting in almost \(6\) million quadrature points for the entire domain. 
\color{mygreenD}
In comparison, the NNMF approach was applied to the cell model obtained after four adaptive refinement cycles, at which the prescribed maximum overlay refinement depth of $k=4$ had been reached. This model required \(213\,309\) quadrature points, corresponding to a reduction by almost a factor of \(28\) per iteration step. Subsequent adaptive cycles further adapt the locally refined regions while retaining the maximum refinement depth of \(k=4\).
\color{black}

The corresponding quadrature point distributions are shown in Fig.~\ref{PorousPlateGP}, illustrating the substantial reduction achieved with NNMF.

\begin{figure}[htb!]
	\centering 
	\resizebox{\textwidth}{!}{\input{quad_points_porous_plateR3.tex}}
	\caption{Porous plate: Porous plate: quadrature point distributions obtained with AST and NNMF for the discretization \color{mygreenD}after four adaptive refinement cycles, at which the prescribed maximum overlay refinement depth of $k=4$ had been reached.} 
	\label{PorousPlateGP}
\end{figure}

The computational effort associated with the numerical integration schemes, AST and NNMF, is assessed in Fig.~\ref{PorousPlateRuntime} in terms of the runtime contributions of the different stages of the analysis. The left figure shows that the computational cost of the full NNMF analysis increases moderately over the adaptive cycles, from approximately \(80\,\mathrm{s}\) to \(280\,\mathrm{s}\), with assembly and the linear solve dominating the computational effort as the model size increases. In comparison, the time required for the Kelly error estimation and overlay refinement remains almost negligible. In contrast, the space-tree quadrature results shown in the right figure exhibit substantially higher computational costs. For the first four adaptive cycles, NNMF achieves a consistently high speedup of approximately \(20\times\), requiring only \(84\)--\(116\,\mathrm{s}\) compared with \(1669\)--\(2331\,\mathrm{s}\) for the corresponding AST computations.

Moreover, the significantly smaller number of quadrature points required by NNMF allows the quadrature information to be generated and stored once, avoiding repeated setup during the incremental analysis and thereby further reducing the computational effort.
 
\begin{figure}[htb!]
	\centering 
	\pgfplotsset{
  paired/.style={
    width=15.0cm, height=9.2cm,
    ybar stacked, bar width=0.50cm,
    xmin=-0.62, xmax=3.62, ymin=0, ymax=2700,
    xtick={0,1,2,3},
    axis line style={oiGrey, line width=0.5pt},
    tick style={oiGrey, line width=0.5pt},
    label style={font=\small}, tick label style={font=\footnotesize},
    grid style={black!12, line width=0.35pt},
  },
}
\newcommand{\segs}[9]{%
  \addplot[fill=oiGreen,      draw=oiGrey, line width=0.25pt, bar shift=#1] coordinates {#2};
  \addplot[fill=oiBlue,       draw=oiGrey, line width=0.25pt, bar shift=#1] coordinates {#3};
  \addplot[fill=oiSky,        draw=oiGrey, line width=0.25pt, bar shift=#1] coordinates {#4};
  \addplot[fill=oiPurple,     draw=oiGrey, line width=0.25pt, bar shift=#1] coordinates {#5};
  \addplot[fill=oiOrange,     draw=oiGrey, line width=0.25pt, bar shift=#1] coordinates {#6};
  \addplot[fill=oiYellow,     draw=oiGrey, line width=0.25pt, bar shift=#1] coordinates {#7};
  \addplot[fill=oiVermillion, draw=oiGrey, line width=0.25pt, bar shift=#1] coordinates {#8};
  \addplot[fill=oiLightGrey,  draw=oiGrey, line width=0.25pt, bar shift=#1] coordinates {#9};}

\begin{tikzpicture}
\begin{axis}[
  width=0.475\textwidth, height=9.0cm,
  ybar stacked,
  bar width=0.3cm,
  title = {Full NNMF run},
  xlabel={Kelly adaptive cycle},
  ylabel={runtime [\si{\second}]},
  xmin=-0.7, xmax=12.7, ymin=0, ymax=318,
  xtick={0,1,2,3,4,5,6,7,8,9,10,11,12},
  ytick={0,50,100,150,200,250,300},
  ymajorgrids=true,
  legend cell align=left,
  legend columns=4,
  legend style={
    at={(-0.1cm, 10.15cm)}, anchor=north west,
    draw=oiGrey, line width=0.4pt, fill=white,
    font=\footnotesize, row sep=1.6pt, inner sep=4.5pt,
    /tikz/every even column/.append style={column sep=8pt},
    /tikz/every odd column/.append style={column sep=3pt},
  },
]

\addplot[fill=oiGreen, draw=oiGrey, line width=0.25pt] coordinates {
  (0,11.51) (1,11.30) (2,11.48) (3,12.18) (4,11.56) (5,11.32) (6,11.76)
  (7,11.93) (8,12.56) (9,12.43) (10,12.84) (11,13.94) (12,14.45)
};
\addlegendentry{quadrature setup}

\addplot[fill=oiBlue, draw=oiGrey, line width=0.25pt] coordinates {
  (0,29.89) (1,31.98) (2,37.46) (3,42.43) (4,48.12) (5,53.99) (6,60.71)
  (7,67.04) (8,89.38) (9,95.94) (10,113.99) (11,133.61) (12,158.14)
};
\addlegendentry{assembly}

\addplot[fill=oiSky, draw=oiGrey, line width=0.25pt] coordinates {
  (0,22.10) (1,21.09) (2,25.19) (3,25.50) (4,27.95) (5,28.11) (6,30.57)
  (7,32.58) (8,40.07) (9,42.48) (10,47.07) (11,53.56) (12,58.27)
};
\addlegendentry{linear solve}

\addplot[fill=oiPurple, draw=oiGrey, line width=0.25pt] coordinates {
  (0,8.77) (1,9.39) (2,10.91) (3,11.65) (4,12.64) (5,13.65) (6,15.05)
  (7,15.95) (8,20.68) (9,22.06) (10,25.33) (11,28.84) (12,33.24)
};
\addlegendentry{material (return map)}

\addplot[fill=oiOrange, draw=oiGrey, line width=0.25pt] coordinates {
  (0,4.31) (1,2.15) (2,3.03) (3,5.99) (4,4.20) (5,4.61) (6,5.21)
  (7,5.23) (8,5.74) (9,6.30) (10,7.25) (11,8.39) (12,9.39)
};
\addlegendentry{setup}

\addplot[fill=oiYellow, draw=oiGrey, line width=0.25pt] coordinates {
  (0,1.90) (1,1.48) (2,2.02) (3,3.29) (4,2.09) (5,2.17) (6,2.25)
  (7,2.26) (8,2.40) (9,2.49) (10,2.57) (11,2.72) (12,0.00)
};
\addlegendentry{Kelly estimate $+$ refinement}

\addplot[fill=oiVermillion, draw=oiGrey, line width=0.25pt] coordinates {
  (0,0.94) (1,1.21) (2,2.46) (3,7.52) (4,7.53) (5,7.77) (6,7.88)
  (7,7.83) (8,8.03) (9,8.36) (10,8.46) (11,8.77) (12,9.01)
};
\addlegendentry{Dirichlet projection}

\addplot[fill=oiLightGrey, draw=oiGrey, line width=0.25pt] coordinates {
  (0,1.62) (1,0.82) (2,0.89) (3,0.93) (4,0.09) (5,0.03) (6,0.04)
  (7,0.11) (8,0.04) (9,0.04) (10,0.04) (11,0.08) (12,0.08)
};
\addlegendentry{unaccounted}

\end{axis}

\begin{scope}[xshift = 0.5\textwidth]
\begin{axis}[paired, name=main,
  title = {NNMF (left) vs. Space-Tree (right)},
  xlabel={Kelly adaptive cycle},
  ylabel={runtime [\si{\second}]},
  ytick={0,400,800,1200,1600,2000,2400},
  ymajorgrids=true,
  width = 0.5\textwidth,
  height = 9.0cm
]
\segs{-0.28cm}
 {(0,11.40) (1,11.32) (2,11.54) (3,12.31)}
 {(0,33.33) (1,35.42) (2,42.85) (3,48.93)}
 {(0,22.41) (1,21.24) (2,24.96) (3,26.79)}
 {(0,9.82)  (1,10.76) (2,12.41) (3,13.66)}
 {(0,3.91)  (1,2.08)  (2,3.00) (3,5.93)}
 {(0,1.91)  (1,1.47)  (2,2.05) (3,0.00)}
 {(0,0.96)  (1,1.24)  (2,2.44) (3,7.57)}
 {(0,1.59)  (1,0.82)  (2,0.95) (3,0.93)}
\end{axis}

\begin{axis}[paired, at={(main.south west)}, anchor=south west,
  axis lines=none, xtick=\empty, ytick=\empty,
  width = 0.5\textwidth,
  height = 9.0cm
]
\segs{+0.28cm}
 {(0,2.58)    (1,2.68)    (2,2.87) (3,3.27)}
 {(0,1169.36) (1,1206.88) (2,1448.10) (3,1677.89)}
 {(0,20.02)   (1,19.23)   (2,21.48) (3,23.87)}
 {(0,464.27)  (1,477.11)  (2,551.75) (3,606.61)}
 {(0,7.83)    (1,6.09)    (2,7.35) (3,10.58)}
 {(0,1.91)    (1,1.50)    (2,2.12) (3,0.00)}
 {(0,0.95)    (1,1.19)    (2,2.43) (3,7.38)}
 {(0,2.25)    (1,1.51)    (2,1.59) (3,1.65)}

\node[font=\scriptsize, text=oiGrey, anchor=south] at (axis cs:0.145,1690) {\SI{1669}{\second}};
\node[font=\scriptsize, text=oiGrey, anchor=south] at (axis cs:1.145,1737) {\SI{1716}{\second}};
\node[font=\scriptsize, text=oiGrey, anchor=south] at (axis cs:2.145,2059) {\SI{2038}{\second}};
\node[font=\scriptsize, text=oiGreen!60!black, anchor=south, align=center]
  at (axis cs:-0.245,120) {\SI{85}{\second}};
\node[font=\scriptsize, text=oiGreen!60!black, anchor=south, align=center]
  at (axis cs:0.755,120) {\SI{84}{\second}};
\node[font=\scriptsize, text=oiGreen!60!black, anchor=south, align=center]
  at (axis cs:1.755,140) {\SI{100}{\second}};
\node[font=\footnotesize, text=oiBlue, anchor=south, align=center, inner sep=2pt]
  at (axis cs:0.145,1880) {$\mathbf{19.6\times}$};
\node[font=\footnotesize, text=oiBlue, anchor=south, align=center, inner sep=2pt]
  at (axis cs:1.145,1927) {$\mathbf{20.3\times}$};
\node[font=\footnotesize, text=oiBlue, anchor=south, align=center, inner sep=2pt]
  at (axis cs:2.145,2249) {$\mathbf{20.3\times}$};
\node[font=\scriptsize, text=oiGrey, anchor=south] at (axis cs:3.145,2352) {\SI{2331}{\second}};
\node[font=\footnotesize, text=oiBlue, anchor=south, align=center, inner sep=2pt]
  at (axis cs:3.145,2530) {$\mathbf{20.1\times}$};
\node[font=\scriptsize, text=oiGreen!60!black, anchor=south, align=center]
  at (axis cs:2.755,155) {\SI{116}{\second}};

\end{axis}
\end{scope}
\end{tikzpicture}
	\caption{Porous plate: decomposition of the total runtime into the individual stages of the full NNMF analysis over all adaptive cycles (left) and corresponding runtime comparison between NNMF and AST for the first four adaptive cycles (right). The reported speedup refers to the total runtime of NNMF relative to AST.}
	\label{PorousPlateRuntime}
\end{figure}

\begin{figure}[htb!]
	\centering 
	\begin{tabular}{@{} c @{\space} c @{\space} c @{}}
		\includegraphics[width = 0.32\textwidth]{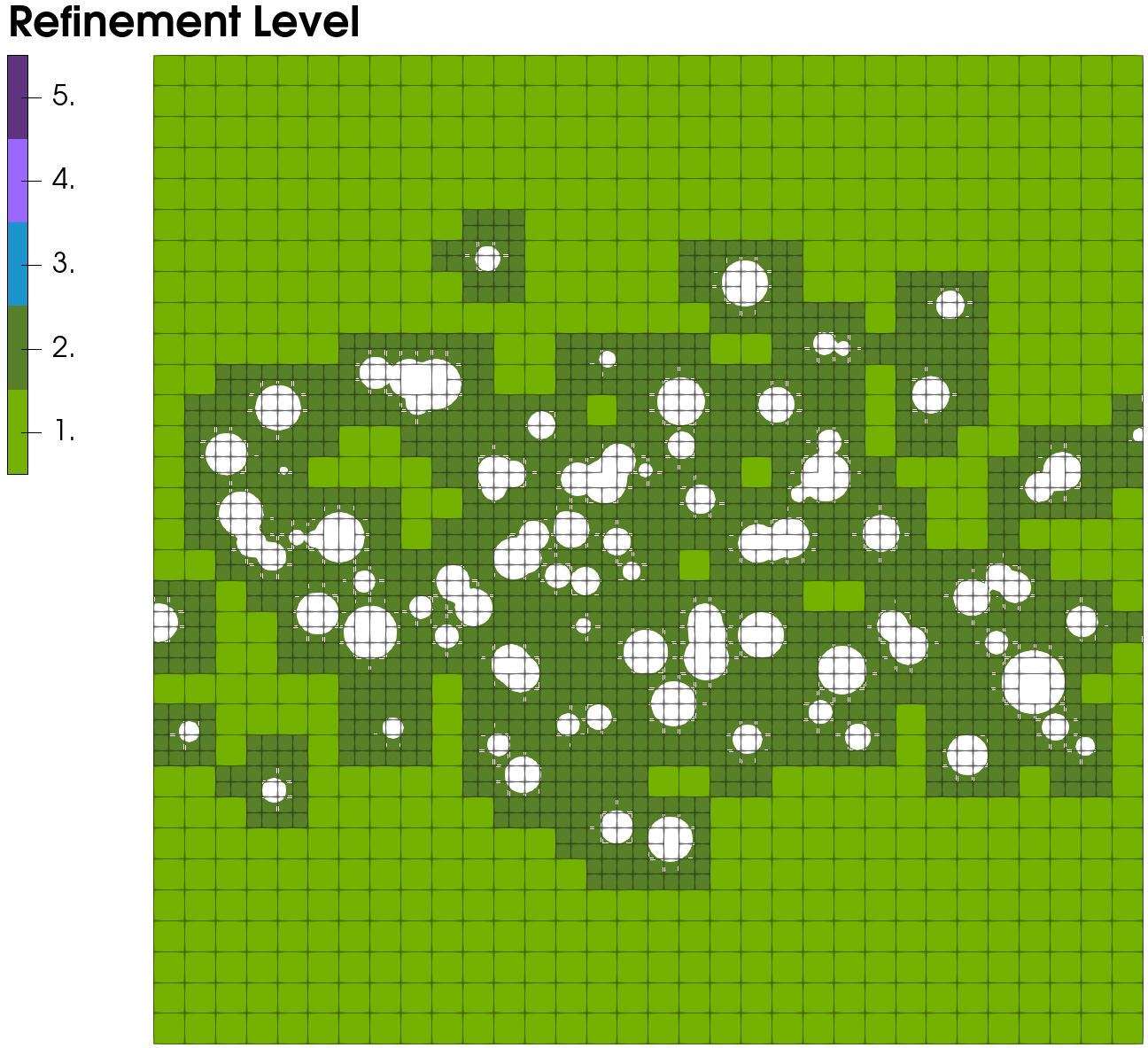} &
		\includegraphics[width = 0.32\textwidth]{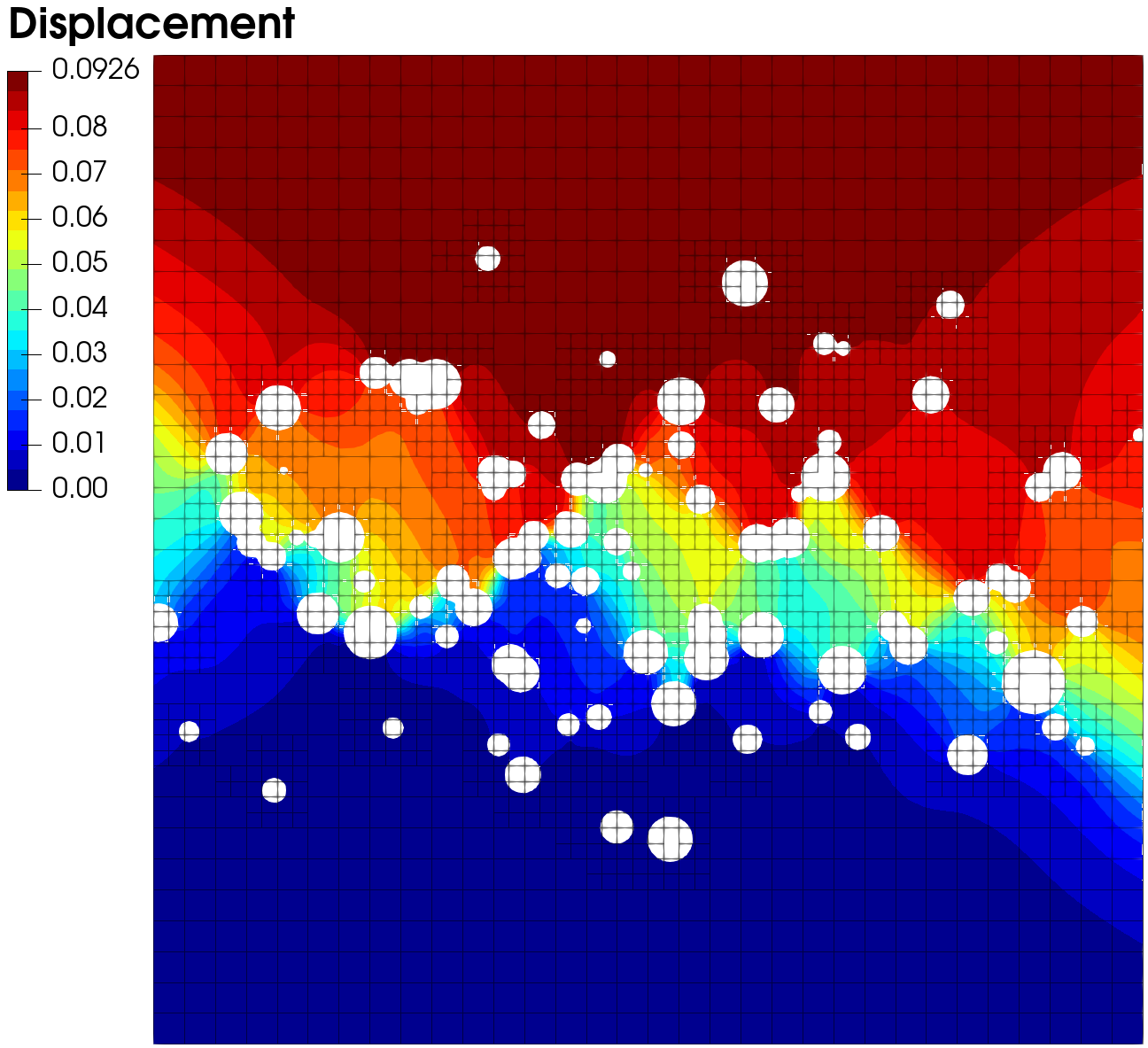} &
		\includegraphics[width = 0.32\textwidth]{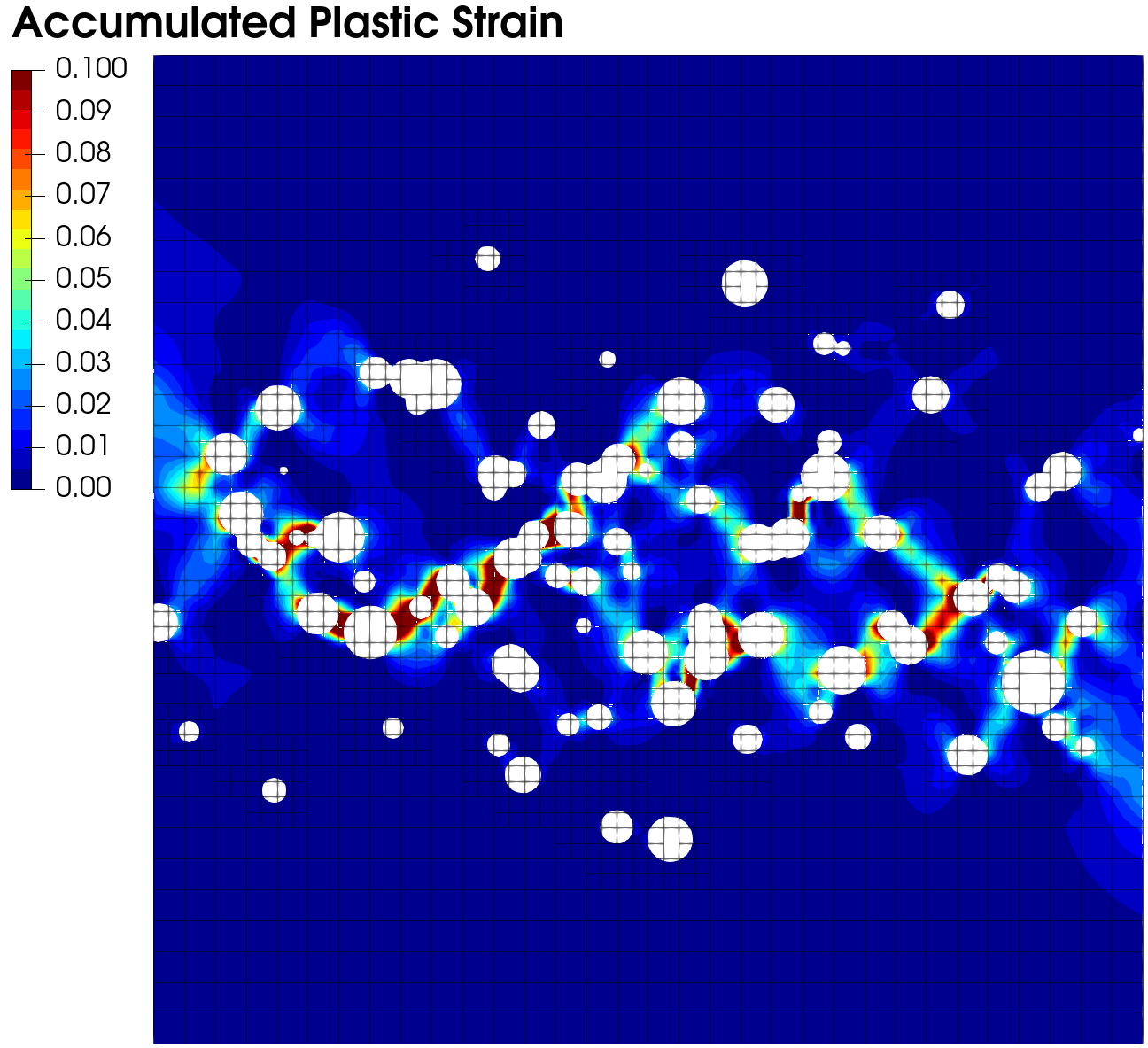} \\
		\includegraphics[width = 0.32\textwidth]{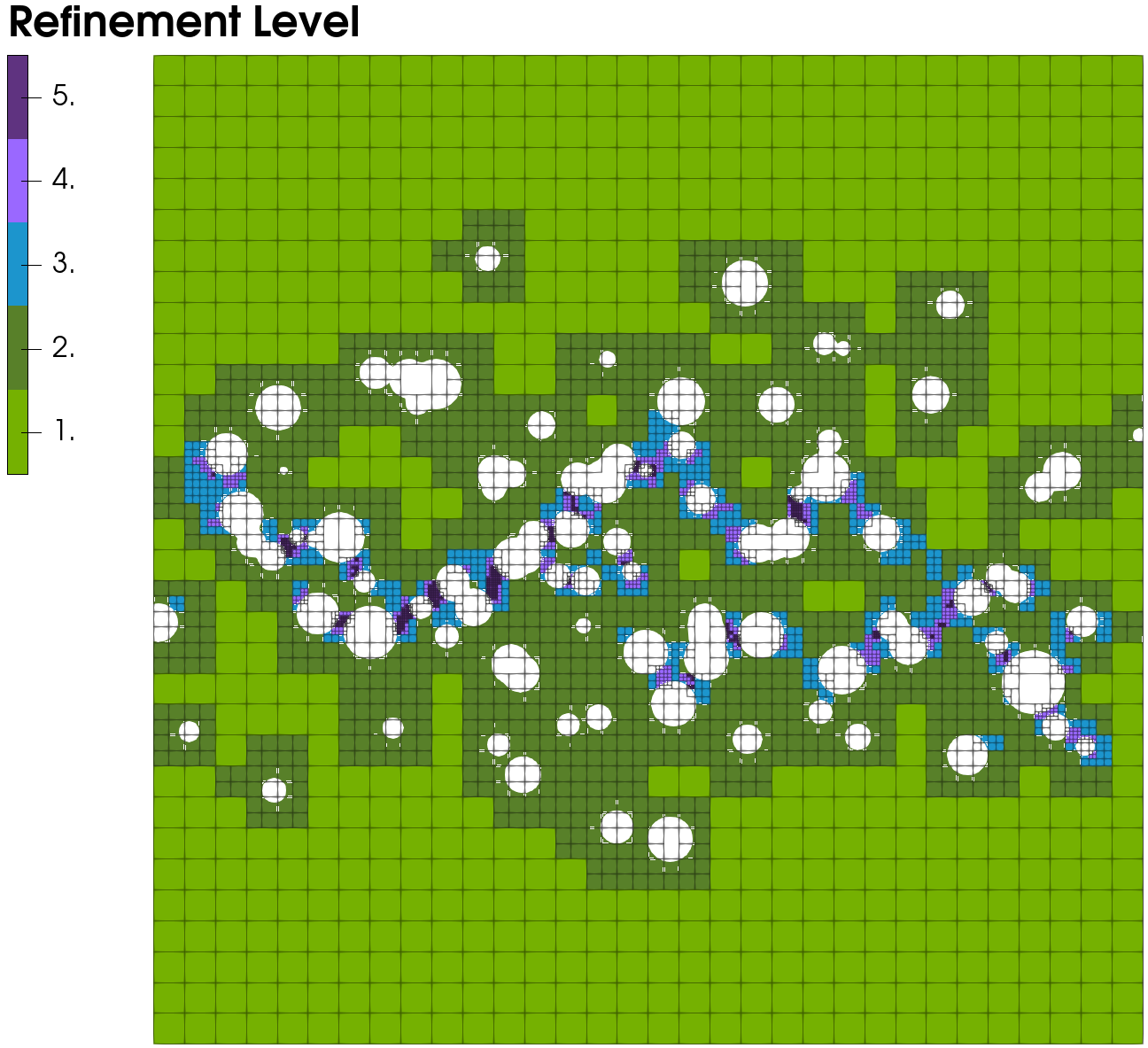} &
		\includegraphics[width = 0.32\textwidth]{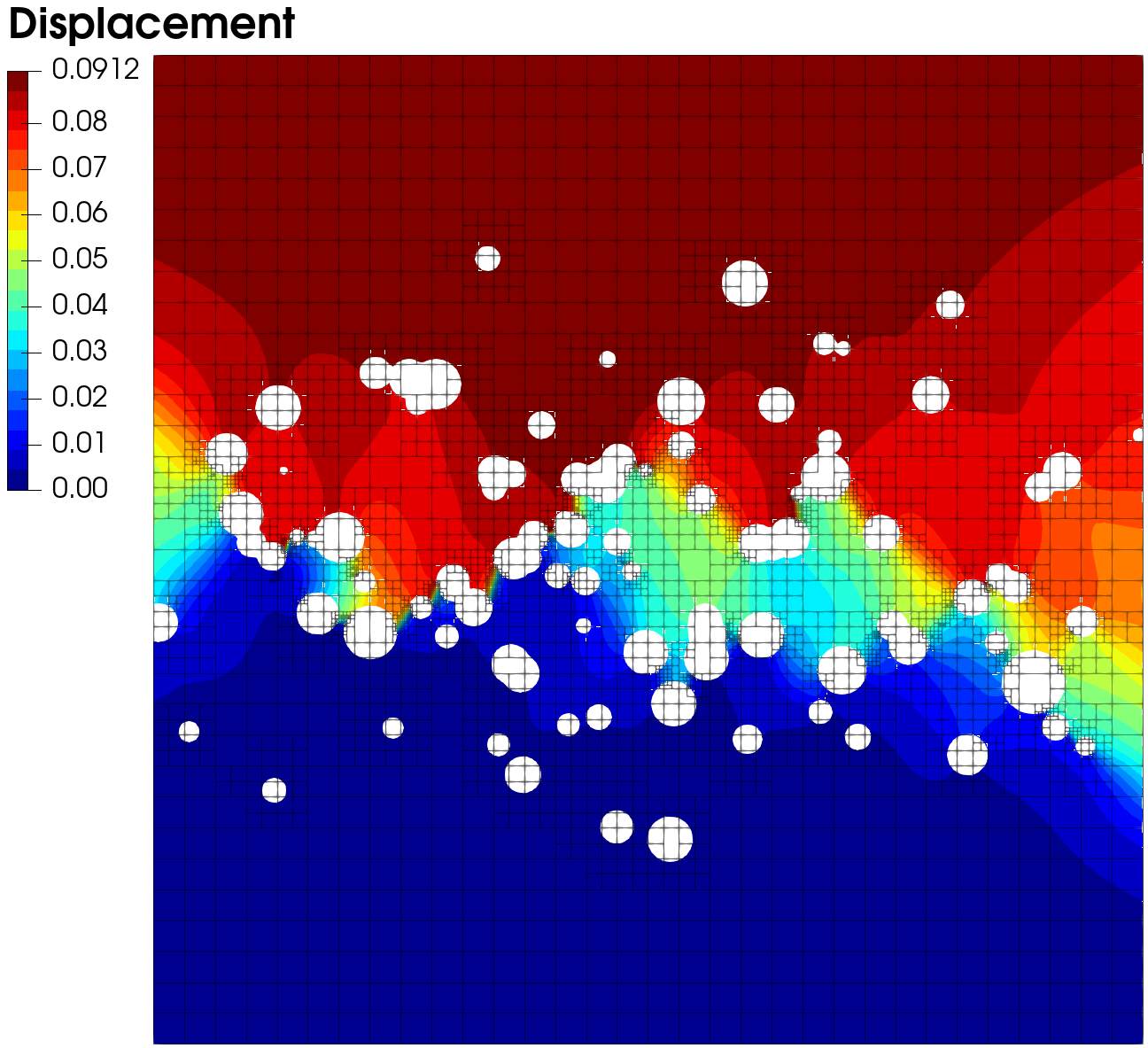} &
		\includegraphics[width = 0.32\textwidth]{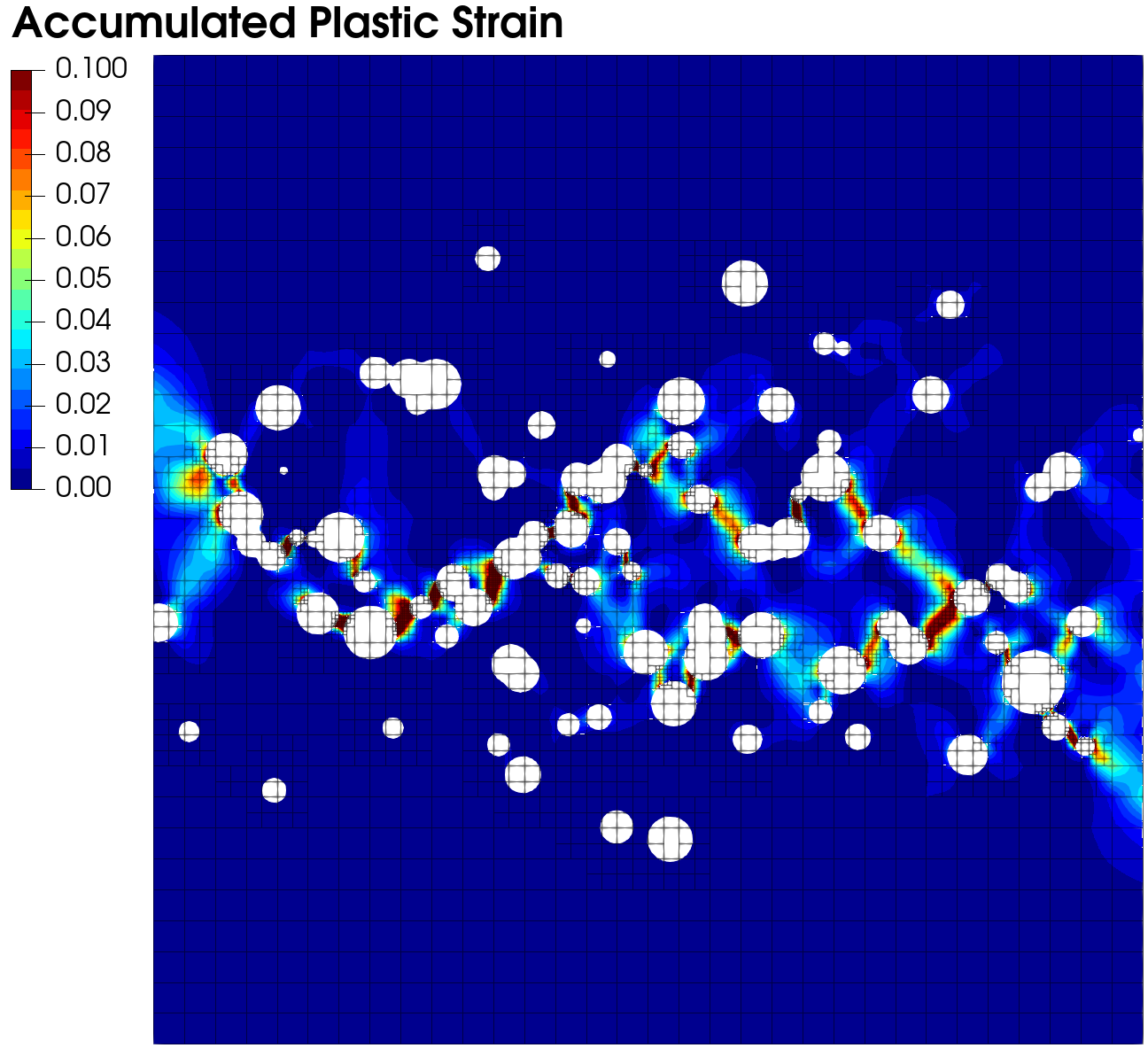} \\
		\includegraphics[width = 0.32\textwidth]{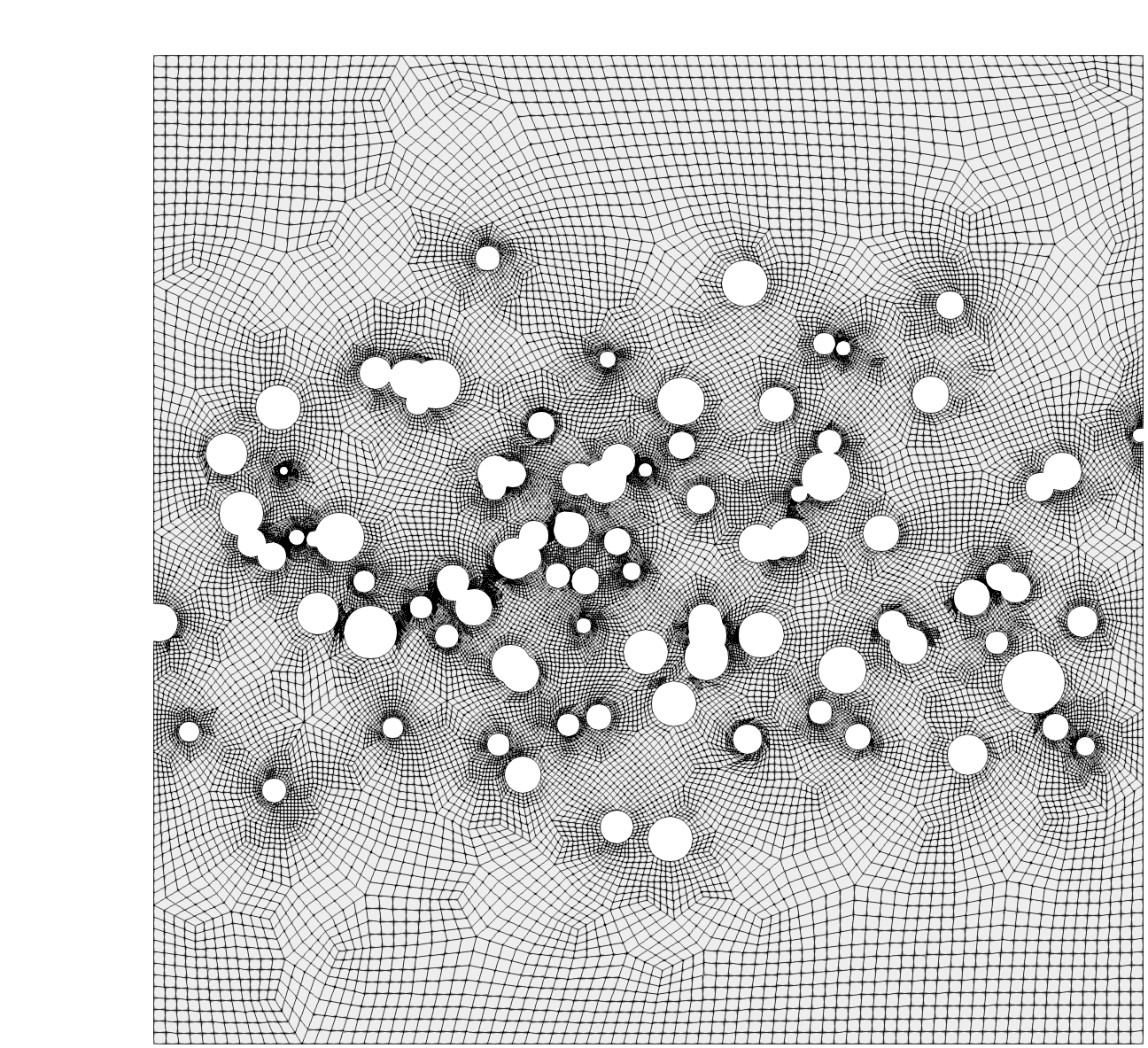} &
		\includegraphics[width = 0.32\textwidth]{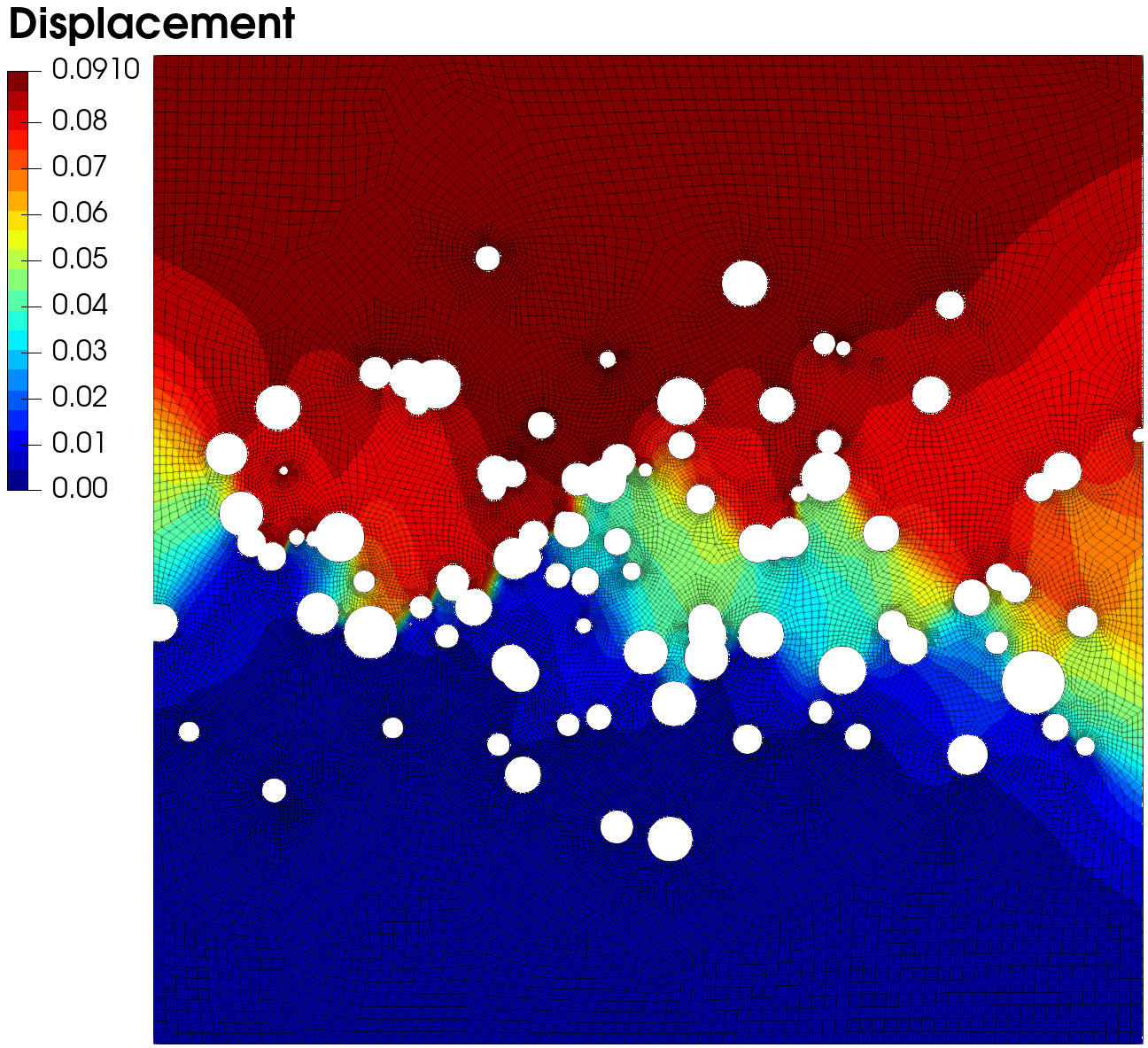} &
		\includegraphics[width = 0.32\textwidth]{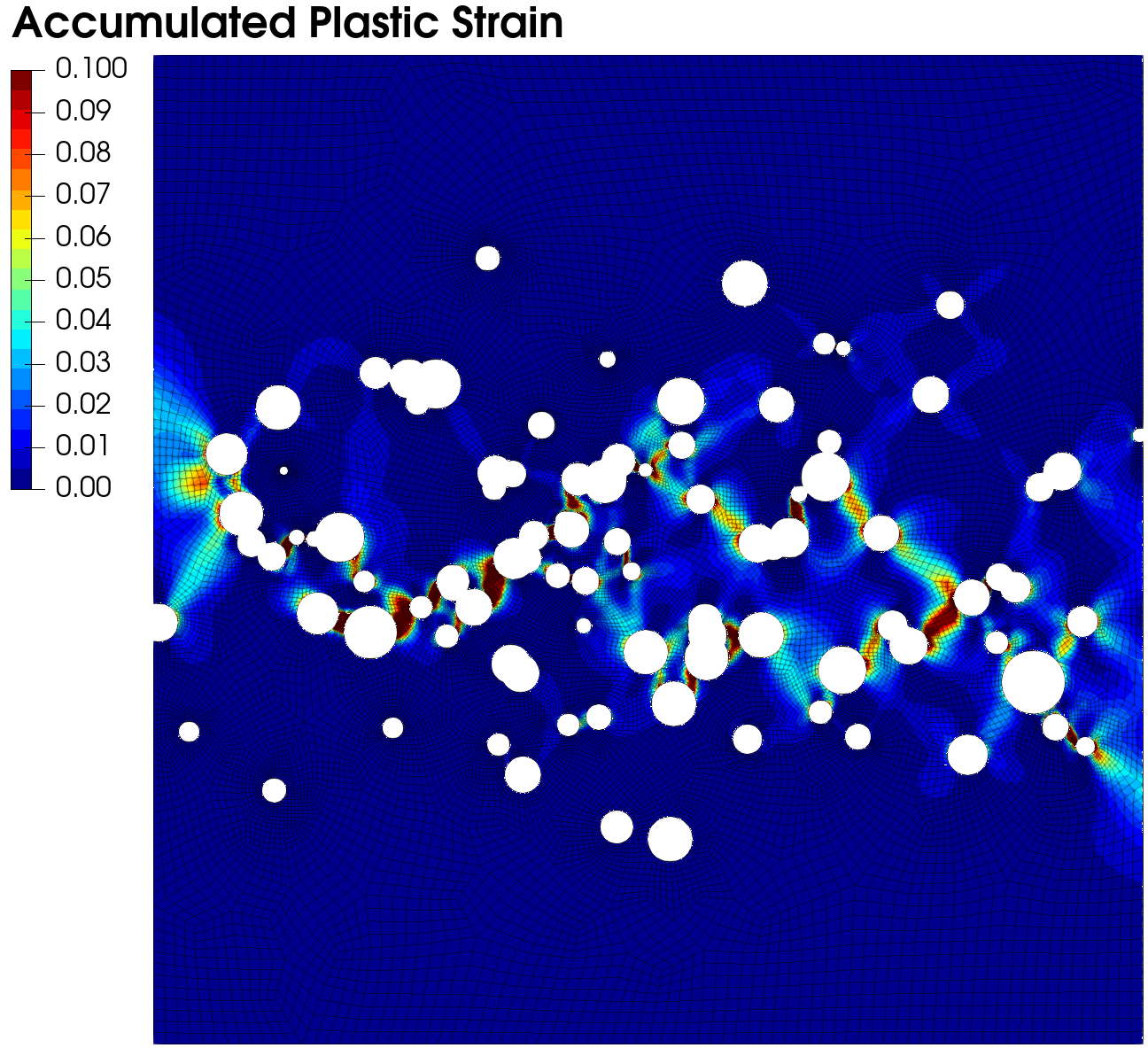} \\
	\end{tabular}
	\caption{Porous plate: refinement levels (left), displacement field (mid), and accumulated plastic strain field \(\bar{\varepsilon}^{p}(\mathbf{x})\) (right). The upper row shows the initial discretization with geometry-induced refinement \color{mygreenD}(i.e., refinement required to resolve the immersed geometry) \color{black}but without adaptive refinement. The middle row shows the solution after \(12\) cycles, with refinement levels ranging from \(1\) to \(5\), where level \(1\) denotes the base mesh. The bottom row shows the corresponding body-fitted reference FEM solution.}
	\label{PorousPlateRefinement}
\end{figure}

\begin{figure}[htb!]
	\centering
	\begin{tabular}{@{} c @{\space} c @{}}
		\includegraphics[width=0.475\textwidth]{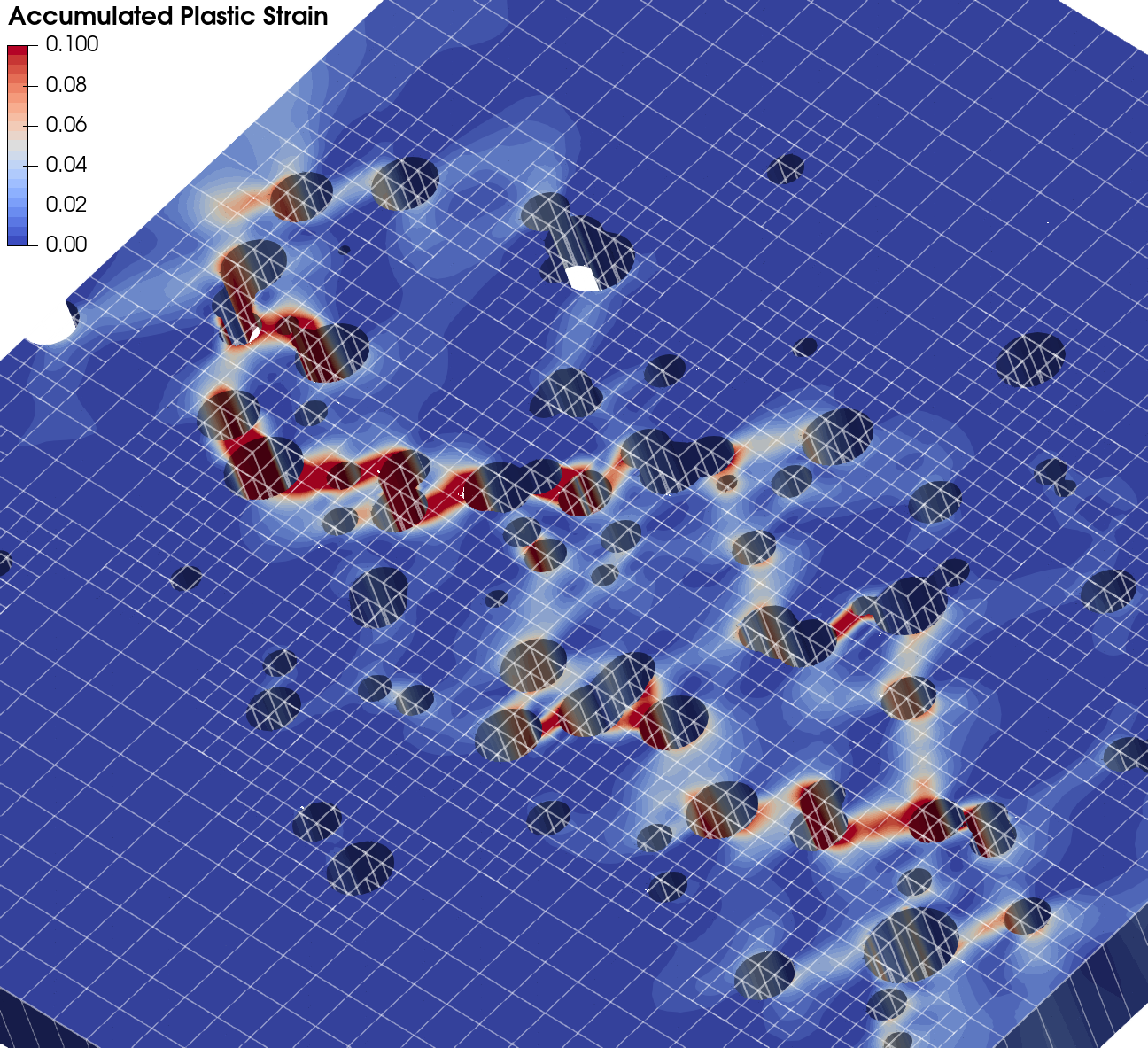} &
		\includegraphics[width=0.475\textwidth]{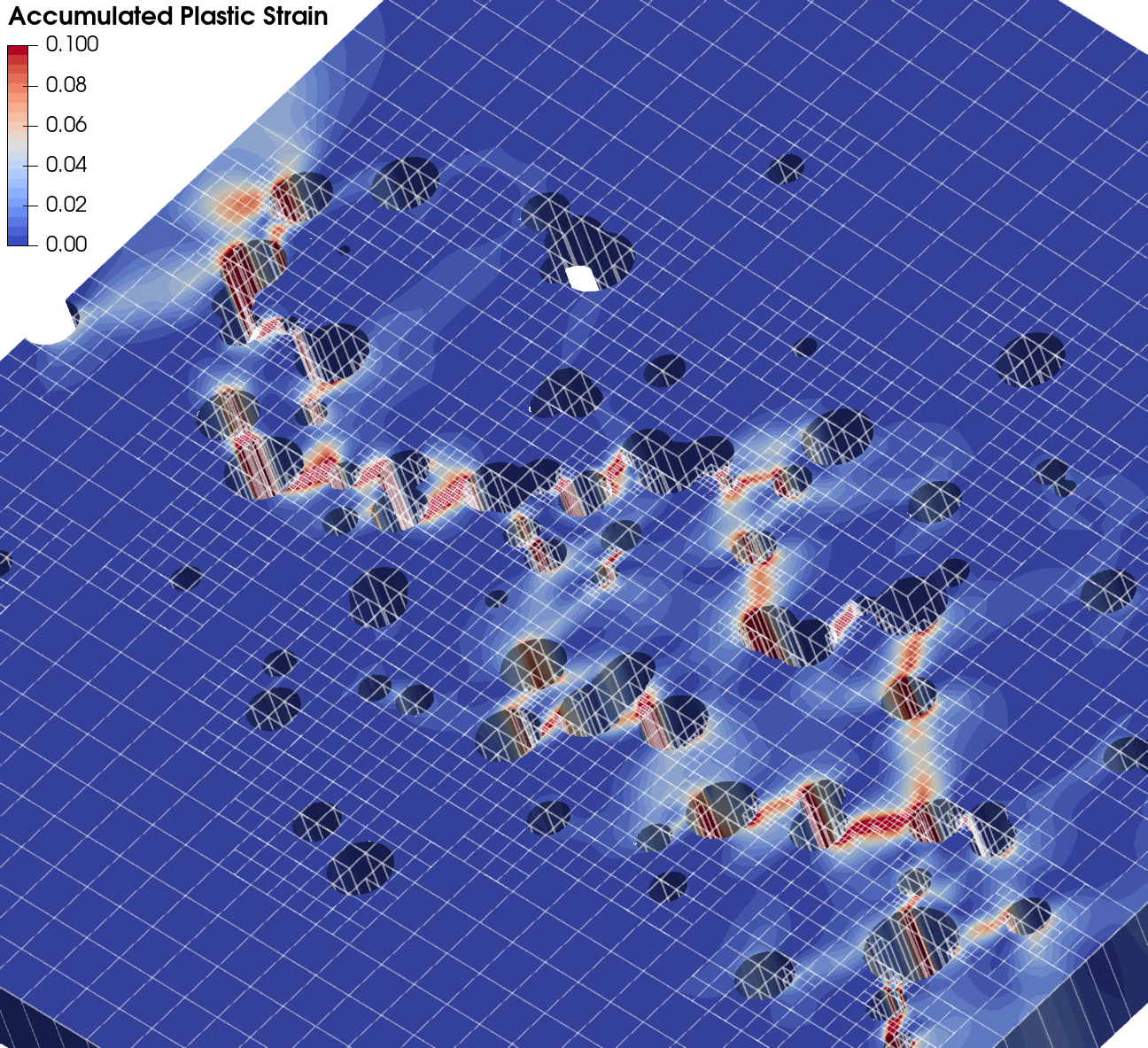} \\
	\end{tabular}
	\caption{Porous plate: comparison of the accumulated plastic strain \(\bar{\varepsilon}^{p}\) in cycle \(0\) (left) and cycle \(12\) (right). The final discretization has a maximum refinement depth of four overlay levels (base mesh + four overlays). The deformation is shown with a scale factor of \(30\).}
	\label{PorousPlate_AccumStrain_B}
\end{figure}

\paragraph{Overlay refinement}
The evolution of the adaptive discretization together with the corresponding displacement and accumulated plastic strain fields is illustrated in Figure~\ref{PorousPlateRefinement}. The error indicator selectively concentrates overlay refinement along the developing deformation band and around pores interacting with this region, while large parts of the domain remain at the initial refinement level, cf. Fig.~\ref{PorousPlate_AccumStrain_B}. Successive refinement increasingly resolves the narrow plastic-strain bands connecting neighboring pores and yields close qualitative agreement with the body-fitted reference FEM solution shown in the bottom row. 

In the present work, a single shared mesh is used for both the thermal and mechanical subproblems. Independent field-wise refinement is also possible and may further improve efficiency, since different solution fields can exhibit substantially different spatial regularity~\cite{Solin:2010}. In particular, the comparatively smooth temperature field may require less spatial resolution than the localized mechanical response observed here. Such field-specific adaptivity is left for future investigation.

The convergence behavior of the reaction force is assessed in Fig.~\ref{PorousPlate_Conv}. While the unrefined FCM solution noticeably overestimates the reaction force, successive Kelly-based refinement approaches both the uniformly refined FCM and the body-fitted reference solution over the entire loading path. The final adaptive solution requires only \(67\,803\) unknowns, compared with approximately \(493\,000\) (factor \(7.3\)) and \(814\,239\) (factor \(12\)) unknowns for the uniform FCM and body-fitted reference solutions, respectively. More importantly, the diagram on the right reports the error in the reaction force at load step~$30$. \color{mygreenD} The extrapolated value \(F_{30}^{\infty}=0.8401535087\,\mathrm{kN}\) used for the convergence analysis was obtained by second-order Richardson extrapolation from two uniformly refined FCM solutions with \(N_1=115\,434\) and \(N_2=447\,957\) degrees of freedom and corresponding reaction forces \(F_{30}^{(1)}=0.9640641128\,\mathrm{kN}\) and \(F_{30}^{(2)}=0.8903307468\,\mathrm{kN}\), respectively. Using the effective refinement ratio \(r=(N_2/N_1)^{1/3}=1.5714510621\), the extrapolated value follows from
\begin{equationarray*}{rcl}
	F_{30}^{\infty} &=& F_{30}^{(2)}+\frac{F_{30}^{(2)}-F_{30}^{(1)}}{r^2-1}\; .

\end{equationarray*} \color{black} The adaptive FCM exhibits a markedly steeper error reduction, with an observed convergence rate of approximately \(2.24\), compared with \(0.63\) for uniform FCM refinement and \(0.28\) for the body-fitted \(h\)-refinement sequence. These results demonstrate that the error-indicator-driven refinement substantially improves the accuracy per degree of freedom by selectively resolving the regions governing the non-linear response.

\begin{figure}[htb!]
	\centering 
	\begin{tikzpicture}

\begin{axis}[
  name=main,
  width=0.475\textwidth, height=9cm,
  xlabel={prescribed top displacement $\overline{u}_{y}$ [\si{\milli\metre}]},
  ylabel={reaction force $F_y$ [\si{\kilo\newton}]},
  xmin=0, xmax=0.0945, ymin=0, ymax=1.10,
  xtick={0,0.015,0.03,0.045,0.06,0.075,0.09},
  ytick={0,0.2,0.4,0.6,0.8,1.0},
  scaled x ticks=false,
  xticklabel style={/pgf/number format/fixed, /pgf/number format/precision=3},
  grid=both, minor tick num=1,
  legend cell align=left,
  legend style={
    at={(0.985,0.030)}, anchor=south east,
    draw=oiGrey, line width=0.4pt, fill=white, fill opacity=0.95,
    text opacity=1, font={\fontsize{8pt}{8pt}\selectfont}, row sep=1.6pt, inner sep=5pt,
  },
]

\addplot[oiVermillion, line width=2.4pt, mark=none] coordinates {
      (0.003000,0.162161574126) (0.006000,0.415649859896) (0.009000,0.566778259934) (0.012000,0.633232682325) (0.015000,0.674084499973) (0.018000,0.703212200428)
      (0.021000,0.725807846186) (0.024000,0.744252471294) (0.027000,0.759715674206) (0.030000,0.772939587757) (0.033000,0.784369640362) (0.036000,0.794400017608)
      (0.039000,0.803351757700) (0.042000,0.811455788815) (0.045000,0.818879330262) (0.048000,0.825736962784) (0.051000,0.832110975164) (0.054000,0.838066394482)
      (0.057000,0.843654597587) (0.060000,0.848916562908) (0.063000,0.853886447674) (0.066000,0.858591048729) (0.069000,0.863054399605) (0.072000,0.867296525504)
      (0.075000,0.871336208492) (0.078000,0.875188908801) (0.081000,0.878868913946) (0.084000,0.882388374193) (0.087000,0.885757857271) (0.090000,0.888986455084)
};
\addlegendentry{FCM, uniform $168^2$, $p{=}2$, \num{492974} dofs}

\addplot[oiBlue, line width=1.05pt, mark=*, mark size=1.35pt,
         mark options={fill=oiBlue, draw=white, line width=0.3pt}, mark repeat=2]
  coordinates {
      (0.003000,0.163006308467) (0.006000,0.416699813427) (0.009000,0.568627767191) (0.012000,0.636087837548) (0.015000,0.676924480764) (0.018000,0.705585957643)
      (0.021000,0.727706499311) (0.024000,0.745723905571) (0.027000,0.760780764124) (0.030000,0.773657120057) (0.033000,0.784827924806) (0.036000,0.794633910873)
      (0.039000,0.803395154535) (0.042000,0.811342580310) (0.045000,0.818633243489) (0.048000,0.825387332972) (0.051000,0.831681402588) (0.054000,0.837571572998)
      (0.057000,0.843107074808) (0.060000,0.848321772034) (0.063000,0.853250964789) (0.066000,0.857923931938) (0.069000,0.862357618780) (0.072000,0.866574825013)
      (0.075000,0.870593577469) (0.078000,0.874429058497) (0.081000,0.878094018810) (0.084000,0.881598670740) (0.087000,0.884954398925) (0.090000,0.888169427857)
};
\addlegendentry{FCM, Kelly $\theta{=}0.4$, cycle 12, \num{67803} dofs}

\addplot[oiGreen, line width=1.05pt, dash pattern=on 5pt off 2.4pt,
         mark=triangle*, mark size=1.9pt,
         mark options={solid, fill=oiGreen, draw=white, line width=0.3pt}, mark repeat=2]
  coordinates {
      (0.003000,0.165002911764) (0.006000,0.434536942645) (0.009000,0.608399172548) (0.012000,0.686199181019) (0.015000,0.732552813985) (0.018000,0.765908041817)
      (0.021000,0.792365385083) (0.024000,0.814156119467) (0.027000,0.832650107646) (0.030000,0.848727078653) (0.033000,0.862996912504) (0.036000,0.875827074240)
      (0.039000,0.887495250158) (0.042000,0.898217433098) (0.045000,0.908134006973) (0.048000,0.917369711732) (0.051000,0.926018335606) (0.054000,0.934157191176)
      (0.057000,0.941839498288) (0.060000,0.949112200936) (0.063000,0.956015232605) (0.066000,0.962578778323) (0.069000,0.968834200618) (0.072000,0.974804459905)
      (0.075000,0.980509356557) (0.078000,0.985967616029) (0.081000,0.991195752024) (0.084000,0.996209788152) (0.087000,1.001027306470) (0.090000,1.005658309050)
};
\addlegendentry{FCM, Kelly $\theta{=}0.4$, cycle 0, \num{39328} dofs}

\addplot[oiPurple, line width=1.05pt, dash pattern=on 3.2pt off 1.6pt on 0.8pt off 1.6pt,
         mark=square*, mark size=1.5pt,
         mark options={solid, fill=oiPurple, draw=white, line width=0.3pt}, mark repeat=2]
  coordinates {
      (0.003000,0.162643) (0.006000,0.414382) (0.009000,0.563715) (0.012000,0.629235) (0.015000,0.669305) (0.018000,0.697751)
      (0.021000,0.719751) (0.024000,0.737691) (0.027000,0.752726) (0.030000,0.765568) (0.033000,0.776654) (0.036000,0.786390)
      (0.039000,0.795084) (0.042000,0.802968) (0.045000,0.810200) (0.048000,0.816887) (0.051000,0.823108) (0.054000,0.828925)
      (0.057000,0.834386) (0.060000,0.839532) (0.063000,0.844393) (0.066000,0.848997) (0.069000,0.853366) (0.072000,0.857520)
      (0.075000,0.861477) (0.078000,0.865252) (0.081000,0.868858) (0.084000,0.872308) (0.087000,0.875611) (0.090000,0.878776)
};
\addlegendentry{body-fitted, $h$-scale ladder, \num{814239} dofs}

\end{axis}

\end{tikzpicture}
	\begin{tikzpicture}
\begin{axis}[
  width=0.475\textwidth, height=9cm,
  xmode=log, ymode=log,
  xlabel={unknowns $N$},
  ylabel={$|F_{30}-F_{30}^{\infty}|$ [\si{\kilo\newton}]},
  xmin=15500, xmax=1000000, ymin=0.030, ymax=0.80,
  xtick={20000,50000,100000,200000,500000,1000000},
  xticklabels={$2\!\cdot\!10^{4}$,$5\!\cdot\!10^{4}$,$10^{5}$,$2\!\cdot\!10^{5}$,$5\!\cdot\!10^{5}$,$10^{6}$},
  ytick={0.015,0.02,0.03,0.05,0.07,0.1,0.2,0.3,0.5},
  yticklabels={0.015,0.02,0.03,0.05,0.07,0.10,0.20,0.30,0.50},
  minor xtick={30000,40000,60000,70000,80000,90000,300000,400000,600000,700000,800000,900000},
  minor ytick={0.04,0.06,0.08,0.09,0.4},
  grid=both,
  legend cell align=left,
  legend style={
    at={(0.978,0.965)}, anchor=north east,
    draw=oiGrey, line width=0.4pt, fill=white, fill opacity=0.95,
    text opacity=1, font={\fontsize{8pt}{8pt}\selectfont}, row sep=1.6pt, inner sep=5pt,
  },
]

\addplot[oiVermillion, line width=1.3pt, mark=square*, mark size=2.3pt,
         mark options={fill=oiVermillion, draw=white, line width=0.4pt}]
  coordinates {
      (18943,0.366189452054) (29434,0.30514333058) (74463,0.164694887324) (115434,0.124064112775)
      (139121,0.107562091286) (289728,0.073282578897) (447957,0.050330746769) (492974,0.048986455084)
};
\addlegendentry{FCM, uniform ladder, $p{=}2$}

\addplot[oiBlue, line width=1.3pt, mark=*, mark size=2.1pt,
         mark options={fill=oiBlue, draw=white, line width=0.4pt}]
  coordinates {
      (39328,0.16565830905) (39933,0.156819613657) (41040,0.149107784965) (42369,0.135048868604)
      (43985,0.110619957366) (44907,0.097022525308) (46764,0.092668402225) (47972,0.087774904535)
      (51014,0.08113942735) (53488,0.076069222222) (57745,0.065857262688) (63205,0.052539056075)
      (67803,0.048169427857)
};
\addlegendentry{FCM, Kelly $\theta{=}0.4$}

\addplot[oiPurple, line width=1.1pt, mark=triangle*, mark size=2.5pt,
         mark options={fill=oiPurple, draw=white, line width=0.4pt}]
  coordinates {
      (169983,0.060811) (274311,0.054179889539) (443151,0.045662) (814239,0.038776)
  };
\addlegendentry{body-fitted, $h$-scale ladder}

\tikzset{convtri/.style={line width=0.6pt, solid, line join=miter}}
\tikzset{convtrilbl/.style={font=\fontsize{8pt}{8pt}\selectfont, inner sep=2pt}}

\draw[convtri, oiVermillion]
  (axis cs:240000,0.260000) -- (axis cs:240000,0.130204)
  -- (axis cs:720000,0.130204) -- cycle;
\node[convtrilbl, oiVermillion, anchor=east]  at (axis cs:240000,0.183992) {$0.63$};
\node[convtrilbl, oiVermillion, anchor=north] at (axis cs:415692,0.130204) {$1$};

\draw[convtri, oiBlue]
  (axis cs:80000,0.080000) -- (axis cs:112000,0.080000)
  -- (axis cs:112000,0.037608) -- cycle;
\node[convtrilbl, oiBlue, anchor=west]  at (axis cs:112000,0.054851) {$2.24$};
\node[convtrilbl, oiBlue, anchor=south] at (axis cs:94657,0.080000) {$1$};

\draw[convtri, oiPurple]
  (axis cs:360000,0.075000) -- (axis cs:600000,0.075000)
  -- (axis cs:600000,0.064845) -- cycle;
\node[convtrilbl, oiPurple, anchor=west]  at (axis cs:600000,0.069738) {$0.28$};
\node[convtrilbl, oiPurple, anchor=south] at (axis cs:464758,0.075000) {$1$};

\end{axis}
\end{tikzpicture}
	\caption{Porous plate: reaction force--displacement response (left) and convergence of the reaction force at the final load step with respect to the number of unknowns (right).}
	\label{PorousPlate_Conv}
\end{figure}
%

\section{Summary and conclusions}
\label{sec:conclusions}
In this work, an extended multi-level \(hp\)-Finite Cell Method has been developed for the simulation of thermo-viscoplastic problems with temperature-dependent material behavior. The proposed framework combines advanced quadrature and adaptive refinement techniques to address the challenges associated with non-linear, history-dependent constitutive models on complex geometries.

The key aspects can be summarized as follows:
\begin{itemize}
	\item[$\filledsquare$] The \emph{Non-Negative Moment Fitting (NNMF)} quadrature rule outperforms the
established \emph{Adaptive Space-Tree (AST)} quadrature, which is highly reliable and well suited for linear problems, in several respects. The drastic reduction in the number of quadrature points required to accurately capture the physical
solution domain has a significant impact on the numerical complexity and overall solution effort, particularly in non-linear incremental analyses where the integration is repeatedly evaluated. For the plate with a circular hole, the
quadrature efficiency is demonstrated by a reduction in the number of quadrature points of approximately \(75\%\) compared with AST. For the three-dimensional porous plate, the computational performance was assessed directly in terms of
execution time, demonstrating a speedup of approximately \(20\times\) for NNMF over AST. 

	\item[$\filledsquare$] The evolving plasticity behavior requires numerically intensive iterative solution procedures to satisfy static equilibrium at the global level and the constitutive equations at the local level. The reduced integration effort of NNMF therefore becomes particularly beneficial in incremental non-linear analyses, where constitutive evaluations are repeatedly performed throughout the loading history. 

	\item[$\filledsquare$] The NNMF approach offers substantial advantages in storing history-dependent
variables due to the significantly smaller number of required quadrature points. For the porous plate example, NNMF provides comparable solution accuracy to the body-fitted FEM reference solution while requiring approximately \(64\%\) fewer quadrature points than the latter, which employs a \((p+1)^3\) Gauss quadrature rule.

	\item[$\filledsquare$] The hierarchical multi-level \(hp\)-refinement with the overlay concept has been enriched with a gradient-based error indicator that reliably drives adaptive \(h\)-refinement in the non-linear regime, thus facilitating automatic model refinement to capture localized phenomena in the evolving non-linear response.

	\item[$\filledsquare$]  A thermo-viscoplastic material model for the case-hardening steel 16MnCr5  under elevated-temperature conditions was implemented within the FCM framework and applied to examples of complex geometry. In combination with the implemented quadrature scheme and hierarchical adaptive refinement, the model has demonstrated excellent numerical stability and reliability at a significantly reduced numerical cost, with accurate results already obtained at a low to medium polynomial degree.
\end{itemize}

In sum, we have implemented, extended and critically tested an analysis platform for thermo-viscoplastic problems that offers considerable flexibility and depth in modeling the desired response. Moreover, the results indicate that the proposed approach alleviates the computational burden commonly associated with highly resolved cut-cell integration in non-linear three-dimensional analyses.\\

\noindent \textbf{Acknowledgments:} Jan Niklas Schm\"{a}ke and Martin Ruess gratefully acknowledge the support provided by the Deutsche Forschungsgemeinschaft (DFG, German Research Foundation) under the funding code RU 885/4-1, project number 515687474. 

The contribution of Oliver Wege was funded by the Deutsche Forschungsgemeinschaft (DFG, German Research Foundation) under Germany's Excellence Strategy - EXC-2023 Internet of Production - 390621612.\\[6mm]

\noindent\textbf{Appendix A: Analytical Benchmark Solution - One-dimensional Example}\\[-2mm]

\noindent The analytical solution of the one-dimensional benchmark problem of Section \ref{sec:oneDWMaterialJump} reads:

\begin{equation}
	u(x) \;\;=\;\;
		\begin{cases}
			\hat{u}(x) & x<a\\ 
			t_0\,(x-a)/E_0+\hat{u}(a) & x\geq a\\ 
		\end{cases} \hfill \
\end{equation}

\noindent with
\begin{equationarray}{lrclclclclc}
	 \hat{u}(x) &=& \frac{1}{\Delta E}a\,t_0\,\left[\log(a(\Delta E+E_0))-\log(a(\Delta E+E_0)-\Delta E\,x)\right]\nonumber
		 \\[5mm]
		&&\forall x\in \Omega \,=\,[0,a)\cup[a,b] \mbox{ and } 0<a<b \, .\nonumber
\end{equationarray}

\noindent We used the energy norm \(a(u,u)\) to assess the quality of the approximate finite element solution compared to the exact solution:
\begin{equationarray}{lrclclclclc}
	 a(u,u) &=& t_0^2\,\left[\frac{b-a}{E_0}+\frac{a}{\Delta E}\,\log\biggl(\frac{\Delta E+E_0}{E_0}\biggr) \right] \;.
\end{equationarray}
\vspace{3mm}

\noindent\textbf{Appendix B: Analytical Benchmark Solution - Two-dimensional Example}\\[-2mm]

\noindent The analytical solution of the two-dimensional benchmark problem of Section \ref{sec:twoDContinuousMaterial} is expressed through the piecewise defined equation:
\begin{equationarray}{rclclclclc}
	u_r(r) &=& 
	\begin{cases}
	\hat{u}_r(r) & 0\leq r < a\\
	\frac{1}{2a\,E_0\,r}\bigl[-t_0\,a^3+a^2\,E_0\,\hat{u}_r(a)+a\,r^2t_0+E_0\,r^2\,\hat{u}_r(a) \bigr] \quad \ & a\leq r \leq b
	\end{cases} 
\end{equationarray}
	with 
\begin{equationarray}{rclclclclc}
	\hat{u}_r(r) &=& \frac{1}{\alpha+\beta}\left[3r\,t_0(\Delta E+E_0) {}_2F_1\biggl(\frac{1}{2}\bigl(3-\sqrt{5}\bigr), \frac{1}{2}\bigl(3+\sqrt{5}\bigr); 3; \frac{\Delta E\,r}{a\,\Delta E+a\,E_0}\biggr)\right]
\end{equationarray}
and the constants \(\alpha\) and \(\beta\) defined as
\begin{equationarray}{rclclclclc}
	\alpha &=& 3\,E_0(\Delta E+E_0){}_2F_1\biggl(\frac{1}{2}\bigl(3-\sqrt{5}\bigr),\frac{1}{2}\bigl(3+\sqrt{5}\bigr); 3; \frac{\Delta E}{\Delta E+E_0}\biggr)\\[4mm]
		\beta &=& \Delta E\,E_0\,{}_2F_1\biggl(\frac{1}{2}(5 - \sqrt{5}), \frac{1}{2}(5 + \sqrt{5}); 4 ; \frac{\Delta E}{\Delta E+E_0}\biggr)
	\end{equationarray}
in which \(({}_2F_1(a;b;c;z))\) refers to the hypergeometric function \cite{Aomoto}.\\

\vspace{4mm}

\noindent\textbf{Appendix C: Temperature-dependent Material Parameter}\\[-2mm]

\noindent Temperature-dependent material properties of case hardening steel 16MnCr5, together with the corresponding fitting parameters, obtained from interpolation functions fitted to experimental data \((\times)\) reported in \cite{Fukuhara:1993,Spittel:2009,Liscic,Spittel:2009.2} with an austenitisation temperature \(T_A^\mathrm{16MnCr5}=750^{\circ}\mathrm{C}\) and a solidus temperature \(T_S^\mathrm{Fe}=1538^{\circ}\textrm{C}\). The discontinuity at \(750^{\circ}\mathrm{C}\) reflects the austenitisation property of a change of the crystalline structure at this temperature. Suitable parameters for \(\kappa\) and \(\gamma_m\) are therefore obtained by a piecewise data interpolation, cf. example \ref{sec:porousPlate}.
\begin{figure}[h!]
	\centering
	\tikzsetnextfilename{model_tvp_mechanical_parameters}
\begin{tikzpicture}[spy using outlines = {densely dashed, rectangle, magnification=14, connect spies}]
	\begin{axis} [
		font={\footnotesize},
		axis lines = box,
		title = Elastic moduli \vphantom{$\sigma_\text{y}$} in GPa,
		width = .45\textwidth,
		height = 5cm,
		legend style = {at = {(.96,.98)},nodes = {scale = 0.72, transform shape}},
		legend cell align={left},
		enlarge x limits=0.00,
		enlarge y limits=0.05,
    extra tick style={grid=major, grid style={black}, major tick length=2mm, tick align=center,ticklabel pos=top, tick style={black}},
		title style = {yshift = -2mm},
		xticklabel={$ \pgfmathprintnumber[fixed,precision=0,1000 sep={}]{\tick}^\circ\textrm{C}$},
		grid = major,
		grid style = {densely dashed, line width = 0.1pt}
	]

	  \addplot+[mark = none, blue, thick, each nth point=5, filter discard warning=false, unbounded coords=discard] table [col sep = comma]
			{model_tvp_emodR2.csv};

		\pgfplotsset{cycle list shift = -1}
		\addplot+[mark = x, blue!70!black, draw=none, mark size = 2.2pt] table [col sep = comma]
			{model_tvp_emod_expR2.csv};

		\addplot+[mark = none, orange, thick, each nth point=5, filter discard warning=false, unbounded coords=discard] table [col sep = comma]
			{model_tvp_kmodR2.csv};

		\pgfplotsset{cycle list shift = -2}
		\addplot+[mark = x, orange!70!black, draw=none, mark size = 2.2pt] table [col sep = comma]
			{model_tvp_kmod_expR2.csv};

		\addplot+[mark = none, green, thick, each nth point=5, filter discard warning=false, unbounded coords=discard] table [col sep = comma]
			{model_tvp_gmodR2.csv};

		\pgfplotsset{cycle list shift = -3}
		\addplot+[mark = x, green!70!black, draw=none, mark size = 2.2pt] table [col sep = comma]
			{model_tvp_gmod_expR2.csv};

		\legend{$E$, ,$K$, ,$G$, };
	\end{axis}

	\begin{axis} [
		shift = {(0.47\textwidth, 0cm)},
		font={\footnotesize},
		axis lines = box,
		title = Yield stress $\sigma_\text{y}$ in  GPa,
		width = .45\textwidth,
		height = 5cm,
		legend style = {at = {(.96,.98)},nodes = {scale = 0.72, transform shape}},
		legend cell align={left},
		enlarge x limits=0.00,
		enlarge y limits=0.05,
    extra tick style={grid=major, grid style={black}, major tick length=2mm, tick align=center,ticklabel pos=top, tick style={black}},
		title style = {yshift = -2mm},
		xticklabel={$ \pgfmathprintnumber[fixed,precision=0,1000 sep={}]{\tick}^\circ\textrm{C}$},
		grid = major,
		grid style = {densely dashed, line width = 0.1pt}
	]

	  \addplot+[mark = none, blue, thick, each nth point=5, filter discard warning=false, unbounded coords=discard] table [col sep = comma]
			{model_tvp_y0R2.csv};
		\pgfplotsset{cycle list shift = -1}
		\addplot+[mark = x, blue!70!black, draw=none, mark size = 2.2pt] table [col sep = comma]
			{model_tvp_y0_expR2.csv};
		\addplot+[mark = none, orange, thick, each nth point=5, filter discard warning=false, unbounded coords=discard] table [col sep = comma]
			{model_tvp_yinfR2.csv}; 

		\legend{Initial, , Saturated};

	\end{axis}
\end{tikzpicture} 
	\caption{Temperature-dependent elastic moduli and yield stress $\sigma_y(T).$}
	\label{ConstitutivePropertiesI}
\end{figure}

\vspace{3mm}

\begin{figure}[h!]
	\centering
	\tikzsetnextfilename{materialInterpolationThermalProps}
\begin{tikzpicture}[spy using outlines = {densely dashed, rectangle, magnification=14, connect spies}]

  \node at (-.5,0) {\,};

	\begin{axis} [
		font={\footnotesize},
		axis lines = box,
		title = Heat conductivity $\kappa$ in $\sfrac{\text{W}}{\text{m}\cdot\text{K}}$ \vphantom{$\sigma_\text{y}$} ,
		width = .45\textwidth,
		height = 5cm,
		ytick distance = 5,
		legend style = {at = {(.96,.98)},nodes = {scale = 0.72, transform shape}},
		legend cell align={left},
		enlarge x limits=0.00,
		enlarge y limits=0.05,
    extra tick style={grid=major, grid style={black}, major tick length=2mm, tick align=center,ticklabel pos=top, tick style={black}},
		title style = {yshift = -2mm},
		xticklabel={$ \pgfmathprintnumber[fixed,precision=0,1000 sep={}]{\tick}^\circ\textrm{C}$},
		grid = major,
		grid style = {densely dashed, line width = 0.1pt}
	]

	  \addplot+[mark = none, thick] table [col sep = comma]
      {model_tvp_thconR2.csv};

		\addplot+[mark = x, blue!70!black, only marks, mark size = 2.2pt] table [col sep = comma]
      {model_tvp_thcon_expR2.csv};
	\end{axis}

	\begin{axis} [
		shift = {(0.47\textwidth, 0cm)},
		font={\footnotesize}, 
		axis lines = box,
		title = Expansion coefficient $\gamma_m$ in  $\sfrac{\text{10}^\text{-6}}{\text{K}}$\vphantom{$\sigma_\text{y}$} ,
		width = .45\textwidth,
		height = 5cm,
		legend style = {at = {(.93,.98)},nodes = {scale = 0.72, transform shape}},
		legend cell align={left},
		enlarge x limits=0.00,
		enlarge y limits=0.05,
    extra tick style={grid=major, grid style={black}, major tick length=2mm, tick align=center,ticklabel pos=top, tick style={black}},
		title style = {yshift = -2mm},
		xticklabel={$ \pgfmathprintnumber[fixed,precision=0,1000 sep={}]{\tick}^\circ\textrm{C}$},
		grid = major,
		grid style = {densely dashed, line width = 0.1pt}
	]

	  \addplot+[mark = none, thick] table [col sep = comma]
      {model_tvp_thexpR2.csv};

		\addplot+[mark = x, blue!70!black, only marks, mark size = 2.2pt] table [col sep = comma]
      {model_tvp_thexp_expR2.csv};

	\end{axis}
\end{tikzpicture}
	\caption{Temperature-dependent heat conductivity $\kappa$ and expansion coefficient $\gamma_m(T).$}
	\label{ConstitutivePropertiesII}
\end{figure}
\vspace{3mm}

\noindent\textbf{Appendix D: Thermo-viscoplastic Model - Material Constants }\\[-2mm]

\color{mygreenD}
\noindent The temperature dependencies of the mean thermal expansion coefficient $\gamma_m(T)$ and the heat conductivity $\kappa(T)$ are represented by piecewise interpolation functions with a transition at the austenitisation temperature $T_A$ (cf. Table \ref{vpParameters}). The mean thermal expansion coefficient is given by
\begin{eqnarray*}
    \gamma_m(T) &=&
    \begin{cases}
        \gamma_{m_{01}} \left[1-\omega_{\gamma_{m1}}(T-T_r)\right], & T \leq T_A, \\[1mm]
        \gamma_{m_{02}} \left[1-\omega_{\gamma_{m2}}(T-T_A)\right], & T > T_A .
    \end{cases}
\end{eqnarray*}
The heat conductivity is interpolated according to
\begin{eqnarray*}
    \kappa(T) &=&
    \begin{cases}
        \kappa_{0_{1}}\,\varphi_{\kappa_{1}}^{p_4}(T-T_r),  & T \leq T_A, \\[1mm]
        \kappa_{0_{2}}\,\varphi_{\kappa_{2}}^{\exp}(T-T_A), & T > T_A ,
    \end{cases}
\end{eqnarray*}
where the polynomial contribution is defined as
\begin{eqnarray*}
  \varphi^{p_4}_{\kappa_1}(\Delta T) &=&  1-\omega_{\kappa_{11}}\Delta T\,
  \Bigl[1-\omega_{\kappa_{21}}\Delta T\,
  \bigl(1-\omega_{\kappa_{31}}\Delta T\,
  (1-\omega_{\kappa_{41}}\Delta T)\bigr)\Bigr]
\end{eqnarray*}
and the exponential contribution reads
\begin{eqnarray*}
    \varphi_{\kappa_{2}}^{\exp}(\Delta T) &=& 1-\omega_{\kappa_{2}} \left[\exp(-\chi_{\kappa_{2}}\Delta T)-1\right].
\end{eqnarray*}
The corresponding interpolation parameters are listed in Table~\ref{vpParameters}.\\
\color{black}

The material constants used in the thermo-viscoplastic model for steel grade 16MnCr5, together with the fitting parameters of the temperature-dependent material properties
\([\cdot](T)\in\{K(T),\,G(T),\,\sigma_{y0}(T),\,\Phi_\chi(T)\}\) are listed in Table~\ref{vpParameters}. The corresponding temperature dependencies are interpolated according to \[ [\cdot](T)=[\cdot]_0\left[1-\frac{2\,\omega_{[\cdot]}}{\pi-2\,\varphi_{[\cdot]}}\left(\mathrm{tan}^{-1}\!\bigl(\chi_{[\cdot]}(T-T_r)+\tan(\varphi_{[\cdot]})
\bigr)-\varphi_{[\cdot]}\right)\right],\] cf.~\cite{Oppermann:2022}.
\noindent 
\def\macro#1#2#3{$#1$ & $#2$ & $#3$}

\vspace{4mm}
\begin{table}[h!]
	\centering
	\begin{tabular}{lrllrllrl}
	\toprule
		\multicolumn{3}{w{c}{5cm}}{\colorbox{gray!20}{\makebox[4.8cm]{Temperatures}} }
		&
		\multicolumn{3}{w{c}{5cm}}{\colorbox{gray!20}{\makebox[4.8cm]{Heat flow}} }
		&
		\multicolumn{3}{w{c}{5cm}}{\colorbox{gray!20}{\makebox[4.8cm]{Hardening}} } \\[1mm]

		\macro{T_{\mathrm{ref}}}{273.15}{\mathrm{K}}
		& \macro{\gamma_{m_{01}}}{1.108358}{\sfrac{10^{-5}}{\mathrm{K}}}
		& \macro{\Delta\sigma_{y_{\infty 0}}}{206.3425}{\mathrm{MPa}} \\[1mm]

		\macro{T_{r}}{293.15}{\mathrm{K}}
		& \macro{\omega_{\gamma_{m1}}}{-4.989882}{\sfrac{10^{-4}}{\mathrm{K}}}
		& \macro{\sigma_{y_{00}}}{302.3605}{\mathrm{MPa}} \\[1mm]

		\macro{T_{A}}{1023.15}{\mathrm{K}}
		& \macro{\gamma_{m_{02}}}{1.232077}{\sfrac{10^{-5}}{\mathrm{K}}}
		& \macro{\bar{\varepsilon}_0^p}{2.784959}{10^{-2}} \\[1mm] \cline{7-9}

		\multicolumn{3}{w{c}{5cm}}{\colorbox{gray!20}{\makebox[4.8cm]{Elasticity}} }
		& \macro{\omega_{\gamma_{m2}}}{-6.655428}{\sfrac{10^{-4}}{\mathrm{K}}}
		& \macro{\omega_{\chi}}{-39.22882}{-} \\[1mm] \cline{4-6}

		\macro{K_0}{163.1952}{\mathrm{GPa}}
		& \macro{\kappa_{0_1}}{43.34265}{\sfrac{\mathrm{W}}{\mathrm{m}\cdot \mathrm{K}}}
		& \macro{\chi_{\chi}}{-9.042314}{\sfrac{10^{-3}}{\mathrm{K}}} \\[1mm]

		\macro{\omega_K}{8.592410}{10^{-1}}
		& \macro{\omega_{\kappa_{11}}}{-5.154848}{\sfrac{10^{-4}}{\mathrm{K}}}
		& \macro{\varphi_{\chi}}{-1.461198}{-} \\

		\macro{\chi_K}{4.509029}{\sfrac{10^{-3}}{\mathrm{K}}}
		& \macro{\omega_{\kappa_{21}}}{5.937507}{\sfrac{10^{-3}}{\mathrm{K}}}
		& \multicolumn{3}{w{c}{5cm}}{\colorbox{gray!20}{\makebox[4.8cm]{Viscosity}} } \\[1mm]

		\macro{\varphi_K}{-1.338357}{-}
		& \macro{\omega_{\kappa_{31}}}{1.071197}{\sfrac{10^{-3}}{\mathrm{K}}}
		& \macro{\mathrm{m}}{1.0}{-} \\[1mm] \cline{1-3}

		\macro{G_0}{80.89153}{\mathrm{GPa}}
		& \macro{\omega_{\kappa_{41}}}{3.689398}{\sfrac{10^{-4}}{\mathrm{K}}}
		& \macro{\mu}{0.001-1.0}{s} \\[1mm]\cline{4-6}

		\macro{\omega_G}{1.265276}{-}
		& \macro{\kappa_{0_2}}{22.25796}{\sfrac{\mathrm{W}}{\mathrm{m}\cdot \mathrm{K}}}
		& \multicolumn{3}{w{c}{5cm}}{\colorbox{gray!20}{\makebox[4.8cm]{Physical}} } \\[1mm]

		\macro{\chi_G}{2.136860}{\sfrac{10^{-3}}{\mathrm{K}}}
		& \macro{\omega_{\kappa_{2}}}{4.963377}{10^{-1}}
		& \macro{\rho}{7821.8}{\sfrac{\mathrm{kg}}{\mathrm{m}^3}} \\[1mm]

		\macro{\varphi_G}{-1.071127}{-}
		& \macro{\chi_{\kappa_{2}}}{2.782682}{\sfrac{10^{-3}}{K}}
		& \macro{}{}{} \\

	\bottomrule
	\end{tabular}
	\caption{Visco-plasticity model parameters according to \cite{Oppermann:2022}.}
	\label{vpParameters}
\end{table}

\vspace{3mm}

\noindent\textbf{Appendix E: Geometry of the porous plate}\\[-2mm]

\noindent The porous plate considered in Section \ref{sec:porousPlate} has dimensions \(10\,\mathrm{mm}\times 10 \,\mathrm{mm}\times 0.75 \,\mathrm{mm}\) and contains $96$ cylindrical through-holes. The geometry of each hole $i$ is defined by its in-plane center coordinates
$\mathbf{c}_i=(x_i, y_i)$ and radius $r_i, i=1,\hdots, 96$, see Table \ref{tab:porous_plate_geometry}. Thus, the physical domain is defined as
\begin{equationarray*}{rclclclc}
	\Omega &=& \bigl\{\mathbf{x} \in [0,10] \times [0,10] \times [0,0.75] \,:\,
	\left\|(x,y)-\mathbf{c}_i \right\| \geq r_i, \; i=1,\hdots, 96 \bigr\}.\\[-6mm]
\end{equationarray*}
\begin{longtable}{rccc|rccc|rccc}
\caption{Center coordinates and radii (in mm) of the cylindrical through-holes defining the porous plate.}
\label{tab:porous_plate_geometry}\\
\hline
$i$ & $x_i$ & $y_i$ & $r_i$ & $i$ & $x_i$ & $y_i$ & $r_i$ & $i$ & $x_i$ & $y_i$ & $r_i$\\
\hline
1  & 4.9749 & 3.9639 & 0.2270 & 33 & 4.6988 & 5.8980 & 0.1752 & 65 & 7.6312 & 4.0323 & 0.1977\\
2  & 5.2211 & 2.0733 & 0.2297 & 34 & 6.7948 & 5.7365 & 0.2519 & 66 & 7.1195 & 3.1088 & 0.1341\\
3  & 5.3312 & 6.4993 & 0.2427 & 35 & 8.2714 & 4.5149 & 0.1868 & 67 & 5.2588 & 3.4422 & 0.2291\\
4  & 9.3880 & 4.2739 & 0.1589 & 36 & 9.4187 & 3.0116 & 0.0982 & 68 & 7.4691 & 4.2323 & 0.1560\\
5  & 4.2877 & 5.7112 & 0.1702 & 37 & 2.2435 & 6.7872 & 0.1650 & 69 & 4.0847 & 4.7374 & 0.1290\\
6  & 3.7085 & 3.7499 & 0.1935 & 38 & 8.7106 & 4.6198 & 0.1597 & 70 & 4.5004 & 3.3071 & 0.1315\\
7  & 2.9638 & 4.1244 & 0.1224 & 39 & 5.5287 & 5.5094 & 0.1511 & 71 & 4.3656 & 4.6829 & 0.1437\\
8  & 4.6858 & 5.0811 & 0.1404 & 40 & 3.4410 & 5.6204 & 0.1236 & 72 & 1.1962 & 4.9312 & 0.1472\\
9  & 1.2589 & 6.4348 & 0.2302 & 41 & 3.7700 & 4.9130 & 0.1538 & 73 & 9.9669 & 6.1559 & 0.0729\\
10 & 2.8557 & 6.6737 & 0.2538 & 42 & 8.5503 & 4.7232 & 0.1423 & 74 & 8.2281 & 2.9265 & 0.2072\\
11 & 1.8794 & 5.1246 & 0.2533 & 43 & 1.6589 & 4.3555 & 0.2147 & 75 & 6.4253 & 5.1196 & 0.2053\\
12 & 6.9608 & 3.7824 & 0.2464 & 44 & 8.0523 & 7.4731 & 0.1456 & 76 & 5.2689 & 6.4508 & 0.1078\\
13 & 4.6817 & 2.1969 & 0.1687 & 45 & 4.8315 & 4.7830 & 0.0953 & 77 & 8.9560 & 5.6298 & 0.1480\\
14 & 2.1916 & 4.1643 & 0.2711 & 46 & 8.8931 & 3.6604 & 0.3221 & 78 & 6.5234 & 5.5631 & 0.0857\\
15 & 1.4468 & 5.1245 & 0.0816 & 47 & 2.5806 & 6.7272 & 0.2011 & 79 & 5.3378 & 6.0551 & 0.1403\\
16 & 3.6119 & 3.8439 & 0.2000 & 48 & 3.3776 & 7.9440 & 0.1283 & 80 & 6.0030 & 3.0845 & 0.1511\\
17 & 0.3589 & 3.1601 & 0.1067 & 49 & 2.1289 & 4.6781 & 0.1131 & 81 & 2.6675 & 6.4812 & 0.1222\\
18 & 0.8857 & 5.3639 & 0.2261 & 50 & 2.4192 & 3.1977 & 0.1087 & 82 & 9.1164 & 3.2062 & 0.1385\\
19 & 3.0249 & 4.6822 & 0.1695 & 51 & 3.9202 & 6.2580 & 0.1407 & 83 & 0.0519 & 4.2626 & 0.1954\\
20 & 4.5877 & 6.9234 & 0.0862 & 52 & 8.5228 & 4.0593 & 0.1189 & 84 & 3.0264 & 4.6273 & 0.1487\\
21 & 3.6265 & 4.8109 & 0.1125 & 53 & 9.1853 & 5.7906 & 0.1962 & 85 & 3.6426 & 4.9105 & 0.2112\\
22 & 4.2237 & 5.2000 & 0.1832 & 54 & 5.5751 & 4.2843 & 0.1724 & 86 & 6.7404 & 3.3576 & 0.1224\\
23 & 6.2974 & 6.4652 & 0.1824 & 55 & 3.6301 & 5.7641 & 0.1385 & 87 & 5.5882 & 3.9045 & 0.2270\\
24 & 4.1565 & 5.2867 & 0.1029 & 56 & 3.4525 & 5.7718 & 0.1831 & 88 & 5.9758 & 7.6916 & 0.2349\\
25 & 5.6008 & 4.1322 & 0.1972 & 57 & 4.5604 & 5.6902 & 0.2232 & 89 & 4.3464 & 4.2306 & 0.0792\\
26 & 0.9883 & 5.0705 & 0.1485 & 58 & 3.7326 & 2.7253 & 0.1862 & 90 & 3.8442 & 5.1425 & 0.1535\\
27 & 0.7380 & 5.9654 & 0.2140 & 59 & 7.3545 & 5.1649 & 0.1870 & 91 & 7.8553 & 6.5651 & 0.1897\\
28 & 6.8300 & 6.0894 & 0.1254 & 60 & 4.9720 & 5.8032 & 0.0732 & 92 & 6.2640 & 5.0263 & 0.1379\\
29 & 5.4991 & 4.2273 & 0.1029 & 61 & 1.3148 & 5.7985 & 0.0449 & 93 & 3.2354 & 4.4172 & 0.1909\\
30 & 1.6226 & 5.1073 & 0.0849 & 62 & 6.7798 & 7.0816 & 0.1155 & 94 & 2.7031 & 4.4177 & 0.1190\\
31 & 6.9721 & 7.0344 & 0.0768 & 63 & 4.1918 & 3.2298 & 0.1169 & 95 & 3.4872 & 3.0282 & 0.1175\\
32 & 6.1050 & 5.0637 & 0.1967 & 64 & 1.2190 & 2.5669 & 0.1254 & 96 & 6.1402 & 4.1366 & 0.2341\\

\hline
\end{longtable}

\newpage

\bibliographystyle{unsrt}
\bibliography{bibliographyR3}

\end{document}